\documentclass[a4paper,12pt]{article}
\usepackage{amsmath,amssymb}
\usepackage[dvips]{graphicx}
\makeatletter
 
  \@addtoreset{equation}{section}
 \makeatother
\newtheorem{theorem}{Theorem}
\newtheorem{proposition}[theorem]{Proposition}

\newtheorem{remark}[theorem]{Remark}

\begin{document}
\title{The Painlev\'e transcendents with solvable monodromy}
\author{Kazuo Kaneko \\ 
\small  Graduate School of Information Science and Technology,\\
\small Osaka University\\}
\date{}
\maketitle

\noindent{\bfseries Abstract}: We will study special solutions of the fourth, fifth 
and sixth Painlev\'e equations with generic values of parameters whose linear monodromy 
can be calculated explicitly. 
We will show the relation between Umemura's classical solutions and our solutions.\\ 

\section{Introduction}

\quad The Painlev\'e equation can be represented by an isomonodromic deformation of a
linear equation
\begin{eqnarray}
\frac{\partial \Psi}{\partial x}&=A(x,y(t),t)\Psi, \label{0:l}\\
\frac{\partial \Psi}{\partial t}&=B(x,y(t),t)\Psi, \label{0:d}
\end{eqnarray}
where $A(x,y,t)$ and $B(x,y,t)$ are $2\times2$ matrices.
The integrability condition of \eqref{0:l} and \eqref{0:d}  gives the Painlev\'e equation for $y(t)$. 
We call the linear equation \eqref{0:l} the {\it linearization} of the Painlev\'e equation.
 We call the monodromy data of the linear equation \eqref{0:l}   {\it  a linear monodromy}
of the Painlev\'e function $y(t)$. The linear monodromy cannot be calculated except for 
special cases. One exceptional case is Umemura's classical solutions. Umemura showed 
that there exist two kinds of special solutions for the Painlev\'e equations,
algebraic solutions and the Riccati solutions \cite{HU1}, which are called classical solutions 
of the Painlev\'e equations.
For most of all Umemura's classical solutions, the linear monodromy can be calculated, but there 
exist some Painlev\'e functions which are not included in Umemura's classical solutions,~such that 
the linear monodromy can be calculated. 
If we can determine the linear monodromy of a Painlev\'e function exactly,
we call the  Painlev\'e function  {\it monodromy\ solvable}.  

 It was R.~Fuchs who found a monodromy solvable solution at first, which is not included
in Umemura's classical solutions \cite{RF2}. He calculated the linear monodromy of
 so-called Picard's solutions,~which satisfies the sixth Painlev\'e equation with a special parameter.
 This result was rediscovered recently \cite{AVKK}, \cite{MMA}. 

\par\bigskip
The first, second and fourth Painlev\'e equations have the following 
simple symmetries which do not change the parameters in the Painlev\'e equations.
\begin{eqnarray*}
&P_{I}\qquad &y \to \zeta^3 y, \quad t \to \zeta t,  \quad (\zeta^5=1)\\
&P_{II}\qquad &y \to \omega y, \quad t \to \omega^2 t, \quad  (\omega^3=1)\\
&P_{IV}\qquad &y \to - y, \quad t \to - t, 
\end{eqnarray*}
There exist finite number of solutions which are invariant under the simple symmetries above. 
We call such a solution  as a {\it symmetric solution}.
A.~V.~Kitaev showed that the symmetric solutions for the first and the second 
Painlev\'e equation are monodromy solvable \cite{AVK}.  We remark that 
Kitaev's symmetric solution for the second Painlev\'e equation exists for 
any parameter of the equation.

We will show that the symmetric solutions  for the fourth Painlev\'e equation 
is monodromy solvable \cite{KK1} in section \ref{P4}. 
Umemura's special solutions exist only for 
special values of parameters but our new special solution exists for
any value of parameters and the associated linear equation can be reduced to the Whittaker
equation for the special initial condition at $t=0$. 
We describe the relations between the symmetric solution and Umemura's
classical solutions. The symmetric solution includes the rational solution
$y=-2t/3$ for parameters $(\alpha,\ \beta)=(0,-2/9)$. The symmetric solution also includes 
one of the Riccati solutions.  

\par\bigskip

We will study special solutions  which are meromorphic  at the origin   of  
the fifth and the sixth Painlev\'e equations in sections \ref{P5} and \ref{P6}.
In general the Painlev\'e transcendent has an essential singularity at the fixed 
singular points. 
But it is known that there exist three meromorphic solutions at $t=0$ for the 
fifth Painlev\'e equation (see \S 37 in \cite{GLS}). 
There exist four meromorphic solutions at $t=0$ for the sixth Painlev\'e equation. 
For the sixth  Painlev\'e equation, we also have four meromorphic solutions 
at $t=1$ and $t=\infty$, respectively \cite{KK2}. 
We remark that any meromorphic solution  $y(t)$ of the fifth and 
the sixth Painlev\'e equations  at $t=0$ becomes holomorphic. 

We can represent the linear monodromy by using asymptotic expansions 
of generic Painlev\'e functions.  For the fifth Painlev\'e  equation, 
such correspondence was given by Andreev and Kitaev  \cite{AK97}  \cite{AK00} 
using WKB analysis.  Although the connection formula by Andreev and Kitaev are 
very complicated, we determine the monodromy data for special solutions 
which are analytic at the origin by an elementary method.  
In section \ref{P5} we will show that the Stokes multiplier of the linear 
monodromy  of one of such solutions is zero since linearization is 
reduced to the Gauss hypergeometric equation at $t=0$.

Umemura's special solutions exist only for special values of parameters but our  solutions exist for
generic value of parameters. One of our special solutions includes the algebraic solution $y\equiv-1$ 
for special parameters $(\alpha+\beta=0, \gamma=0)$
and  also includes one point of the Riccati solution. 
We will transform  Miwa-Jimbo's linearization to a simple equation  without 
the deformation parameter for a  rational solution of the fifth Painlev\'e equation. 
The idea to calculate such a transformation of the independent variable is due to K.~Okamoto. 

For algebraic solutions, we can take a suitable transformation $z=z(x,t)$  and a 
gauge transformation $\tilde{\Psi}=R(z,t)\Psi$     such that \eqref{0:l} is transformed to 
\begin{equation}\label{0:r.fuchs}
\frac{\partial \tilde{\Psi}}{\partial z} =\tilde{A}(z)\tilde{\Psi},
\end{equation}
which does not contain the deformation parameter $t$. This fact is observed by 
R.~Fuchs \cite{RF2} at first. He gave such transformation for some algebraic solutions of 
the sixth Painlev\'e equation.   We can take a similar transformation 
for most of all algebraic solutions of Painlev\'e equations 
\cite{OO}. For the sixth  Painlev\'e equation, we do not know all of algebraic solutions, 
but many algebraic solutions are constructed by such transformations by Kitaev \cite{AVKK}.
In \cite{OO}, Ohyama and Okumura constructed such transformations for the first to the 
fifth Painlev\'e equations.  

 The meromorphic solution around the origin of the fifth Painlev\'e  equation
 appeared in \cite{JMMS}. 
They study a special fifth Painlev\'e equation with the parameter 
$\alpha=1/2, \beta=-1/2, \gamma=-2i, \delta=2$ (See the equation (2.20) in \cite{JMMS}).  
This solution is not algebraic. We think that our solutions  may have application 
to mathematical physics although they are special. 

\par\bigskip

For the sixth Painlev\'e equation Jimbo \cite{MJ2} gave  a correspondence between 
 the linear monodromy  and a local expansion of the 
Painlev\'e function. 
But if the Painlev\'e function is meromorphic around a fixed singularity,
 we can calculate the linear monodromy easily. 
We will give  twelve sets of 
meromorphic solutions around a fixed singularity for the sixth Painlev\'e equation. 
For these meromorphic solutions, we can consider confluence of singularities 
of \eqref{0:l}.  If a sixth Painlev\'e function $y(t)$ is holomorphic at $t=0$, 
we can take a limit $t\to 0$ in \eqref{0:l}. The equation \eqref{0:l} 
is still Fuchsian after we take the limit and is reduced to the Gauss hypergeometric equation.  
We can take another limit $x=1$ to $x=\infty$. In this case  \eqref{0:l} 
is reduced to the Heun equation, which is solved by elementary functions. 
Therefore we can determine the monodromy of both reduced equations and we can also
determine the monodromy of \eqref{0:l}.

A.~D.~Bryuno and I.~V.~Goryuchkina construct
asymptotic solutions around the fixed singularity \cite{BR} and D.~Guzzetti presents the leading 
term of the critical behavior at the fixed singularity \cite{DG} \cite{DG1} \cite{DG2}.
Another type of local behavior for the sixth Painlev\'e equation are studied by 
K.~Takano and S.~Shimomura    \cite{KT}, \cite{SS1}, \cite{SS2}.

R.~Fuchs showed that Picard's solution is monodromy solvable \cite{RF2}.
Picard's solution is expressed in terms of Weierstrass' $\wp$ function:
\begin{eqnarray}
y(t)&=&4\wp(c_1\omega_1(t)+c_2\omega_2(t)\, |\, g_2, g_3)+\frac{t+1}{3},\label{Pi:1}\\
&&g_2=\frac{1}{12}(t^2-t+1),\quad g_3=\frac{1}{432}(t+1)(2t-1)(t-2),\nonumber
\end{eqnarray}
where $c_1$ and $c_2$ are arbitrary constants and $\omega_1$ and $\omega_2$
are a pair of fundamental period 
\cite{EP}. This satisfies the sixth Painlev\'e equation
with special parameters:
$\alpha=\beta=\gamma=0$ and $\delta=1/2$. 
For the rational numbers $c_1$ and $c_2$, R.~Fuchs determined the monodromy invariant:
\begin{eqnarray}\label{Picard:monodromy}
p_{0t}=-2\cos2c_2\pi,\quad p_{1t}=-2\cos2c_1\pi,\quad p_{01}=-2\cos2(c_1-c_2)\pi.
\end{eqnarray}
Here $p_{ij}$=tr$M_iM_j$\ and $M_j$ is a monodromy matrix (See section \ref{P6} and \cite{MJ2}).
He showed that
Picard's solution \eqref{Pi:1} is expanded at $t=0$ as
\begin{eqnarray}
y(t)=-4\frac{e^{2\pi ik/n}}{2^{8l/n}}t^{\frac{2l}{n}}+a_1t^{\frac{2l+1}n}+\cdots
 \end{eqnarray}
for $c_1=  \frac{k}{n}, \,c_2=\frac{l}{n}$ and $\frac{l}{n}<\frac{1}{2}$. In the case 
$\frac{l}{n}>\frac{1}{2}$, we have similar expansion of $y(t)$.
Therefore we can take a limit $t\to 0$ in \eqref{0:l} and \eqref{0:l} is
also reduced to the Gauss hypergeometric equation. Since \eqref{Pi:1} has 
a similar expansion at $t=1$, we can take another limit $t\to 1$ in \eqref{0:l}. 
Because the limit of \eqref{0:l} is reduced to the Gauss hypergeometric equation again,
we obtain the monodromy invariant \eqref{Picard:monodromy}. 

The  paper \cite{RF2} was completely forgotten for 
long years. The author thanks Professor Y.~Ohyama who introduced him the paper \cite{RF2}. 

\par\bigskip

We will list up all of known monodromy solvable solutions: 
\begin{itemize}
 \item Umemura's classical solutions  \cite{AVK1}, \cite{OO} 
 \item Symmetric solutions (\cite{AVK}, \cite{KOK}, section \ref{P4})
 \item meromorphic solutions around fixed singularities (\cite{KOH}, \cite{KK2}, sections \ref{P5} and \ref{P6})
 \item Picard's solution \cite{RF2}
\end{itemize}
We do not have a rigorous definition of monodromy solvability. One may  think any Painlev\'e
 function is monodromy solvable. For example, 
Jimbo \cite{MJ2} gave a correspondence between local expansion of generic solutions 
of the sixth Painlev\'e equation  at $t=0$  and linear monodromy. In this sense,
a generic sixth Painlev\'e function is monodromy solvable. Similarly, there exist 
a correspondence between local expansion of Painlev\'e functions and linear monodromy
for other types of Painlev\'e equations \cite{AK00}. But our monodromy solvable solutions 
listed above are more special, since we can determine the linear monodromy by reducing 
classical special functions, such as the Gauss hypergeometric function or the 
Kummer confluent hypergeometric function. 

\par\bigskip

In section \ref{painleve}, we review the Painlev\'e equations. In subsection 
\ref{isomonodromic} we list the "Lax form" of the Painlev\'e equations.  
In subsection \ref{backlund}, we review the B\"acklund transformation groups of the Painlev\'e 
equations.

In section \ref{P4}, we show that the symmetric solution of the fourth 
Painlev\'e equation is monodromy solvable. This section is based on the paper 
\cite{KK1}. In section \ref{P5}, we show that meromorphic solutions at $t=0$
 of the fifth Painlev\'e equation are monodromy solvable. This section is based on the paper 
\cite{KOH}. In section \ref{P6}, we show that meromorphic solutions at $t=0$ of the sixth 
Painlev\'e equation are monodromy solvable. This section is based on the paper 
\cite{KK2}. 

\par\bigskip
The author wishes to thank Professor Y.~Ohyama 
for his constant guidance and encouragements over time to complete this work.  
The author also thanks  Dr.~D.~Guzzetti, Professor K.~Okamoto and 
Professor S.~Shimomura for fruitful discussions.


\section{The Painlev\'e equations}\label{painleve}
The Painlev\'e equations was found by Paul Painlev\'e about one hundred years
ago \cite{P:06}. He and his pupil Gambier classified second order nonlinear equations 
without movable singularities \cite{Gambier}.  After they removed equations which can be solved by 
known functions, the following six equations are remained.

\begin{eqnarray}
P_{I})&  y^{\prime\prime}=&6y^2+t,\\
P_{II})&  y^{\prime\prime} =&2y^3+ty+\alpha,\\
P_{III})&  y^{\prime\prime} =& \frac{1}{ y}{y^{\prime}}^2-\frac{y^{\prime}}{ t}+
\frac{\alpha y^2 + \beta}{t}+\gamma y^3 + \frac{\delta }{ y},\\
P_{IV})&  y^{\prime\prime}=& \frac{1}{2y}{y^{\prime}}^2+\frac{3}{ 2}y^3+ 4 t y^2 
+2(t^2-\alpha)y +\frac{\beta}{y},\\
P_{V})&  y^{\prime\prime}=& \left(\frac{1}{2y}+\frac{1}{y-1} \right){y^{\prime}}^2-\frac{1}{t}
{y^{\prime}}+\frac{(y-1)^2}{t^2}\left(\alpha y +\frac{\beta}{y}\right) +\gamma {y\over t}+
\delta {y(y+1)\over y-1},\label{p5}\\
P_{VI})&  y^{\prime\prime} =& {1\over 2}\left({1\over y}+{1\over y-1}+{1 \over y-t}\right){y^{\prime}}^2
-\left({1\over t}+{1\over t-1}+{1 \over y-t}\right)y^{\prime}\nonumber\\
& &+ {y(y-1)(y-t)\over t^2(t-1)^2}\biggl[\alpha+\beta {t\over y^2}+\gamma {t-1\over (y-1)^2}
+ \delta {t(t-1)\over (y-t)^2}\biggr].\label{p6}
\end{eqnarray}
Here $\alpha, \beta, \gamma$ and $\delta$ are complex parameters.
They are called the Painlev\'e equations.  It is known that generic solutions 
of these six equations are transcendental functions  and they are called
the Painlev\'e transcendents.  

We may use a different type of the third Painlev\'e equation $P_{\rm III}^\prime$
\begin{equation*}\label{pd}
q''=\frac{{q'}^2}{q}  - \frac{q'}{x} 
+\frac{\alpha {q}^2}{4x^2}   +\frac{\beta}{4\,x} +\frac{\gamma{q}^3}{4x^2} +
 \frac{\delta}{4q}.
\end{equation*}
instead of $P_{\rm III}$, since it is easy to study isomonodromic deformation 
for $P_{\rm III}^\prime$.  
$P_{\rm III}^\prime$ is equivalent to $P_{\rm III}$ by 
\begin{equation*}
x=t^2, \quad y=t q.
\end{equation*}

The third Painlev\'e equation is divided into three types:
\begin{itemize}
	\item $D_8^{(1)}$ if $\alpha\not=0,\ \beta\not=0$,\ $\gamma=0,\ \delta=0$,
	\item $D_7^{(1)}$ if $\delta=0,\ \beta\not=0$ or $\gamma=0,\  \alpha\not=0$,
	\item $D_6^{(1)}$ if $\gamma\delta\not=0$.
\end{itemize}
In the case $\beta=0,\ \delta=0$ (or $\alpha=0,\ \gamma=0$), the third Painlev\'e equation 
is a quadrature, and we exclude this case from the Painlev\'e family.  $D_j^{(1)} (j=6,7,8)$ mean 
the affine Dynkin diagrams corresponding to Okamoto's initial value spaces. 
 By   suitable scale transformations
$t \to ct, y \to dt$, we may fix $\gamma=4, \delta=-4$ for $D_6^{(1)}$, 
and $\gamma=2$ for $D_7^{(1)}$. 

For the fifth Painlev\'e equation, we assume that $\delta\not=0$. 
When $\delta=0, \gamma\not=0$, the fifth equation is equivalent to the third equation of 
the $D_6^{(1)}$ type.  When $\delta=0, \gamma=0$, the fifth equation is
quadrature and we exclude this case from the Painlev\'e family.
By a suitable scale transformation $t \to ct$, we can fix $\delta=-1/2$ for the fifth equation.

\subsection{Isomonodromic deformation equations}\label{isomonodromic}
In 1905, R.~Fuchs showed that the sixth Painlev\'e equation is an isomonodromic 
deformation equation of a second order Fuchsian linear differential equation \cite{RF05} \cite{RF1}. 
Later Garnier showed that other Painlev\'e equations are  also isomonodromic 
deformation equations of a second order linear differential equation with 
irregular singularities \cite{Garnier:1902}.  

We will list up the isomonodromic  deformation equations for 
all Painlev\'e equations. We use Miwa-Jimbo's form \cite{MJ}, which is isomonodromic  deformation
of $2\times 2$ matrix type linear equations
\begin{equation}\label{isomono}\begin{aligned}
\frac{\partial Y}{\partial x}&=A(x,t)Y\\
\frac{\partial Y}{\partial t}&=B(x,t)Y.
\end{aligned}\end{equation}
For a  suitable pair $A$ and $B$, the integrability condition 
\begin{equation}\label{Lax}
\frac{\partial A}{\partial t}(x,t)-\frac{\partial B}{\partial x}(x,t)+[A(x,t),B(x,t)]=0
\end{equation}
gives the Painlev\'e equations.


\subsubsection{The first Painlev\'e equation}
We take
\begin{equation}\begin{aligned}
A(x,t)&=\begin{pmatrix}0&1\\ 0&0 \end{pmatrix}x^2 +
\begin{pmatrix}0&y\\ 4&0\end{pmatrix}x+\begin{pmatrix}-z&y^2+t/2 \\-4y&z \end{pmatrix},\\
B(x,t)&=\begin{pmatrix} 0 & 1/2 \\ 0 & 0 \end{pmatrix}x +
\begin{pmatrix} 0 & y \\ 2 & 0 \end{pmatrix}.
\end{aligned}\end{equation}
The  integrability condition \eqref{Lax} is
$$
\frac{dy}{dt}=z, \quad \frac{dz}{dt}=6y^2+t,
$$ 
and we obtain the first Painlev\'e equation
$$\frac{d^2y}{dt^2}=6y^2+t.$$

\noindent By a transformation
$$x=\zeta^2, \quad Y(x)=\begin{pmatrix}1&0\\ 0&\zeta^{-1}\end{pmatrix}
\begin{pmatrix}1&1\\ 2&-2\end{pmatrix}Z(\zeta),$$
\eqref{isomono} is transformed into 
\begin{equation}\label{PI_Lax}\begin{aligned}
\frac{dZ}{d\zeta}&=A_0(\zeta,t)Z,\\
\frac{dZ}{dt}&=B_0(\zeta,t)Z,
\end{aligned}\end{equation}
where
$$A_0=
\begin{pmatrix}4&0\\ 0&-4\end{pmatrix}\zeta^4+
\begin{pmatrix}0&-4y\\4y&0\end{pmatrix}\zeta^2+
\begin{pmatrix}0&-2z\\2z&0\end{pmatrix}\zeta+
\begin{pmatrix}1&-1\\ 1&-1 \end{pmatrix}(2y^2+t)+
\begin{pmatrix}1&-1\\ -1&1\end{pmatrix}\frac{1}{2\zeta},
$$
$$B_0=
\begin{pmatrix}1&0\\ 0&-1\end{pmatrix}\zeta+
\begin{pmatrix}y&-y\\ y&-y\end{pmatrix}\frac{1}{ \zeta}.
$$
The first equation in \eqref{PI_Lax} has a regular singular point at $\zeta=0$, where 
the local exponents are 0 and 1. Since solutions have no logarithmic terms at $\zeta=0$,
 $\zeta=0$ is an apparent singularity.
A formal solution is given by
$$Z(\zeta)=\left(1+\frac{Z_1}{\zeta}+\frac{Z_2}{\zeta^2}+\cdots \right)e^T(\zeta),$$
$$T(\zeta)=\frac45\begin{pmatrix}1&0\\ 0&-1\end{pmatrix} \zeta^5 +
\begin{pmatrix} t&0\\ 0&-t\end{pmatrix} \zeta +\frac12\log \zeta,\quad 
Z_1=\begin{pmatrix}-H_I&0\\ 0&-H_I\end{pmatrix},$$
where 
$$H_I=\frac12 z^2-(2y^3+ty) $$
is a Hamiltonian of the first Painlev\'e equation.

\subsubsection{The second Painlev\'e equation }
We take
\begin{equation}\begin{aligned}
A(x,t)&=\begin{pmatrix}1&0\\ 0&-1 \end{pmatrix}x^2 +
\begin{pmatrix}0 & u \\ -2u^{-1}z & 0\end{pmatrix}x+
\begin{pmatrix} z + t/2 & -u y \\ -2u^{-1}(\theta + y z) & -z -t/2 \end{pmatrix},\\
B(x,t)&=\frac12 \begin{pmatrix} 1 & 0 \\ 0 & -1 \end{pmatrix}x +
\frac12\begin{pmatrix} 0 & u \\ -2u^{-1} z & 0 \end{pmatrix}.
\end{aligned}\end{equation}
Here $\alpha = \frac12 -\theta$.

\noindent The integrability condition \eqref{Lax} is
\begin{eqnarray*}
\frac{dy}{dt}&=&y^2 +z+\frac t2, \quad \frac{dz}{dt}= - 2 yz -\theta,\\
\frac{du}{dt}&=&-uy,
\end{eqnarray*}
and we obtain the second Painlev\'e equation
$$\frac{d^2y}{dt^2}=2y^3+ty+\left(\frac12 -\theta\right).$$

\noindent At $x=\infty$, a formal solution is given by
$$Y(x)=\left(1+\frac{Y_1}{x}+\frac{Y_1}{x^2}+\cdots \right)e^T(x),$$
$$T(x)=\begin{pmatrix}1&0\\ 0&-1\end{pmatrix} \frac{x^3}3 +
\begin{pmatrix} t&0\\ 0&-t\end{pmatrix}\frac{x}{2}+
\begin{pmatrix} \theta&0\\ 0&-\theta\end{pmatrix}\log(\frac1x),$$
$$Y_1=\begin{pmatrix}-H_{II}&-u/2\\ -z/u& H_{II}\end{pmatrix},\quad
Y_2=\begin{pmatrix}H_{II}^2/2 +(z-t\theta)/4& uy/2-uH_{II}/2\\ 
-(\theta+yz)/u+zH_{II}/u&H_{II}^2/2 +(z+t\theta)/4\end{pmatrix},
$$
where 
$$H_{II}=\frac12 z^2+(y^2+\frac t2)z+\theta y $$
is a Hamiltonian of the second Painlev\'e equation.

\subsubsection{The third  Painlev\'e equation of type $D_8^{(1)}$}
We take 
\begin{equation}\label{d8:imd}\begin{aligned}
A(x,t)&= \begin{pmatrix}
          z & x-y \\
            \rule[8mm]{0mm}{0mm}
          - \dfrac{t^2}{4 y x^3} - \dfrac{y z^2 + z - (1/4)}{x^2}
          - \dfrac{z^2}{x} &
           \rule[-4mm]{0mm}{0mm}
          - \dfrac{2}{x} - z
          \end{pmatrix},\\
 B(x,t)&= \begin{pmatrix}
           \dfrac{2y z}t &  \dfrac{2y x}t \\
          \dfrac{s}{2 y x^2}-\dfrac{2y z^2  +t z' }{ t x}
         &\  -\dfrac{2yz}t
          \end{pmatrix}.
\end{aligned}\end{equation}
\noindent Then the integrability condition \eqref{Lax} is 
\begin{equation}\label{d8:hamil}
ty'=4 y^2z+2y, \quad tz'= - 4 yz^2  -2z -\frac{t^2}{2y^2} +\frac12.
\end{equation}
This is a Hamiltonian system with the Hamiltonian
$$
tH=2y^2z^2+2yz -\dfrac y2 -\dfrac{t^2}{2y}. $$
From \eqref{d8:hamil} we obtain
$$q^{\prime\prime}=\frac{(q')^2}q-\frac{q'}s
+\frac{2q^2}{ s^2}-2.
$$
Changing the variable $s=t^2$, we have $P_{III}'(\alpha=2, \beta=-2, \gamma=0, \delta=0)$
 $$q^{\prime\prime}=\frac{(q')^2}q-\frac{q'}s
+\frac{q^2}{ 2s^2}-\frac1{2s}.
$$

\noindent By a transformation
$$x=\zeta^2, \quad 
Y=\begin{pmatrix} 2\zeta^3& 2\zeta^3\\-2 z\zeta + 1 & -2 z\zeta - 1\end{pmatrix}Z,$$
the Lax form is changed into
\begin{equation*}\begin{aligned}
\frac{dZ}{d\zeta }&=\left( A_0+\frac1\zeta  A_1 +\frac1{\zeta ^2} A_2\right)Z,\quad
\frac{\partial Z}{\partial t}&=-\frac{A_2}{t\zeta } Z,
\end{aligned}\end{equation*}
where
$$A_0=\begin{pmatrix}
     1 & 0 \\
     0& -1 \end{pmatrix},\, 
A_1=\frac12 \begin{pmatrix}
     -7 & 4yz+1 \\
     4yz+1 & -7 \end{pmatrix},\, 
A_2=\frac1{2y}\begin{pmatrix}
     -y^2-t^2 & y^2-t^2 \\
     -y^2+t^2 & y^2+t^2 \end{pmatrix}.$$
 
\subsubsection{The third  Painlev\'e equation of type $D_7^{(1)}$}
We take
\begin{equation}\begin{aligned}
A(x,t)&= 
 \begin{pmatrix}0 &0 \\ -yz^2+\theta_0 z+1/4 & -t \end{pmatrix}\frac1{x^2}+ \begin{pmatrix} 0& 0 \\ -z^2 & t \end{pmatrix}\frac1{x}+\begin{pmatrix}z &-y \\0 & -z \end{pmatrix}
+\begin{pmatrix}0 &1 \\0 & 0 \end{pmatrix}x,\\
B(x,t)&=\begin{pmatrix}  \frac{yz}t & \frac{y x}t \\
     -\frac{z^2q+tz'}{tx} & \frac {t-yz x}{tx}
    \end{pmatrix},
\end{aligned}\end{equation}
\noindent Then the integrability condition \eqref{Lax} is 
\begin{equation}\label{d7:hamil}
ty'=  2 y^2 z-\theta_0 y  +t, \quad tz'= -2yz^2  +\theta_0 z +\frac14.
\end{equation}
This is a Hamiltonian system with the Hamiltonian
$$
tH=y^2z^2 +(-\theta_0 y+t)z -\frac y4. $$
From \eqref{d7:hamil} we obtain $P_{III}'(\alpha=2, \beta=4(\theta_0+1), \gamma=2, \delta=0)$
$$q^{\prime\prime}=\frac{(q')^2}q-\frac{q'}t
+\frac{q^2}{2t^2}-\frac1q+\frac{1+\theta_0}t.
$$

\noindent By a transformation
$$x=\zeta^2, \quad Y(x)=\begin{pmatrix}\zeta^2&0\\ 0&1\end{pmatrix}
\begin{pmatrix}1&0\\ -z&1\end{pmatrix}\begin{pmatrix}\zeta&0\\ 0&1\end{pmatrix}
\begin{pmatrix}2&2\\ 1&-1\end{pmatrix}Z(\zeta),$$
\eqref{isomono} is transformed into 
\begin{equation}\label{PD7_Lax}\begin{aligned}
\frac{dZ}{d\zeta}&=A_0(\zeta,t)Z,\\
\frac{dZ}{dt}&=B_0(\zeta,t)Z,
\end{aligned}\end{equation}
where
\begin{equation*}\begin{aligned}
A_0=&
   \begin{pmatrix}
                  1 &  0 \\
                  0 & -1
                 \end{pmatrix} + 
   \frac{1}{2}
   \begin{pmatrix}
    2 \theta_0 - 5 & 4 y z - 2 \theta_0 - 1 \\
    4 y z - 2 \theta_0 - 1 & 2 \theta_0 - 5
   \end{pmatrix}\frac1{\zeta} \\
&+ 
   \frac{1}{2}
   \begin{pmatrix}
      4 z t - y &   4 z t + y \\
    - 4 z t - y & - 4 z t + y
   \end{pmatrix}\frac1{\zeta^2} + 
\begin{pmatrix}
           -t &  t \\
            t & -t
           \end{pmatrix}\frac1{\zeta^3},\\
B_0=&
   \begin{pmatrix}
     - z   + y/4t &  - z - y/4t \\
       z   + y/4t &    z - y/4t
   \end{pmatrix}\frac1{\zeta } +
\begin{pmatrix}
            1/2 & -1/2 \\
           -1/2 &  1/2
           \end{pmatrix}\frac1{\zeta^2}.
\end{aligned}\end{equation*}

\subsubsection{The third Painlev\'e equation of type $D_6^{(1)}$ }
We take
\begin{equation}\begin{aligned}
A(x,t)&= \frac12 \begin{pmatrix}t&0\\ 0&-t \end{pmatrix} +
\frac1x \begin{pmatrix}-\theta_\infty/2 & u \\ v & \theta_\infty/2 \end{pmatrix}+
\frac1{2x^2} G \begin{pmatrix} -t & 0 \\ 0 & t \end{pmatrix}G^{-1},\\
B(x,t)&=\frac12 \begin{pmatrix} 1 & 0 \\ 0 & -1 \end{pmatrix}x +
\frac1t\begin{pmatrix} 0 & u \\ v & 0 \end{pmatrix} +
\frac1{2x}G \begin{pmatrix} 1& 0 \\ 0 & -1 \end{pmatrix} G^{-1},
\end{aligned}\end{equation}
where $G= \begin{pmatrix} a & b  \\ c & d \end{pmatrix}\in SL(2, \mathbb{C})$.  We set
\begin{equation}\begin{aligned}
 \frac12 G \begin{pmatrix} -t & 0 \\ 0 & t \end{pmatrix}G^{-1} &= 
 \begin{pmatrix} z-t/2 & -wz \\ w^{-1}(z-t) & -z+t/2 \end{pmatrix}, \\
G^{-1} \begin{pmatrix}-\theta_\infty/2 & u \\ v & \theta_\infty/2 \end{pmatrix}&G=
\begin{pmatrix} \theta_0/2 & \bar{u} \\ \bar{v} & -\theta_0/2 \end{pmatrix}.
\end{aligned}\end{equation}
These parameters satisfy the following constraints:
$$ ad-1=bc =-z/t. \qquad ab=-wz/t, \qquad cd=-(z-t)/tw,$$
$$ (\theta_0+\theta_\infty)a-2c u+2b  \bar{v}=0, \quad
 (\theta_0 - \theta_\infty) b - 2 a \bar{u}+2du=0, \quad
 (\theta_0+\theta_\infty)c -2av+2d\bar{v}=0.$$

\noindent The integrability condition \eqref{Lax} is
\begin{equation*}\begin{aligned}
t\frac{dG}{dt}&=\begin{pmatrix} 0 & u \\ v & 0 \end{pmatrix}G 
+G \begin{pmatrix} 0 & \bar{u} \\ \bar{v} & 0 \end{pmatrix}, \\
\frac{du}{dt}&= \frac{\theta_\infty}t u+ 2tab,\quad
\frac{dv}{dt} = -\frac{\theta_\infty}t v+ 2tcd,
\\
\frac{d\bar{u}}{dt}&= \frac{\theta_0}t \bar{u}+ 2tbd,\quad
\frac{d\bar{v}}{dt} = -\frac{\theta_0}t \bar{v}+ 2tac.
\end{aligned}\end{equation*}
We set $y=-u/zw$.  Then we have
\begin{equation*}\begin{aligned}
t\frac{dy}{dt}&=4zy^2-2ty^2+(2\theta_\infty-1)y+2t, \\
t\frac{dz}{dt}&= -4yz^2+(4ty-2\theta_\infty+1)z+(\theta_0+\theta_\infty)t,\\
t\frac{d }{dt}&\log w = -\frac{(\theta_0+\theta_\infty)t}z-2ty+\theta_\infty.
\end{aligned}\end{equation*}
We obtain the third Painlev\'e
$$
y^{\prime\prime}= {1\over y}{y^{\prime}}^2-{y^{\prime}\over t}+
{\alpha y^2 +\beta\over t}+\gamma y^3 + {\delta \over y},
$$
$$\alpha=4\theta_0, \ \beta =4(1-\theta_\infty), \
\gamma=4, \delta=-4.$$
At $x=\infty$, a formal solution  is
$$Y(x)=\left(1+\frac{Y_1}{x}+\frac{Y_2}{x^2}+\cdots \right)e^T(x),$$
$$T(x)=\frac12\begin{pmatrix}t&0\\ 0&-t\end{pmatrix} x
       + \frac12 \begin{pmatrix}\theta_\infty&0\\ 0&-\theta_\infty\end{pmatrix}\log\left( \frac1x\right), $$
$$Y_1=\begin{pmatrix}-uv/t-z+t/2 & -u/t\\  v/t&  uv/t+z-t/2\end{pmatrix}.$$
At $x=0$, a formal solution  is
$$\bar{Y}(x)=\left(1+\frac{\bar{Y}_1}{x}+\frac{\bar{Y}_2}{x^2}+\cdots \right)e^{\bar{T}(x)},$$
$$\bar{T}(x)=\frac12\begin{pmatrix}t&0\\ 0&-t\end{pmatrix} \frac1x
       + \frac12 \begin{pmatrix}\theta_0&0\\ 0&-\theta_0 \end{pmatrix}\log x, $$
$$\bar{Y}_1=\begin{pmatrix}-\bar{u}\bar{v}/t-z+t/2 & -\bar{u}/t\\  \bar{v}/t& \bar{u}\bar{v}/t+z-t/2\end{pmatrix}.$$

\noindent The Hamiltonian of the third Painlev\'e equation is
\begin{equation*}\begin{aligned}
t H_{III}&
=2y^2z^2 +2(-ty^2+\theta_\infty y+t)z-(\theta_0+\theta_\infty)ty
-t^2-\frac{\theta_0^2-\theta_\infty^2}4,
\end{aligned}\end{equation*}
and
$$2 H_{III}={\rm Tr}\left( Y_1+\bar{Y}_1\right)\begin{pmatrix}-1&0\\ 0&1 \end{pmatrix}.$$


 \subsubsection{The fourth Painlev\'e equation}
We take 
\begin{eqnarray}
A(x,t)&=&
  \left(
    \begin{array}{cc}
    1 &0 \\
    0 &-1
    \end{array}
\right)
x
+
\left(
    \begin{array}{cc}
    t &u \\
    \frac{2}{u}(z-\theta_0-\theta_\infty) &-t
    \end{array}
\right)
+\frac{1}{x}
\left(
    \begin{array}{cc}
    -z+\theta_0 &-\frac{uy}{2} \\
    \frac{2z}{uy}(z-2\theta_0) &z-\theta_0
    \end{array}
\right),    \nonumber \\
B(x,t)&=&
\left(
    \begin{array}{cc}
    1 &0 \\
    0 &-1 
    \end{array}
\right)
x
+
\left(
    \begin{array}{cc}
    0 &u \\
    \frac{2}{u}(z-\theta_0-\theta_\infty) &0
    \end{array}
\right), \nonumber
\end{eqnarray}
where $y,z$ and $u$ are functions of $t$, and $\theta_0$ and $\theta_\infty$ are constants
\begin{eqnarray}
\alpha=2\theta_\infty-1,\qquad  \beta=-8 \theta_0^2.
\end{eqnarray}
Setting $w=z/y$, the integrability condition \eqref{Lax} gives
\begin{equation}\label{p4:ham}\begin{aligned}
 \frac{dy}{dt}&= -4yw+y^2+2ty+4\theta_0,  \\
 \frac{dw}{dt}&= 2w^2-2yw-2tw+(\theta_0+\theta_\infty),\\
 \frac{d\log u}{dt}&= -y-2t.
\end{aligned}\end{equation}
The system \eqref{p4:ham} is the Hamiltonian system with the polynomial Hamiltonian $H_4$:  
\begin{eqnarray}      
H_4=-2yw^2+y^2w+2tyw+4\theta_0w-(\theta_0+\theta_\infty)y.
\end{eqnarray}
The function $u$ can be obtained from \eqref{p4:ham} by a quadrature.

\subsubsection{The fifth Painlev\'e equation}
We take
\begin{eqnarray}
&A(x,t)&=
  \frac{1}{2}\left(
    \begin{array}{cc}
    t &0 \\
    0 &-t
    \end{array}
\right)
+
\frac{1}{x}\biggl(
    \begin{array}{cc}
    z+\frac{\theta_0}{2} &-u(z+\theta_0) \\
    u^{-1}z &-z-\frac{\theta_0}{2}
    \end{array}
\biggr)\nonumber \\
& &+\frac{1}{x-1}
\left(
    \begin{array}{cc}
    -z-\frac{\theta_0+\theta_\infty}{2} & uy\left( z+\frac{\theta_0-\theta_1+\theta_\infty}{2} \right) \\
    -\frac{1}{uy}\left(z+\frac{\theta_0+\theta_1+\theta_\infty}{2}\right)
 & z+\frac{\theta_0+\theta_\infty}{2}
    \end{array}
    \right), \label{5:16}    \\
 &B(x,t)&=
 \frac{1}{2}\left(
    \begin{array}{cc}
    1 &0 \\
    0 &-1 
    \end{array}
\right)
x \nonumber \\
 &&+\frac{1}{t}
 \left(
    \begin{array}{cc}
    0 &-u\left[z+\theta_0-y\left(z+\frac{\theta_0-\theta_1+\theta_\infty}{2}\right)\right] \\
    \frac{1}{u}\left[z-\frac{1}{y}\left(z+\frac{\theta_0+\theta_1+\theta_\infty}{2}\right)\right]
    &0
    \end{array}
\right),\label{5:17}
\end{eqnarray}
where $y,z$ and $u$ are functions of $t$, and $\theta_0,\theta_1$ and $\theta_\infty$ are parameters. 
From the integrability condition \eqref{Lax}, we have
\begin{eqnarray*}
t\frac{dy}{dt}&=&ty-2z(y-1)^2-(y-1)\left(\frac{\theta_0-\theta_1+\theta_\infty}{2}y
-\frac{3\theta_0+\theta_1+\theta_\infty}{2}\right),         \\
t\frac{dz}{dt}&=&yz\left(z+\frac{\theta_0-\theta_1+\theta_\infty}{2}\right)-\frac{z+\theta_0}{y}
\left(z+\frac{\theta_0+\theta_1+\theta_\infty}{2}\right),      \\
t\frac{d\log u}{dt}&=&-2z-\theta_0+y\left(z+\frac{\theta_0-\theta_1+\theta_\infty}{2}\right)
+\frac{1}{y}\left(z+\frac{\theta_0+\theta_1+\theta_\infty}{2}\right).                    
\end{eqnarray*}
Eliminating $z$, we have the fifth Painlev\'e equation for
\begin{eqnarray*}
\alpha=\frac{1}{2}\left(\frac{\theta_0-\theta_1+\theta_\infty}{2}\right)^2,  
\beta=\frac{-1}{2}\left(\frac{\theta_0-\theta_1-\theta_\infty}{2}\right)^2,
\gamma=1-\theta_0-\theta_1,\delta=\frac{-1}{2}.
\end{eqnarray*}
Putting
\begin{eqnarray*}
w=\frac{1}{y}\left(z+\frac{\theta_0+\theta_1+\theta_\infty}{2}\right),  
\end{eqnarray*}
we have      
\begin{eqnarray}
 t\frac{dy}{dt}&=&ty-2y^3w+y^2\left(\frac{\theta_0+3\theta_1+\theta_\infty}{2}\right) \nonumber\\
 & & +4y^2w-y(2\theta_1+\theta_\infty)-2yw-\frac{\theta_0-\theta_1-\theta_\infty}{2}, \label{2:y} \\
 t\frac{dw}{dt}&=&3y^2w^2-yw(\theta_0+3\theta_1+\theta_\infty)+w^2-wt-4yw^2 \nonumber \\
 & & +w(2\theta_1+\theta_\infty)+\frac{\theta_1(\theta_0+\theta_1+\theta_\infty)}{2}, \label{2:z}\\
t\frac{d\log u}{dt}&=&-2yw+\theta_1+\theta_\infty+y\left(yw-\theta_1\right)+w.   \label{2:u}
\end{eqnarray}
 The system \eqref{2:y}  and \eqref{2:z} are the Hamiltonian system with the polynomial Hamiltonian $H_5$ as shown below: \\
\begin{eqnarray*}
tH_5&=&-y(y-1)^2 w^2+ \left[ \left(\frac{-\theta_0+ \theta_1+\theta_\infty}{2}\right)(y-1)^2
+ ( \theta_0+\theta_1)y(y-1) +t y +(\theta_0+\theta_1)  \right] w \nonumber\\
&&-\frac{\theta_1(\theta_0+\theta_1+\theta_\infty)}{2}y.
\end{eqnarray*}
The function $u$ can be obtained from \eqref{2:u} by a quadrature.

\subsubsection{The sixth Painlev\'e equation}\label{isomonodromic:p6}

We take
\begin{equation*}
A(x,t)=\sum_{j=0,1,t}
  \frac{A_j}{x-j}
=\left(
    \begin{array}{cc}
    a_{11}(x,t) &a_{12}(x,t) \\
    a_{21}(x,t) &a_{22}(x,t)
    \end{array}
\right), \quad
B(x,t)=-\frac{A_t}{x-t}, \label{J:1}
\end{equation*}
where
\begin{equation}
A_j=\left(
    \begin{array}{cc}
    z_j+\theta_j &-u_jz_j \\
    u_j^{-1}(z_j+\theta_j) &-z_j
    \end{array}
\right)\quad(j=0,1,t).  \nonumber
\end{equation}
We define $A_{\infty},\ y$ and $z$ as follows:
\begin{eqnarray*}
A_\infty&=&-\sum_{j=0,1,t}A_j
=\left(\begin{array}{cc} \frac{1}{2}(\theta_{\infty}-\sum_{j=0,1,t}\theta_j) 
& 0 \\ 0 & -\frac{1}{2}(\theta_{\infty}+\sum_{j=0,1,t}\theta_j) 
\end{array}\right),\\
a_{12}(x,t)&=&-\sum_{j=0,1,t}\frac{u_jz_j}{x-j}
=\frac{k(x-y)}{x(x-1)(x-t)},\\
z&=&-a_{11}(y,t)=\sum_{j=0,1,t}\frac{z_j+\theta_j}{y-j},
\end{eqnarray*}
where $y, z, z_j, u_j$ and $k$ are functions of $t$ and $\theta_j\,(j=0,1,t,\infty)$ 
are parameters. \\
We then have
\begin{eqnarray*}
 \sum_{j=0,1,t}z_j=-\frac{1}{2}\left(\sum_{i=0, 1, t, \infty} \theta_i\right),\quad
\sum_{j=0,1,t}u_jz_j=0,\\
\sum_{J=0,1,t}u_j^{-1}(z_j+\theta_j)=0,\quad 
(t+1)u_0z_0+tu_1z_1+ u_tz_t=k.
\end{eqnarray*}
In what follows, instead of $\theta_j$, we mainly use the parameters $\alpha_j\, (j=0,1,2,3,4)$ defined by
 the following relations:
\begin{equation}
\theta_0=\alpha_4,\quad \theta_1=\alpha_3,\quad \theta_t=\alpha_0,
\quad \theta_{\infty}=1-\alpha_1\quad
(\alpha_0+\alpha_1+2\alpha_2+\alpha_3+\alpha_4=1).\label{A:L}
\end{equation}
From the integrability condition  \eqref{Lax}, we have   
\begin{eqnarray}  
t(t-1)\frac{dy}{dt}&=& 2zy(y-1)(y-t)-\alpha_4(y-1)(y-t)-\alpha_3y(y-t) \nonumber\\ 
&&\hskip 1cm-(\alpha_0-1)y(y-1),   \label{S:1}     \\
t(t-1)\frac{dz}{dt}&=& \left(-3y^2+2(1+t)y-t\right)z^2+\biggl[(2y-1-t)\alpha_4
+(2y-t)\alpha_3  \nonumber \\
& &\hskip 1cm +(2y-1)(\alpha_0-1)\biggr]z-\alpha_2(\alpha_1+\alpha_2). \label{S:2}
 \end{eqnarray} 
Eliminating $z$, we have the sixth Painlev\'e equation with
\begin{eqnarray} 
\alpha=\frac{\alpha_1^2}{2}=\frac{(1-\theta_\infty)^2}{2},\ 
\beta=\frac{-\alpha_4^2}{2}=\frac{-1}{2}\theta_0^2,
\ \gamma=\frac{\alpha_3^2}{2}=\frac{1}{2}\theta_1^2,\
\delta=\frac{1-\alpha_0^2}{2}=\frac{1-\theta_t^2}{2} . \label{p:r} 
\end{eqnarray}
The system   \eqref{S:1} and \eqref{S:2}  can be written as a Hamiltonian system 
with the polynomial Hamiltonian $H_{VI}$ given by 
\begin{eqnarray}      
t(t-1)H_{VI}&=&y(y-1)(y-t)z^2-\biggl[\alpha_4(y-1)(y-t)+\alpha_3y(y-t)
+(\alpha_0-1)y(y-1)\biggr]z \nonumber\\
&& +\alpha_2(\alpha_1+\alpha_2)(y-t).
\end{eqnarray}

\begin{remark}
This polynomial Hamiltonian system is the same 
as Garnier-Okamoto's Hamiltonian system. Putting $Y={}^t(\psi_1,\psi_2)$ and
eliminating $\psi_2$ from \eqref{isomono}, we have the same second order single equation 
as Garnier-Okamoto's equation \cite{KO}.
\end{remark}

\subsubsection{ Normalized   form of the sixth Painlev\'e equation}
In this section, we give the normalized Jimbo-Miwa's isomonodromic deformation equations 
 whose linear monodromy belongs to $SL(2,\mathbb {C})$. We use this system 
for the calculation of the monodromy data and the asymptotic expansion of $\tau$-function.
 Put 
\begin{eqnarray}
\bar{Y}=x^{-\frac{\theta_0}{2}}(x-1)^{-\frac{\theta_1}{2}}(x-t)^{-\frac{\theta_t}{2}}
Y,
\end{eqnarray}
 in \eqref{isomono}.
Then we have
\begin{eqnarray}
\frac{\partial \bar{Y}(x,t)}{\partial x}&=&\bar A(x,t)\bar{Y}(x,t),  \
\bar A(x,t)=
  \sum_{j=0,1,t}\frac{\bar A_j}{x-j}=\left(
    \begin{array}{cc}
    \bar a_{11}(x,t) &\bar a_{12}(x,t) \\
    \bar a_{21}(x,t) &\bar a_{22}(x,t)
    \end{array}
\right),\label{J:3}
\end{eqnarray}
\begin{eqnarray}
\bar A_j&=&\left(
    \begin{array}{cc}
    \bar z_j+\frac{\theta_j}{2} &-\bar u_j\bar z_j \\
    \bar u_j^{-1}(\bar z_j+\theta_j) &-\bar z_j-\frac{\theta_j}{2}
    \end{array}
\right)\quad (j=0,1,t),\nonumber\\
 \frac{\partial \bar{Y}(x,t)}{\partial t}&=&\bar B(x,t)\bar{Y}(x,t),  \quad
\bar B(x,t)=-\frac{\bar 
A_t}{x-t}.\label{J:4}
\end{eqnarray}
We define $\bar A_{\infty}, y$ and $\bar z$ as follows:
\begin{eqnarray}
\bar A_\infty&=&-\sum_{j=0,1,t}\bar A_j=\left(\begin{array}{cc} \frac{\theta_{\infty}}{2} & 0 \\ 0 & 
-\frac{\theta_{\infty}}{2}   
\end{array}\right),\quad
\bar a_{12}(x,t)=
\frac{\bar k(x-y)}{x(x-1)(x-t)},\\
\bar z&=&\bar a_{11}(y,t)=\sum_{j=0,1,t}\frac{\bar z_j+\frac{\theta_j}{2}}{y-j},
\end{eqnarray}
where $y, \bar z, \bar z_j, \bar u_j, \bar k$ 
are functions of $t$ and $\theta_j, \theta_\infty$ are parameters.\\ 
We then have
\begin{eqnarray}
 \sum_{j=0,1,t}\bar z_j=-\frac{1}{2}\left(\sum_{i=0,1,t,\infty} \theta_i\right),\quad
\sum_{j=0,1,t}\bar u_j\bar z_j=0,\\
\sum_{J=0,1,t}\bar u_j^{-1}(\bar z_j+\theta_j)=0,\quad 
(t+1)\bar u_0\bar z_0+t\bar u_1\bar z_1+\bar u_t \bar z_t=\bar k.
\end{eqnarray}
We can solve as
follows: \\
\begin{eqnarray}
\bar u_0&=&\frac{\bar ky}{t\bar z_0},\quad\bar u_1=-\frac{\bar k(y-1)}{(t-1)\bar z_1},
\quad\bar u_t=\frac{\bar k(y-t)}{t(t-1)\bar z_t},\\
\bar z_0&=&\frac{1}{t\theta_\infty}\biggl[y^2(y-1)(y-t)\bar z^2
+\theta_{\infty}y(y-1)(y-t)\bar z+\frac{\theta_{\infty}^2}{4}(y-1)(y-t)\nonumber\\
&&-\frac{1}{4}(\theta_{\infty}+\theta_0)^2t
+\frac{\theta_1^2}{4}\cdot\frac{y}{y-1}(t-1)-\frac{\theta_t^2}{4}t(t-1)\frac{y}{y-t} 
\biggr],\\
\bar z_1&=&\frac{-1}{(t-1)\theta_\infty}\biggl[y(y-1)^2(y-t)\bar z^2
+\theta_\infty y(y-1)(y-t)\bar z\nonumber\\
&&+\frac{\theta_{\infty}^2}{4}y(y-t)+\frac{1}{4}(\theta_{\infty}+\theta_1)^2(t-1)
-\frac{\theta_0^2}{4}\cdot\frac{y-1}{y}t-\frac{\theta_t^2}{4}t(t-1)\frac{y-1}{y-t} 
\biggr],
\end{eqnarray}
\begin{eqnarray}
\bar z_t&=&\frac{1}{t(t-1)\theta_\infty}\biggl[y(y-1)(y-t)^2\bar z^2
+\theta_\infty y(y-1)(y-t)\bar z\nonumber\\
&&+\frac{\theta_{\infty}^2}{4}y(y-1)-\frac{1}{4}(\theta_{\infty}+\theta_t)^2t(t-1)
-\frac{\theta_0^2}{4}\cdot\frac{y-t}{y}t+\frac{\theta_1^2}{4}(t-1)\frac{y-t}{y-1}
\biggr].
\end{eqnarray} Hereinafter we 
use $\alpha_j$ which are defined by \eqref{A:L}.
From the integrability condition of \eqref{J:3} and \eqref{J:4}, we have
\begin{eqnarray}
t(t-1)\frac{dy}{dt}&=&2y(y-1)(y-t)\bar z+y(y-1), \label{y:0}  \\
t(t-1)\frac{d \bar z}{dt}&=&\left[-3y^2+2(1+t)y-t\right]\bar z^2-(2y-1)\bar z\nonumber \\
&&+\biggl[-\frac{1-\alpha_1^2}{4}-\frac{\alpha_4^2}{4}\cdot\frac{t}{y^2}
+\frac{\alpha_3^2}{4}\cdot\frac{t-1}{(y-1)^2}-\frac{\alpha_0^2}{4}\cdot\frac{t(t-1)}{(y-t)^2} 
\biggr].\label{z:0}
\end{eqnarray}
Eliminating $\bar z$, we again obtain the sixth Painlev\'e equation with \eqref{p:r}.
The system of equations \eqref{y:0} and \eqref{z:0} is a rational Hamiltonian system 
with the Hamiltonian $\bar H_{VI}$ defined by 
\begin{eqnarray}      
t(t-1)\bar H_{VI}&=&y(y-1)(y-t)\bar z^2 +y(y-1)\bar z \nonumber\\
&&-\biggl[\frac{1-\alpha_1^2}{-4}y+\frac{\alpha_4^2}{4}\cdot\frac{t}{y}
-\frac{\alpha_3^2}{4}\cdot\frac{t-1}{(y-1)}+\frac{\alpha_0^2}{4}\cdot\frac{t(t-1)}{(y-t)} 
\biggr].
\end{eqnarray}
From \eqref{S:1} and \eqref{y:0}, we have
\begin{eqnarray}
2(z-\bar z)=\frac{\alpha_4}{y}+\frac{\alpha_3}{y-1}+\frac{\alpha_0}{y-t}.\label{ct:1}
\end{eqnarray}
\begin{remark} The transformation \eqref{ct:1} gives the following canonical transformation 
between two Hamiltonian systems \eqref{S:1}, \eqref{S:2} and \eqref{y:0}, \eqref{z:0}
which keeps $y$ invariant: 
\begin{eqnarray}
   dz\wedge dy -dH_{VI}\wedge dt=d\bar z\wedge dy-d\bar H_{VI}\wedge dt.
\end{eqnarray}
\end{remark}

\subsection{The B\"acklund transformation groups}\label{backlund}
There exist rational transformations which change a Painlev\'e equation 
to another Painlev\'e equation of the same type with different parameters.  
The transformation group of each type of Painlev\'e equations 
is called the B\"acklund transformation group.  The B\"acklund transformation group
is isomorphic to an affine Weyl group.

For a classical root system $R$, we denote the Weyl group by  ${W}(R)$.
 We denote by $P$ and $Q$ the weight lattice and the root lattice of $R$, respectively \cite{Hum}. 
 It is known that the affine Weyl group ${W}(R^{(1)}) \cong Q \ltimes W(R)$.
 We set $\widehat{W}(R^{(1)})= P \ltimes W(R)$. Let $G$ be the Dynkin automorphism group of the 
extended Dynkin diagram. The quotient $P/Q$ is contained in G.
We denote  the extended affine Weyl group by  $\widetilde{W}(R) \cong G \ltimes {W}(R^{(1)})$.
 \par\bigskip

Since the first Painlev\'e equation has no parameter, it does not 
have any B\"acklund transformation.  We will list up all of the 
B\"acklund transformations for the Painlev\'e equation from the second to
the sixth.  

\subsubsection{Simple symmetry}
For the first, second and fourth Painlev\'e equations, there exist
simple transformations which keep the parameters.  
\begin{eqnarray*}
&P_{I}\qquad &y \to \zeta^3 y, \quad t \to \zeta t,  \quad (\zeta^5=1)\\
&P_{II}\qquad &y \to \omega y, \quad t \to \omega^2 t, \quad  (\omega^3=1)\\
&P_{IV}\qquad &y \to - y, \quad t \to - t, 
\end{eqnarray*}
They are not contained in the B\"acklund transformation groups. 
We will use these symmetry to define symmetric solutions of 
the Painlev\'e equations.

\subsubsection{The second Painlev\'e equation}
The Hamiltonian is 
\begin{equation}
H_{II}={1 \over 2}p^2- \left(q^2+{t \over 2} \right) p-\alpha_1 q.
\end{equation}
The equation for $y=q$ is the second Painlev\'e equation:
\begin{equation}
\frac{ d^2 y}{ dt^2} = 2 y^3+t y+\alpha,
\end{equation}
where $\alpha=\alpha_1-{1 \over 2}$.

\par\bigskip

\noindent The B\"acklund transformation  is 
\begin{eqnarray*}
&&\widetilde{W}(A_1^{(1)})=G\ltimes W(A_1^{(1)})= \langle  s_1, \pi\rangle, \\
&& W(A_1^{(1)})=\langle s_0, s_1\rangle, \\
&&G=P/Q={\rm Aut}(E^{(1)}_7)= {\rm Aut}(A^{(1)}_1)
= \langle  \pi\rangle \cong \mathbb{Z}_2.
\end{eqnarray*}
The birational transformations are given by:
\begin{equation*}
\begin{array}{|c||c|c|c|c|c|}
\hline
 & \alpha_0 & \alpha_1 & q & p & t \\
\hline
s_0  & -\alpha_0 & \alpha_1 +2 \alpha_0
& q+{\alpha_0 \over f}
& p+{4 \alpha_0 q \over f}+{2 {\alpha_0}^2 \over f^2} 
& t \\
s_1  & \alpha_0+2 \alpha_1 & -\alpha_1 
& q+{\alpha_1 \over p} & p & t \\
\hline
\pi & \alpha_1 & \alpha_0 & -q & -f  & t \\
\hline
\end{array},
\end{equation*}
where $\alpha_0=1-\alpha_1$ and $f=p-2 q^2-t$.

\subsubsection{The third Painlev\'e equation of $D_{8}^{(1)}$ type}
The Hamiltonian is:
\begin{equation}
tH_{D_8} =q^2 p^2+qp-\frac{1}{2}\left( q+\frac{t}{q}\right).
\end{equation}
The equation for $y=q/\tau$,  $t=\tau^2$  is the special case of
the third Painlev\'e equation:
\begin{equation}
{ d^2 y \over  d\tau^2}=
{1 \over y} \left({ dy \over  d\tau} \right)^2
-{1 \over \tau}{ dy \over  d\tau}
+{4 \over \tau}(  y^2-1)+4 y^3-{4\over y}.
\end{equation} 

\noindent The symmetry of the equation is:
\begin{eqnarray*}
&G=\langle \pi\rangle \cong \mathbb{Z}_2.
\end{eqnarray*}
The birational transformations are given by:
\par\medskip

\begin{center}
  \begin{tabular}[h]{|c||c|c|c|}\hline
  & $q$ & $p$ & $t$ \\ \hline
$\pi$  &  ${\renewcommand\arraystretch{1.2}\begin{array}{c} \displaystyle \frac{t}{q} \end{array}}$     &
 $\displaystyle -\frac{q(2q p + 1 )}{2t}$ & $t$ \\ 
\hline
  \end{tabular}.\end{center}

\subsubsection{The third Painlev\'e equation of $D_{7}^{(1)}$ type}
The Hamiltonian is 
\begin{equation}
tH_{D_7} =q^2 p^2+\alpha_1 qp+tp+ q.
\end{equation}
The equation for $y=q/\tau$, $t=\tau^2$  is the special case of
the third Painlev\'e equation:
\begin{equation}
{ d^2 y \over  d\tau^2}=
{1 \over y} \left( { dy \over  d\tau} \right)^2
-{1 \over \tau}{ dy \over  d\tau}
+{1 \over \tau}(-8 y^2+\beta)- {4 \over y},
\end{equation}
with
\begin{equation}
\beta=4(1-\alpha_1).
\end{equation}

\noindent The symmetry of the equation is:
\begin{eqnarray*}
& \widetilde{W}(A_1^{(1)})=\langle s_1, \sigma\rangle, \\
& G=\langle \pi\rangle \cong {\mathbb Z},
\end{eqnarray*}
where $\pi=\sigma \circ s_1$.
The birational transformations are given by:
\begin{equation}
\begin{array}{|c||c|c|c|c|c|}
\hline
 & \alpha_0 & \alpha_1 & q & p & t  \\
\hline
s_0  & -\alpha_0 & \alpha_1 +2 \alpha_0
& q & p+{\alpha_0 \over q}-{t \over q^2}  & -t  \\
s_1  & \alpha_0+2 \alpha_1 & -\alpha_1 & -q + \frac{\alpha_1}{p} + \frac{1}{p^2}  & -p & -t  \\ \hline
\sigma  & \alpha_1 & \alpha_0 & t p & - \frac{q}{t}  & -t  \\
\hline
\end{array},
\end{equation}
where $\alpha_0=1-\alpha_1$.

\noindent Any element of $G$ has no fixed value of parameters.

\subsubsection{The third Painlev\'e equation of $D_{6}^{(1)}$ type}
The Hamiltonian is 
\begin{equation}
t H_{D_6} =q^2 p^2-(q^2-(\alpha_1+\beta_1) q-t)p - \alpha_1 q.
\end{equation}
The equation for $y= q/\tau$, $t=\tau^2$ is 
the third Painlev\'e equation:
\begin{equation}
{ d^2 y \over  d\tau^2}
=
{1 \over y} \left( { dy \over  d\tau} \right)^2
-{1 \over \tau}{  dy \over  d\tau}
+{1 \over \tau}(\alpha y^2+\beta)+4 y^3-{4 \over y},
\end{equation}
with
\begin{equation}
\alpha=4 (\alpha_1-\beta_1), \quad
\beta=-4(\alpha_1+\beta_1-1).
\end{equation}

\noindent The symmetry of the equation is:
\begin{eqnarray*}
&& \widetilde{W}((2A_1)^{(1)})=G\ltimes W((2A_1)^{(1)})= \langle s_0, s_1, s_0^{\prime}, s_1^{\prime}, \pi_1, \pi_2, \sigma\rangle, \\
&&\widehat{W}((2A_1)^{(1)}) = \langle s_0, s_1, s_0^{\prime}, s_1^{\prime}, \pi_1, \pi_2 \rangle, \\
&& W((2A_1)^{(1)})=\langle s_0, s_1, s_0^{\prime}, s_1^{\prime}\rangle, \\
&&G={\rm Aut}(D^{(1)}_6)={\rm Aut}((2A_1)^{(1)})
=\langle \pi_1, \pi_2, \sigma\rangle \cong  {\mathfrak D}_8, \\
&&P/Q =\langle \pi_1, \pi_2\rangle \cong \mathbb{Z}_2\times \mathbb{Z}_2.
\end{eqnarray*}
Since the table of  birational transformations are too long, we split into two parts:
\begin{flushleft}$
\begin{array}{|c||c|c|c|c|c }
\hline
 & \alpha_0 & \alpha_1 & \beta_0  & \beta_1 & \hskip 5mm {} \\
\hline
s_0  & -\alpha_0 & \alpha_1+2\alpha_0 & \beta_0 & \beta_1 
&   \\
s_1  & \alpha_0+2\alpha_1 & -\alpha_1 & \beta_0 & \beta_1 
&  \\
s_0^{\prime}  & \alpha_0 & \alpha_1 & -\beta_0 & \beta_1+2\beta_0  
&    \\
s_1^{\prime}  & \alpha_0 & \alpha_1 & \beta_0+2\beta_1 & -\beta_1 
&  \\
\hline
\pi_1  & \alpha_1 & \alpha_0 & \beta_0 & \beta_1 
&  \\
\pi_2  & \alpha_0 & \alpha_1 & \beta_1 & \beta_0 
& \\
\hline
\sigma_1  & \beta_0 & \beta_1 & \alpha_0 & \alpha_1 
&   \\
\sigma_2   & \beta_1 & \beta_0 & \alpha_1 & \alpha_0 
& \\
\hline
\end{array} $
\end{flushleft}
\begin{flushright}
$ 
\begin{array}{|c|| c|c|c|c|}
\hline
 & \hskip 5mm {}  & q & p  & t \\
\hline
s_0  &  
&  q + \frac{ \alpha_0 }{ (p-1)+ \frac{\alpha_1 + \beta_1 - 1}{q} + \frac{t}{q^2}  }
&      p - \frac{ \alpha_0 (  2 q(p-1)+ \alpha_1 + \beta_1 - 1) }{f_1} - \frac{ \alpha_0^2 \, t}{ f_1^2}  
&   t   \\
s_1  &  
& q + \frac{\alpha_1}{p}   
& p 
&   t \\
s_0^{\prime}  &   
& q + \frac{ \beta_0 }{ p + \frac{\alpha_1 + \beta_1 - 1}{q} + \frac{t}{q^2}  }   
&    p - \frac{ \beta_0 ( 2 q p + \alpha_1 + \beta_1 - 1 ) }{f_2} - \frac{ \beta_0^2 \, t}{ f_2^2 } 
&   t   \\
s_1^{\prime}  &  
& q + \frac{\beta_1}{p-1} & p & t \\
\hline
\pi_1  &  
& - \frac{t}{q} 
& \displaystyle{\frac{q}{t}(q(p-1)+\beta_1)+1}   
&  t \\
\pi_2  &  
& \frac{t}{q}
& \displaystyle{-\frac{q}{t}(q p + \alpha_1)}  
& t \\
\hline
\sigma_1  &  
& -q & 1-p& -t  \\
\sigma_2   &  
& q & \displaystyle{ p+ \frac{\alpha_1 + \beta_1-1}{q} + \frac{t}{q^2}} 
& -t\\
\hline
\end{array},
$\end{flushright}
where $\alpha_0=1-\alpha_1$, \  $\beta_0=1-\beta_1$ ,\  $f_1 
= \beta_1 q + (p-1)q^2 + t$ \  and \  $f_2 = \alpha_1 q + p q^2 + t $.
\\


\subsubsection{The fourth Painlev\'e equation}
The Hamiltonian is 
\begin{equation}
H_{IV}=(p-q-2 t) p q-2 \alpha_1 p-2 \alpha_2 q.
\end{equation}
The equation for $y=q$ is the fourth Painlev\'e equation:
\begin{equation}
{ d^2 y \over  dt^2}
=  \frac{1}{2 y} \left(\frac{ d y }{  dt } \right)^2
+{3 \over 2}y^3 +4 t y^2+2 (t^2 -\alpha) y+{\beta \over y},
\end{equation}
where
\begin{equation}
\alpha=2\theta_{\infty}-1=\alpha_0-\alpha_2, \quad
\beta=-8\theta_0^2=-2 {\alpha_1}^2.
\end{equation}
 
\noindent The symmetry of the equation is:
\begin{eqnarray*}
&& \widetilde{W}(A_2^{(1)})
=G\ltimes W(A_2^{(1)})= \langle s_0, s_1, s_2,   \sigma_1, \sigma_2\rangle, \\
&& \widehat{W}(A_2^{(1)})=\langle s_0, s_1, s_2, \pi \rangle,\\
&&  {W}(A_2^{(1)})=\langle s_0, s_1, s_2  \rangle, \\
&&G={\rm Aut}(E^{(1)}_6)={\rm Aut}(A^{(1)}_2)
=\langle \sigma_1 , \sigma_2\rangle \cong{\mathfrak S}_3, \\
&& P/Q= \langle \pi\rangle \cong \mathbb{Z}_3.
\end{eqnarray*}
The birational transformations are: 
\begin{equation*}
\begin{array}{|c||c|c|c|c|c|c|}
\hline
  & \alpha_0 &\alpha_1 & \alpha_2 & q & p & t \\
\hline
s_0  & -\alpha_0 & \alpha_1+\alpha_0 & \alpha_2+\alpha_0
& q+{2 \alpha_0 \over f}& p+{2 \alpha_0 \over f}  & t \\
s_1  & \alpha_0+\alpha_1 & -\alpha_1 & \alpha_2+\alpha_1 
& q
& p-{2 \alpha_1 \over q} & t \\
s_2  & \alpha_0+\alpha_2 & \alpha_1+\alpha_2 & -\alpha_2 
& q+{2 \alpha_2 \over p}
& p  & t \\
\hline
\pi  & \alpha_1 & \alpha_2 & \alpha_0
& -p& -f& t \\
\hline
\sigma_1 & \alpha_0 & \alpha_2 & \alpha_1 
 & -\sqrt{-1}p & -\sqrt{-1}q  & \sqrt{-1}t \\
\sigma_2  & \alpha_2 & \alpha_1 & \alpha_0 
& \sqrt{-1}f& \sqrt{-1}p & \sqrt{-1}t \\
\hline
\end{array},
\end{equation*}
where $\alpha_0=1-\alpha_1-\alpha_2$ and $f=p-q-2 t$.

\subsubsection{The fifth Painlev\'e equation}
The Hamiltonian is
\begin{eqnarray*}
t\tilde  H_{V}  = p(p+t)q (q-1)+\alpha_2 q t-\alpha_3 p q- \alpha_1 p(q-1).
\end{eqnarray*}
The Hamiltonian system 
\begin{equation*}\label{2:ha5}
\tilde{\cal H}_{V}: \
\left\{\begin{array}{l}
 t  q^{\prime } =  q (2 pq  -2p + t q  -t-\alpha_1 - \alpha_3 ) +\alpha_1, \\
 t  p^{\prime } =- p (2 p  q - p+ 2 t q-t-\alpha_1 -\alpha_3 ) - \alpha_2 t, \\
\end{array}
\right.
\end{equation*}
is equivalent to the fifth  Painlev\'e equation by $y=1-1/q$  for 
\begin{equation*}
\alpha={{\alpha_1}^2 \over 2}, \quad
\beta=-{{\alpha_3}^2 \over 2}, \quad 
\gamma=\alpha_0-\alpha_2, \quad
\delta=-{1 \over 2},
\end{equation*}
\begin{equation*}
\alpha_0+\alpha_1+\alpha_2+\alpha_3=1.
\end{equation*}
We list the relations between different parameters:
\begin{equation*}
\alpha_1    =\dfrac{\theta_0-\theta_1+\theta_\infty}{2},\quad 
\alpha_3  =\dfrac{\theta_0-\theta_1-\theta_\infty}{2},\quad 
\alpha_0-\alpha_2   =1-\theta_0-\theta_1. 
\end{equation*}
We remark that we fix $\delta=-1/2$.

The transformation 
$$w=-pq^2-\alpha_2q$$
gives a canonical transformation
\begin{eqnarray*}
 dp\wedge dq+dt\wedge d\tilde{H}_{V}=-(dw \wedge dy+dt \wedge dH_5).
\end{eqnarray*}
The B\"acklund transformation group of the fifth Painlev\'e equation is:
\begin{eqnarray*}
&&\widetilde{W}(A_3^{(1)})=G\ltimes W(A_3^{(1)})=\langle s_0, s_1, s_2, s_3, \pi, \sigma \rangle, \\
&&\widehat{W}(A_3^{(1)})=W(A_3)\ltimes P(A_3)=\langle s_0, s_1, s_2, s_3, \pi \rangle, \\
&&W(A_3^{(1)})=W(A_3)\ltimes R(A_3)=\langle s_0, s_1, s_2, s_3\rangle, \\
&&G= {\rm Aut}(A^{(1)}_3)
=\langle\sigma , \pi\rangle\cong{\mathfrak D}_8, \\
&&P /Q  \cong {\mathbb{Z}}_4.
\end{eqnarray*}

The birational transformations are 
\begin{equation*}
\begin{array}{|c||c|c|c|c|c|c|c|}
\hline
x & \alpha_0 &\alpha_1 & \alpha_2 & \alpha_3 & q & p & t \\
\hline
s_0(x) & -\alpha_0 & \alpha_1+\alpha_0 & \alpha_2 
& \alpha_3+\alpha_0 
& q+{\alpha_0 \over p+t} & p & t \\
s_1(x) & \alpha_0+\alpha_1 & -\alpha_1 
& \alpha_2+\alpha_1 & \alpha_3 
& q& p-{\alpha_1 \over q}  & t \\
s_2(x) & \alpha_0  & \alpha_1+\alpha_2 & -\alpha_2 & \alpha_3+\alpha_2 
& q+{\alpha_2 \over p} & p & t \\
s_3(x) & \alpha_0+\alpha_3 & \alpha_1 &\alpha_2+\alpha_3 
& -\alpha_3 
& q& p-{\alpha_3 \over q-1} & t \\
\hline
\pi(x) & \alpha_1 & \alpha_2 & \alpha_3 & \alpha_0
& -{p \over t} & t(q-1) & t \\
\hline
\sigma(x) & \alpha_0 & \alpha_3 & \alpha_2 & \alpha_1 
& 1-q & -p & -t \\
\hline
\end{array}.
\end{equation*}
\par\medskip\noindent


\subsubsection{The sixth Painlev\'e equation}
The Hamiltonian is
\begin{eqnarray}
t(t-1)H_{VI}&=&q(q-1)(q-t)p^2-[\alpha_4(q-1)(q-t)+\alpha_3 q(q-t)
 \nonumber \\
&&+(\alpha_0-1)q(q-1)]p +\alpha_2(\alpha_1+\alpha_2)(q-t).
\end{eqnarray}
The equation for $y=q$ is the sixth Painlev\'e equation:
\begin{eqnarray*}
\frac{  d^2 y}{  d  t^2  }&=&{1 \over 2}\left({1 \over y}+{1 \over y-1}+
{1 \over y-t}\right) \left( \frac{  d y}{  d t  }    \right)^2-\left({1 \over t}+{1
\over t-1}+{1 \over y-t}\right)  \frac{  d y}{  d t  }    \\
&&+{y(y-1)(y-t) \over t^2(t-1)^2}
\left[\alpha +\beta{t \over y^2}+
\gamma{t-1 \over (y-1)^2}+
\delta{t(t-1) \over (y-t)^2}\right], 
\end{eqnarray*}
where
\begin{equation}
\alpha= \frac{{\alpha_1}^2}{2}, \quad
\beta= -\frac{{\alpha_4}^2}{2}, \quad
\gamma= \frac{{\alpha_3}^2}{2}, \quad
\delta= - \frac{{\alpha_0}^2-1}{2},
\end{equation}
and $\alpha_0+\alpha_1+2 \alpha_2+\alpha_3+\alpha_4=1$. 

The symmetry of the equation is described as follows:
\begin{eqnarray*}
&& \widetilde{W}(D_4^{(1)})=G\ltimes W(D_4^{(1)})=\langle s_0, s_1, s_2, s_3, s_4, \sigma_1, \sigma_2, \sigma_3\rangle,\\
&& \widehat{W}(D_4^{(1)})=\langle s_0, s_1, s_2, s_3, s_4, \pi_1, \pi_2\rangle,\\
&& W(D_4^{(1)})=\langle s_0, s_1, s_2, s_3, s_4\rangle, \\
&& G={\rm Aut}(D^{(1)}_4)={\mathfrak S}_4=\langle \sigma_1, \sigma_2, \sigma_3\rangle, \\
&& P/Q=\mathbb{Z}_2\times \mathbb{Z}_2=\langle \pi_1, \pi_2\rangle.
\end{eqnarray*}

\noindent The list of birational transformations are given by the following table:
\begin{equation*}
\begin{array}{|c||c|c|c|c|c|c|c|c|}
\hline
  & \alpha_0 & \alpha_1 & \alpha_2 
& \alpha_3 & \alpha_4 & q & p & t  \\ 
\hline
s_0  & -\alpha_0 & \alpha_1 & \alpha_2+\alpha_0 & \alpha_3 
& \alpha_4 & q & p-{\alpha_0 \over q-t}  & t \\
s_1  & \alpha_0 & -\alpha_1 & \alpha_2+\alpha_1 & \alpha_3 
& \alpha_4 & q & p & t  \\
s_2  & \alpha_0+\alpha_2 & \alpha_1+\alpha_2 & -\alpha_2 
& \alpha_3+\alpha_2 & \alpha_4+\alpha_2 
& q+{\alpha_2 \over p} & p & t  \\
s_3  & \alpha_0 & \alpha_1 & \alpha_2+\alpha_3 & -\alpha_3 & \alpha_4 
& q& p -{\alpha_3 \over q-1} & t  \\
s_4  & \alpha_0 & \alpha_1 & \alpha_2+\alpha_4 & \alpha_3 & -\alpha_4
& q & p-{\alpha_4 \over q} & t  \\
\hline
\pi_1   & \alpha_3 & \alpha_4 & \alpha_2 & \alpha_0 & \alpha_1 
& \frac{t}{q}& - \frac{q(q p + \alpha_2)}{t}   & t \\
\pi_2   & \alpha_1 & \alpha_0 & \alpha_2 & \alpha_4 & \alpha_3 
&  {(q-1)t \over (q-t)} & -{p(q-t)^2+\alpha_2(q-t) \over t(t-1)} 
& t \\
\hline
\sigma_1   & \alpha_0 & \alpha_1 & \alpha_2 & \alpha_4 & \alpha_3 
& 1-q & -p    & 1-t  \\
\sigma_2   & \alpha_0 & \alpha_4 & \alpha_2 & \alpha_3 & \alpha_1 
& \frac{1}{q}& -q(q p + \alpha_2)     & \frac{1}{t}  \\
\sigma_3   & \alpha_4 & \alpha_1 & \alpha_2 & \alpha_3 & \alpha_0 
&  \frac{t-q}{t-1}& -(t-1)p  & \frac{t}{t-1}  \\
\hline
\end{array}
.\end{equation*}
Here $\pi_1= \sigma_2 \sigma_1 \sigma_3 \sigma_1$ and $\pi_2 = \sigma_1 \sigma_2 \sigma_3 \sigma_2$.

\section{Symmetric solutions of the fourth Painlev\'e equation}\label{P4}
In this section, we will determine the linear monodromy of 
symmetric solutions of the fourth Painlev\'e equation.  This section is 
based on the paper \cite{KK1}.  In order to determine the linear monodromy, 
it is sufficient to calculate the monodromy data of the linearization 
for a  special variable $t$, since the Painlev\'e equation is given by 
isomonodromic deformation. For the symmetric solutions, we fix the variable $t=0$.
Then the linearization of the fourth Painlev\'e equation is reduced to 
the Whittaker confluent hypergeometric equation.  Therefore we can determine 
the linear monodromy for the symmetric solution. 

For special parameters, the symmetric solutions become Umemura's classical 
solutions. We will compare the symmetric solutions with classical 
solutions.

\subsection{Symmetric solutions}
The first, the second and the fourth Painlev\'e equations have simple symmetry
explained in subsection 2.2.1. We call a solution of the Painlev\'e equation 
which is invariant under the action of the simple symmetry as a {\it symmetric
solution}.  We list all of symmetric solutions. The symmetric
solutions of the first and the second   Painlev\'e equations was found by
A.~V.~Kitaev \cite{AVK}. 

\begin{proposition}\label{symmetric} 1) For $P_{I}$, we have two symmetric solutions
\begin{eqnarray*}
y &=& \frac16 t^3 +\frac1{336}t^8+\frac1{26208}t^{13}+\frac{95}{224550144}t^{18}+\cdots,\\
y &=& t^{-2} -\frac1{6}t^3+\frac1{264}t^{8}-\frac1{19008}t^{13}+\cdots.
\end{eqnarray*}

2) For $P_{II}$($\alpha$), we have three symmetric solutions
\begin{eqnarray*}
y &=& \frac{\alpha}2 t^2 +\frac{\alpha}{40} t^5  +\frac{10\alpha^3+\alpha}{2280}t^{8}+\cdots,\\
y &=&  t^{-1}-\frac{\alpha+1}{4} t^2   +\frac{(\alpha+1)(3\alpha+1)}{112}t^{5}+\cdots,\\
y &=& - t^{-1} -\frac{\alpha-1}{4} t^2  -\frac{(\alpha-1)(3\alpha-1)}{112}t^{5}+\cdots.
\end{eqnarray*}
They are equivalent to each other by the B\"acklund transformations.

3) For $P_{IV}$($\alpha, -8\theta_0^2$), we have four symmetric solutions
\begin{eqnarray*}
y &=& \pm 4\theta_0\left(t -\frac{2 \alpha}{3} t^3 
     +\frac{2}{15} (\alpha^2+12\theta_0^2\pm 8\theta_0+1 )t^{5}+\cdots\right),\\
y &=& \pm t^{-1}+\frac{2}{3}(\pm\alpha-2) t 
       \mp \frac{2}{45}(-7\alpha^2\pm 16\alpha+36\theta_0^2-4 )t^{3}+\cdots.
\end{eqnarray*}
They are equivalent to each other by the B\"acklund transformations.
\end{proposition}
A.~V.~Kitaev showed that symmetric solutions are monodromy solvable for the
the first and the second   Painlev\'e equations \cite{AVK}.

\noindent From Proposition \ref{symmetric} the solutions of \eqref{p4:ham} with initial data $y(0)=0$ and 
$w(0)=0$ are expanded as follows: 

\begin{eqnarray}
 y&=&4\theta_0t\sum_{k=0}^\infty a_kt^{2k}, \label{3.1} \\
 &&a_0=1, \quad a_1=\frac{-2}{3}(2\theta_\infty-1), \nonumber\\
 &&a_2=\frac{1}{30}\{4(2\theta_\infty-1)^2+3(4\theta_0)^2+8(4\theta_0)+4\}, \cdots ,\nonumber  \nonumber \\
 w&=&(\theta_0+\theta_\infty)t\sum_{k=0}^{\infty}b_{k}t^{2k}, \label{3.2}   \\
 &&b_0=1,\quad b_1=\frac{2}{3}(\theta_\infty-3\theta_0-1), \nonumber\\
 &&b_2=\frac{4}{15}\{(\theta_\infty-3\theta_0-1)^2+4\theta_0(2\theta_\infty-1)\},\cdots .\nonumber
\end{eqnarray} 
We will determine the linear monodromy of the above solution. 

 \subsection{Transformation of the linear equation}
The linearization of the fourth Painlev\'e equation is given by 
\begin{equation}\label{linear:p4}
\frac{\partial Y}{\partial x} =A(x,t)Y,
\end{equation}
where
$$
A(x,t)=
  \left(
    \begin{array}{cc}
    1 &0 \\
    0 &-1
    \end{array}
\right)
x
+
\left(
    \begin{array}{cc}
    t &u \\
    \frac{2}{u}(yw-\theta_0-\theta_\infty) &-t
    \end{array}
\right)
+\frac{1}{x}
\left(
    \begin{array}{cc}
    -yw+\theta_0 &-\frac{uy}{2} \\
    \frac{2w}{u }(yw-2\theta_0) &yw-\theta_0
    \end{array}
\right).    
$$

\noindent By putting $t=0, y=0 $ and  $w=0$, we have 
\begin{eqnarray} \label{3.3}
\frac{d}{dx}\left( \begin{array}{c} y_1 \\ y_2  \end{array}  \right)=
\left(  \begin{array}{cc} x+\frac{\theta_0}{x} &u \\ \frac{-2(\theta_0+\theta_\infty)}{u} &-x-\frac{\theta_0}{x}
  \end{array}  \right)  \left(  \begin{array}{cc} y_1 \\ y_2  \end{array} \right). 
\end{eqnarray}                            
By the transformation $x^2=\xi$\ and\ $y_i=\xi^{\frac{-1}{4}}v_i, (i=1,2)$, we have the 
Whittaker equations:
\begin{eqnarray}
 \frac{d^2v_1}{d\xi^2}+\left[\frac{-1}{4}+\frac{k}{\xi}+\frac{\frac{1}{4}-m^2}{\xi^2} \right]v_1=0,   \\
 \frac{d^2v_2}{d\xi^2}+\left[\frac{-1}{4}+\frac{k+\frac{1}{2}}{\xi}+\frac{\frac{1}{4}-(m+\frac{1}{2})^2}{\xi^2}\right]v_2=0,
  \\   k=\frac{2\theta_\infty-1}{4}, \qquad m=\frac{2\theta_0-1}{4} .
 \end{eqnarray}
Therefore we have
\begin{theorem}
 The symmetric solution \eqref{3.1} and \eqref{3.2} of the fourth Painlev\'e equation 
is monodromy solvable.
 For \eqref{3.1} and \eqref{3.2}, \eqref{3.3} is reduced to the Whittaker equation when $t=0$. 
The solution of \eqref{3.3} is given by
\end{theorem}
\begin{eqnarray}
\left(\begin{array}{c} y_1 \\ y_2 \end{array}\right)
&=&\left( \begin{array}{cc} L_{k,m}(x) & L_{k,-m}(x) \\
       \frac{-2k-2m-1}{u(2m+1)}L_{k+\frac{1}{2},m+\frac{1}{2}}(x) & \frac{-4m}{u}L_{k+\frac{1}{2},-m-\frac{1}{2}}(x)
          \end{array} 
\right),
\end{eqnarray}
where
\begin{eqnarray}
 L_{k,m}(x)&=& x^{2m+\frac{1}{2}} e^{-\frac{x^2}{2}} {} _1F_1 (m-k+\frac{1}{2},2m+1;x^2) \\
 &=&x^{2m+\frac{1}{2}} e^{-\frac{x^2}{2}} \sum_{n=0}^ 
 \infty\frac{\Gamma(2m+1)\Gamma(m-k+\frac{1}{2}+n)x^{2n}}{\Gamma(2m+1+n)\Gamma(m-k+\frac{1}{2})n!} .  
\end{eqnarray}

  \subsection{The linear monodromy }
The equation \eqref{linear:p4} has a regular singular point $x=0$ and an irregular singular point
$x=\infty$ with the Poincar\'e rank 2. We will define the linear monodromy  
$\{M_0,~\Gamma,~G_1,~G_2,~G_3,~G_4,$ $ e^{2\pi iT_{\infty}}\}$ of \eqref{linear:p4} \cite{FZ}, \cite{MJ}.  
\par\bigskip\noindent
1) At the regular singularity $x=0$, the local behavior of $Y(x)$ is given by 
\begin{eqnarray}
Y^{(0)}(x)&=&\left(1+O(x)\right)x^{T_0},
\end{eqnarray}
where
\begin{eqnarray}
 T_0&=&\left(
    \begin{array}{cc}
    \theta_0&0 \\
    0 &-\theta_0 
    \end{array}\right).
\end{eqnarray}
The local monodromy of $Y^{(0)}(x)$ around $x=0$ is 
\begin{eqnarray}
 M_0&=&e^{2\pi iT_0}.
\end{eqnarray}
\par\bigskip\noindent
2) At the irregular singularity $x=\infty$, a formal solution is given by
\begin{eqnarray}
Y^{(\infty)}&=&\left(1+\frac{Y_1}{x}+\cdots\right)e^{T(x)}, \\
T(x)&=& \left(
    \begin{array}{cc}
    1 &0 \\
    0 &-1
    \end{array}
\right)
\frac{x^2}{2}
+
\left(
    \begin{array}{cc}
    t &0 \\
    0 &-t
    \end{array}
\right)x
+\left(
    \begin{array}{cc}
    \theta_\infty &0 \\
    0 &-\theta_\infty
    \end{array}
\right)\log \frac{1}{x}, \\
Y_1&=&\frac{1}{2}\left(
    \begin{array}{cc}
    -H_{IV} &-u \\
    2(z-\theta_0-\theta_\infty)u^{-1} & H_{IV}
    \end{array}
\right), 
\end{eqnarray}
where 
\begin{eqnarray}
H_{IV}&=&\frac{2}{y}z^2-\left(y+2t+\frac{4}{y}\theta_0 \right)z+(\theta_0+\theta_\infty)(y+2t).
\end{eqnarray}
Since $x=\infty$ is an irregular singularity, the actual asymptotic behavior of $Y(x)$ changes the form
in the Stokes region of the complex $x$-plane: 
\begin{eqnarray}  
S_j&=& \left\{x\ \left|\frac{\pi}{2}(j-1)-\epsilon < \arg x < \frac{\pi}{2}j+\epsilon,\ |x|>R \right. \right\}, (j=1,2,3,4,5), 
\end{eqnarray} 
where $\epsilon$ is sufficiently small and $R$ is sufficiently large.

We denote $Y^{(j)}$ is a holomorphic solution in $S_j.$
According to the Stokes phenomenon, if $Y^{(j)} \sim Y^{(\infty)}(x)$ as  
$x \rightarrow \infty$ in $S_j$, then $Y^{(j+1)} = Y^{(j)}G_j$
 and $Y^{(5)}=Y^{(1)}e^{2\pi iT_{\infty}}$, where the matrices $G_j\ (1 \leq j \leq 4)$ are called the Stokes matrices
and $e^{2\pi iT_{\infty}}$ is a formal monodromy around $x=\infty$.

\par\bigskip\noindent
3) Connection matrix $\Gamma$ \\
 Since both $Y^{(0)}$ and $Y^{(1)}$ satisfy \eqref{linear:p4}, they are related by the connection matrix:
\begin{eqnarray}
Y^{(1)}&=&Y^{(0)}\Gamma. 
\end{eqnarray}
\par\bigskip\noindent
4) We have
\begin{eqnarray}
\Gamma^{-1}M_0 \Gamma G_1G_2G_3G_4e^{2i\pi T_{\infty}}&=&I_2. 
\end{eqnarray}
Generally, we cannot calculate $G_i$ and $\Gamma$.
By the isomonodromic condition, the linear monodromy is invariant for any t.
For the symmetric solution of the fourth Painlev\'e equation we can calculate the linear monodromy, 
because \eqref{linear:p4} is reduced to the Whittaker equation when $t=0$.  
\begin{theorem}
 For the symmetric solution \eqref{3.1} and \eqref{3.2} of the fourth Painlev\'e equation,
the linear monodromy is
\begin{eqnarray}
  M_0&=&\left(  \begin{array}{cc}  e^{2i\pi \theta_0} &0 \\  0 &e^{2i\pi(1-\theta_0)}  \end{array} \right)= \left( \begin{array}{cc} -e^{4mi\pi} &0 \\ 0 &-e^{-4mi\pi}  \end{array} \right),    \\
\Gamma&=& \left( \begin{array}{cc} \frac{\Gamma(-2m)}{\Gamma(\frac{1}{2}-m-k)} & \frac{\Gamma(-2m) e^{-i\pi (k+m+\frac{1}{2})}}{\Gamma( \frac{1}{2}-m+k)} \\ \frac{\Gamma(2m)}{\Gamma(\frac{1}{2}+m-k)} & \frac{\Gamma(2m) e^{-i\pi(k-m+\frac{1}{2})}} {\Gamma(\frac{1}{2}+m+k)} \end{array} \right),   \\
 G_1&=&\left( \begin{array}{cc} 1 &0 \\  \frac{2\pi e^{i\pi (\frac{-1}{2}+2k)}}{\Gamma(\frac{1}{2}-m-k) \Gamma(\frac{1}{2}+m-k)} &1 \end{array} \right),   \\
 G_2&=& \left( \begin{array}{cc} 1 &\frac{2\pi e^{i\pi (\frac{-1}{2}-4k)}}{\Gamma(\frac{1}{2}-m+k)\Gamma(\frac{1}{2}+m+k)} \\ 0 &1 \end{array} \right),  \\
G_3&=& \left( \begin{array}{cc} 1 &0 \\ \frac{2\pi e^{i\pi (\frac{-1}{2}+6k)}}{\Gamma(\frac{1}{2}-m-k) \Gamma(\frac{1}{2}+m-k)} &1 \end{array} \right),      \\
  G_4&=& \left( \begin{array}{cc} 1 &\frac{2\pi e^{i\pi (\frac{-1}{2}-8k)}}{\Gamma(\frac{1}{2}-m+k)\Gamma(\frac{1}{2}+m+k)} \\ 0 &1 \end{array} \right),   \\
  e^{2i\pi T_{\infty}}&=& \left( \begin{array}{cc} e^{2i\pi (1-\theta_\infty)} &0 \\ 0 &e^{2i\pi \theta_\infty} \end{array}  \right)=\left( \begin{array}{cc} -e^ {-4ki\pi} &0 \\ 0 & -e^{4ki\pi} \end{array} \right).
\end{eqnarray}
\end{theorem}
{\it Proof.}\ 
1) Two fundamental solutions $X_{k,m}(x)$ and $X_{-k,m}(x e^{\frac{-i\pi}{2}})$
in the Stokes region $S_j$ are expressed in the linear combination of $L_{k,m}(x)$ 
and $L_{k,-m}(x)$ \cite{JH}. 

For $r, s, t \in Z$,
 \begin{align}
  X_{k,m}(xe^{r i \pi})&=&\frac{\Gamma(-2m) e^{r i \pi \theta_0}L_{k,m}(x)}{\Gamma(\frac{1}{2}-m-k)}+\frac{\Gamma(2m) e^{ri \pi(1-\theta_0)}L_{k,-m}(x)}{\Gamma(\frac{1}{2}+m-k)}, \\
X_{k,m}(xe^{si \pi})&=&\frac{\Gamma(-2m) e^{s i \pi \theta_0}L_{k,m}(x)}{\Gamma(\frac{1}{2}-m-k)}+\frac{\Gamma(2m) e^{si \pi(1-\theta_0)}
L_{k,-m}(x)}{\Gamma(\frac{1}{2}+m-k)}, \\
X_{-k,m}(xe^{t i\pi-\frac{i\pi}{2}})&=&\frac{\Gamma(-2m) e^{i \pi \theta_0(t-\frac{1}{2}})L_{k,m}(x)}{\Gamma(\frac{1}{2}-m+k)}+\frac{\Gamma(2m) e^{i\pi(t-\frac{1}{2})(1-\theta_0)}L_{k,-m}(x)}{\Gamma(\frac{1}{2}+m+k)}
\end{align}
hold. 
Eliminating $L_{k,m}, L_{k,-m}$, and putting $s=0, t=0$ and $x \rightarrow x e^{-ri\pi}$, then we have 
\begin{eqnarray}
 X_{k,m}(x)& \sim& C_r e^{\frac{-x^2}{2}} x^{\theta_\infty-1}+D_r e^{\frac{x^2}{2}} 
 x^{-\theta_\infty},\label{x:1} \\
  && (r-\frac{1}{4}) \pi< \arg x <(r+\frac{3}{4}) \pi,\quad  (r=0,1,2, \cdots). \nonumber 
\end{eqnarray}
Similarly, we have 
\begin{eqnarray}
X_{-k,m}(x e^{\frac{-i \pi}{2}}) &\sim &E_r e^{\frac{-x^2}{2}}x^{\theta_\infty-1}+F_r 
e^{\frac{x^2}{2}}x^{-\theta_\infty},\label{x:2} \\
  &&\left( r-\frac{1}{4} \right) \pi < \arg x < (r+\frac{3}{4} ) \pi ,\quad(r=0,1,2, \cdots), \nonumber
\end{eqnarray}
where 
\begin{eqnarray}
 C_r&=& e^{r(1- \theta_\infty)i \pi} e^{\frac{ri \pi}{2}} \biggl[\frac{\sin 2(r+1)m \pi}{\sin 2m \pi}
+ e^{-2ki \pi} \frac{\sin 2rm \pi}{\sin 2m \pi} \biggr] ,\label{c:r} \\
 D_r&=& e^{(r+\frac{1}{2}) \theta_\infty i\pi} \frac{-2 \pi e^{\frac{\pi}{2} i(r+1)} 
e^{-ki\pi} e^{-\frac{i\pi}{4}} \sin 2rm \pi}
{\Gamma(\frac{1}{2}-m-k)\Gamma(\frac{1}{2}+m-k) \sin 2m \pi},\label{d:r}  \\
 E_r&=& e^{[r(1-\theta_\infty)-\frac{\theta_\infty}{2}]i\pi} \frac{ e^{\frac{\pi}{2}ir} e^{-\frac{i\pi}{4}} e^{-ki\pi}2 \pi \sin 2rm\pi}{\Gamma(\frac{1}{2}-m+k)\Gamma(\frac{1}{2}+m+k) \sin 2m\pi} , \\
 F_r&=&-e^{r \theta_\infty i\pi} e^{\frac{\pi}{2}ir} \left[\frac{\sin 2(r-1)m\pi}{ \sin 2m\pi}+ e^{-2ki\pi} \frac{ \sin 2rm\pi}{ \sin 2m\pi} \right]  . 
\end{eqnarray}
2) Stokes matrices $G_j$ \\
 For $r \pi < \arg x < (r+\frac{1}{2})\pi,(r \in Z)$, we write the coefficient matrix of 
\eqref{x:1} and \eqref{x:2} as
\begin{eqnarray}
      \left( \begin{array}{cc}  C_r & E_r \\ D_r & F_r
             \end{array}
      \right).
    \end{eqnarray}\\
 For $(r+\frac{1}{2})\pi < \arg x < (r+1)\pi$, we have
\begin{eqnarray}
      \left(  \begin{array}{cc} C_r &E_r \\ D_{r+1} & F_{r+1}
              \end{array}
      \right),
    \end{eqnarray}
                                                               
 \begin{eqnarray}
    G_{2r+1}\left( \begin{array}{cc} C_r &E_r \\ D_{r+1} &F_{r+1}
                  \end{array}  \right)
      =    \left( \begin{array}{cc} C_r &E_r \\ D_r &F_r 
                  \end{array}
           \right),
\end{eqnarray} \\
where \\
\begin{eqnarray}
   G_{2r+1}=\left( \begin{array}{cc} 1 &0 \\ T_{2r+1} &1
                   \end{array}
            \right),
\end{eqnarray}
\begin{eqnarray}  
 T_{2r+1}=\frac{D_r-D_{r+1}}{C_r}=\frac{F_r-F_{r+1}}{E_r} .
\end{eqnarray}
Substituting \eqref{c:r} and \eqref{d:r}, we have 
\begin{eqnarray}
T_{2r+1}=\frac{2 \pi e^{i\pi(\frac{-1}{2}+(4r+2)k)}}{\Gamma(\frac{1}{2}+m-k)\Gamma(\frac{1}{2}-m-k)},\quad(r=0,1,2, \cdots).
\end{eqnarray}
 In similar way, we have 
\begin{eqnarray}
   G_{2r}= \left( \begin{array}{cc} 1 & T_{2r} \\ 0 &1
                       \end{array}
                     \right),
\end{eqnarray}
\begin{eqnarray}
T_{2r}=\frac{2\pi e^{i\pi(\frac{-1}{2}-4rk)}}{\Gamma(\frac{1}{2}+m+k)\Gamma(\frac{1}{2}-m+k)},\quad(r=1,2, \cdots).
 \end{eqnarray}

\hfill $\boxed{}$

\noindent For special parameters, we have  
\begin{remark}\label{r:6}
 We set \ $2\theta_\infty-1=\alpha_0-\alpha_2,\ 2\theta_0=-\alpha_1$
and \ $\alpha_0+\alpha_1+\alpha_2=1$.\\
1) In case of $\alpha_0=0$, we have $m+k=-1/2$ and $G_2=G_4=I_2$.\\
2) In case of $\alpha_2=0$, we have $m-k=-1/2$ and $G_1=G_3=I_2$.\\
3) In case of $\alpha_0=0\ and\ \alpha_2=0$, we have $G_1=G_2=G_3=G_4=I_2$.\\
\end{remark}

 \subsection{Comparison with classical solutions}

 Umemura studied special solutions of the Painlev\'e equations \cite{HU1}. Umemura's classical 
solutions are either rational solution or the Riccati solution \cite{YM},\cite{MN1},\cite{HU2}. We show that the symmetric 
solution of the fourth Painlev\'e equation includes rational solutions
and one point of the Riccati solution of Umemura's classical solutions.
\par\bigskip\noindent 
{\bfseries 1) The Riccati solution} 
\par\medskip\noindent 
We set $p=y+2t-2w$. Then the system \eqref{p4:ham} is equivalent to
the following system:
\begin{eqnarray}
 \frac{dy}{dt}&=& 2yp -y^2-2ty+4\theta_0 ,\\
 \frac{dp}{dt}&=& 2yp -p^2+2tp+2(\theta_0-\theta_\infty+1). 
\end{eqnarray}
If $\alpha_2=0, \theta_0-\theta_\infty+1=0$. $p=0$ is a special solution and 
 $y$ satisfies the Riccati equation
\begin{eqnarray}
 \frac{dy}{dt}&=&-y^2-2ty+4\theta_0 ,
\end{eqnarray}
which is solved by the Weber function.
In this case, the linear monodromy is upper triangular matrices by Remark \ref{r:6} (2).  If $y(0)=0$ in (3.1), the Riccati solution is a symmetric solution.
We remark that the Riccati solutions have the same linear monodromy. 
\par\bigskip\noindent 
{\bfseries 2) Rational solutions}  
\par\medskip\noindent 
{\bfseries 2-1)} If $\alpha_0=\alpha_2=0, \theta_0=-1/2$.
    The Riccati equation is
\begin{eqnarray}
 \frac{dy}{dt}&=&y^2+2ty-2 ,
\end{eqnarray}
which has a rational solution $y=-2t$. 
This solution is reduced to the Hermite polynomial. The solution $(y,w)=(-2t,0)$
is a symmetric solution of the fourth Painlev\'e equation. 
In this case, every Stokes matrix is a unit matrix by Remark \ref{r:6} (3). 
\par\bigskip\noindent 
{\bfseries2-2)} If $\alpha_0=\alpha_1=\alpha_2=1/3,$ the fourth Painlev\'e
equation has an rational solution:
\begin{eqnarray}
y=\frac{-2t}{3} ,\qquad w=\frac{t}{3},
\end{eqnarray}
which is a symmetric solution of the fourth Painlev\'e equation. 
Since we have $(k,m)=(0,-1/3)$, \eqref{linear:p4} is reduced to the Airy function.

\subsection{Conclusion}
1) The symmetric solution of the fourth Painlev\'e equation exists for any parameter $\alpha$\
and\ $\beta$.

\par\bigskip\noindent 
2) There exist rational solutions and the Riccati solutions for the fourth Painlev\'e equation for
special parameters. Only for such special parameters, the symmetric solution coincides with 
Umemura's classical solution.
In this sense, the symmetric solution is a new special solution beyond Umemura's class.

\par\bigskip\noindent 
3) Two of four Stokes matrices ($G_1$ and $G_3$ or $G_2$ and $G_4$) become unit matrices
 when $\alpha_0$\ or\ $\alpha_2=0,$ 
 and every Stokes matrix becomes a unit matrix when $\alpha_0 =\alpha_2=0$. 
 Especially when $\alpha_2=0$, the linear monodromy become upper triangular matrices. 
 When $\alpha_0=\alpha_1=\alpha_2=1/3$ \ and\ $y=-2t/3$, the solution of the associated 
linear equation can be solved by the Airy function.

\section{The fifth Painlev\'e equation}\label{P5}
This section is based on the paper \cite{KOH}.
In this section, we will give three holomorphic solutions around $t=0$, which are 
invariant under the action of the B\"acklund transformation group. We will calculate 
the linear monodromy for one of these holomorphic solutions at $t=0$. The linear equation
is reduced  to the Gauss hypergeometric equation when $t=0$. The equation \eqref{5:16} 
has an 
irregular singularity at $x=\infty$ with the Poincar\'e rank 1, which becomes a 
regular singularity when $t=0$. We will show the extension of the isomonodromic 
deformation to $t=0$.

We will transform the linearization \eqref{5:16} to \eqref{0:r.fuchs} 
for the rational solution $y(t) \equiv -1$.  Therefore R.~Fuchs' observation 
is valid for $y(t) \equiv -1$.

\subsection{Meromorphic solutions around $t=0$}
Generic solutions of the fifth Painlev\'e equation have an essential singularity around $t=0$. 
Meromorphic solutions of the fifth equation around $t=0$ is classified in \S 37 \cite{GLS}. 
\begin{theorem}\label{p5:mero} 
1) Under some generic condition, 
the fifth Painlev\'e equation \eqref{p5} has a holomorphic solution around $t=0$: 
$$y=\pm\frac{\theta_0-\theta_1-\theta_{\infty}}{\theta_0-\theta_1+\theta_{\infty}} + \sum_{n=1}^\infty a_n t^n.$$
2) Assume that $\theta_0+\theta_1\notin \boldsymbol{Z} $. 
The fifth Painlev\'e equation \eqref{p5} has a holomorphic solution around $t=0$: 
$$y=1+ \frac {t}{1-\theta_0-\theta_1}+ \sum_{n=2}^\infty a_n t^n.$$

\end{theorem}

\noindent{\it Proof.}\quad  
Putting $t=0$ in  \eqref{2:y}, \eqref{2:z}, we have
the initial conditions $(y(0),w(0))$ as follows: 
\begin{eqnarray*}
 \left(\frac{\theta_0-\theta_1-\theta_\infty}{\theta_0-\theta_1+\theta_\infty},\quad
\frac{(\theta_0-\theta_1+\theta_\infty)(\theta_0+\theta_1+\theta_\infty)}{-4\theta_\infty}\right),\\
\left(\frac{-\theta_0+\theta_1+\theta_\infty}{\theta_0-\theta_1+\theta_\infty},\quad
\frac{\theta_1(\theta_0-\theta_1+\theta_\infty)}{-2(\theta_0-\theta_1)}\right), \quad\\
\left(1,\quad \frac{\theta_1(\theta_0+\theta_1+\theta_\infty)}{2(\theta_0+\theta_1)}\right).\quad\quad\quad\quad
\end{eqnarray*} 
Therefore, any holomorphic solution  $(y(t),w(t))$ has one of the above initial values. 
Since the system \eqref{2:y}, \eqref{2:z} is the Briot-Bouquet type, all these solutions converge. 
We will explain the Briot-Bouquet theorem in section 6.
\hfill $\boxed{}$

By the Hamiltonian form \eqref{2:y}, \eqref{2:z}, the first solution 
$$y(t)=\sum_{k=0}^\infty a_kt^k, \qquad 
w(t)=\sum_{k=0}^\infty b_kt^k, $$
 is expanded as  
\begin{eqnarray}
 y(t)&=&\sum_{k=0}^\infty a_kt^k,\quad a_0=\frac{\theta_0-\theta_1-\theta_\infty}{\theta_0-\theta_1+\theta_\infty},
\nonumber\label{I:y}     \\
&&a_1=\frac{a_0}{\Delta}\biggl[4b_0a_0^2-a_0(\theta_0+3\theta_1+\theta_\infty)
-4a_0b_0+2\theta_1+\theta_\infty-1\biggr],\nonumber\\
&&\cdots, 
\end{eqnarray}
\begin{eqnarray}
w(t)&=&\sum_{k=0}^\infty b_kt^k,\quad b_0=\frac{(\theta_0-\theta_1+\theta_\infty)(\theta_0+\theta_1+\theta_\infty)}{
-4\theta_\infty},\nonumber \label{I:w}    \\
&&b_1=\frac{b_0}{\Delta}\biggl[1-4a_0b_0+2b_0+2\theta_1+\theta_\infty\biggr],\quad \cdots, 
\end{eqnarray}
where
\begin{eqnarray*}
&&\Delta=\biggl[1+6a_0^2b_0-a_0(\theta_0+3\theta_1+\theta_\infty)-8a_0b_0+2\theta_1+\theta_\infty+2b_0\biggr]\nonumber\\
&&\times\biggl[-1+6a_0^2b_0-a_0(\theta_0+3\theta_1+\theta_\infty)-8a_0b_0+2\theta_1+\theta_\infty+2b_0\biggr]\nonumber\\
&&-2a_0(a_0-1)^2b_0\biggl[6a_0b_0-4b_0-(\theta_0+3\theta_1+\theta_\infty)\biggr].
\end{eqnarray*}
We denote this solution as (I). 

The second solution is expanded as
\begin{eqnarray*}
 y(t)&=&\sum_{k=0}^\infty a_kt^k,\quad a_0=\frac{-\theta_0+\theta_1+\theta_\infty}{\theta_0-\theta_1+\theta_\infty},
\nonumber\\
&&a_1=\frac{a_0}{\Delta}\biggl[4b_0a_0^2-a_0(\theta_0+3\theta_1+\theta_\infty)
-4a_0b_0+2\theta_1+\theta_\infty-1\biggr],\nonumber\\
&&\cdots, \\
w(t)&=&\sum_{k=0}^\infty b_kt^k,\quad b_0=\frac{(\theta_0-\theta_1+\theta_\infty)\theta_1}{
-2(\theta_0-\theta_1)},\nonumber \\
&&b_1=\frac{b_0}{\Delta}\biggl[1-4a_0b_0+2b_0+2\theta_1+\theta_\infty\biggr],\quad \cdots. 
\end{eqnarray*}
We denote this solution as (II).  We remark that $\Delta$ is a different function on 
$\theta_0,\theta_1,\theta_\infty$ in (I) and (II) although we use the same notation $\Delta$. 

The third solution
$$y(t)=\sum_{k=0}^\infty a_kt^k, \qquad 
w(t)=\sum_{k=0}^\infty b_kt^k, \qquad (\theta_0+\theta_1\notin \mathbb{Z})$$
is expanded as
\begin{eqnarray*}
&a_0=1,\quad a_1=\frac{1}{1-\theta_0-\theta_1},\quad a_2=\frac{1}{2-\theta_0-\theta_1}\biggl[a_1-2a_1^2b_0
+\frac{a_1^2}{2}(\theta_0+3\theta_1+\theta_\infty)\biggr],\cdots, \\
&b_0=\frac{(\theta_0+\theta_1+\theta_\infty)\theta_1}{2(\theta_0+\theta_1)},
\quad b_1=\frac{b_0}{1+\theta_0+\theta_1}\biggl[2a_1b_0-a_1(\theta_0+3\theta_1+\theta_\infty)-1\biggr],  \nonumber \\
&b_2=\frac{1}{2+\theta_0+\theta_1}\biggl[3a_1^2b_0^2+2a_2b_0^2-(a_2b_0+a_1b_1)
(\theta_0+3\theta_1+\theta_\infty)-b_1-8a_1b_1b_0\biggr],\cdots.
\end{eqnarray*}
We denote this solution as (III).

\begin{theorem} The three holomorphic solutions (I), (II) and (III) are invariant under 
the action of the B\"acklund transformation group. 
\[
  \begin{array}{l|| c| c| c| c | c| c}
     &\ s_0\ &\ s_1\ &\ s_2\ &\ s_3\ &\ \pi\ &\ \sigma \\
\hline 
    I \   &I&II&I&II&III&I \\
\hline 
    II \   &III&I&III&I&II&II \\
\hline 
    III \   &II&III&II&III&I&III \\
  \end{array}\]
\end{theorem}

We can prove the above theorem easily.

\subsection{The linear equation at $t=0$}

For a locally holomorphic solution (I) around $t=0$, we may extend the deformation 
equation to $t=0$ because $B(x,t)$ in \eqref{5:17} is holomorphic at $t=0$. 
Therefore   we can continue Miwa-Jimbo's isomonodromic deformation equation
to $t=0$. We describe more detail in subsection 4.6.

After substituting the solution (I) into the equation \eqref{5:16}, we put $t=0$. Then we have
\begin{eqnarray*}
\frac{\partial \Psi(x,0)}{\partial x}&=&A(x,0)\Psi(x,0),   \\ 
A(x,0)&=&\left. \biggl(\begin{array}{cc}
 \frac{1}{x}\left(
z+\frac{\theta_0}{2}\right)-\frac{1}{x-1}\left(z+\frac{\theta_0+\theta_\infty}{2}\right) 
& -\frac{u}{x}\left(z+\theta_0\right)-u\frac{(\theta_0-\theta_1-\theta_\infty)(\theta_0+\theta_1-\theta_\infty)}{
4\theta_\infty(x-1)} \\
    \frac{z}{u}\left(\frac{1}{x}-\frac{1}{x-1}\right) &-\frac{1}{x}\left(z+\frac{\theta_0}{2}\right)+
\frac{1}{x-1}\left(z+\frac{\theta_0+\theta_\infty}{2}\right)
    \end{array}
\biggr)\right|_{t=0},\nonumber 
 \end{eqnarray*}                            
which is reduced to the hypergeometric equation.  

The above discussion proves the following: 

\begin{theorem}     We can determine the linear monodromy of the   special solution (I).
For the solution (I), \eqref{5:16} is reduced 
to the hypergeometric equation  when $t=0$. 

The fundamental solution matrix is expressed as follows:
\begin{eqnarray*}
\Psi&=&\left(\begin{array}{cc}
     \psi_1^{(1)} & \psi_1^{(2)} \\
     \psi_2^{(1)} &\psi_2^{(2)}
 \end{array}  \right),  
\end{eqnarray*}
where
\begin{eqnarray*}
   \psi_1^{(1)}&=&x^{\frac{-\theta_0}{2}}(x-1)^\frac{-\theta_1}{2}{}_2F_1\left(\frac{\theta_\infty-\theta_0-\theta_1}{2},
1-\frac{\theta_\infty+\theta_0+\theta_1}{2},1-\theta_0;x\right) ,     \\
   \psi_1^{(2)}&=&x^{\frac{\theta_0}{2}}(x-1)^\frac{-\theta_1}{2}{}_2F_1\left(\frac{\theta_\infty+\theta_0-\theta_1}{2},
1-\frac{\theta_\infty-\theta_0+\theta_1}{2},1+\theta_0;x\right) ,     \\
   \psi_2^{(1)}&=&\frac{x^{\frac{-\theta_0}{2}}}{u_0}(x-1)^{\frac{-\theta_1}{2}}{}_2F_1\left(
\frac{\theta_\infty+\theta_0+\theta_1}{-2},
1+\frac{\theta_\infty-\theta_0-\theta_1}{2},1-\theta_0;x\right) ,     \\
   \psi_2^{(2)}&=&\frac{x^{\frac{\theta_0}{2}}}{u_0}(x-1)^{\frac{-\theta_1}{2}}{}_2F_1\left(
\frac{\theta_\infty-\theta_0+\theta_1}{-2},
1+\frac{\theta_\infty+\theta_0-\theta_1}{2},1+\theta_0;x\right).       
\end{eqnarray*}
Here $u_0=u(0)$.
\end{theorem}

Since Miwa-Jimbo's isomonodromic deformation equation can be continued to $t=0$,
the linear monodromy is invariant for any $t\in \mathbb{C}$.

\subsection{The linear monodromy}

\subsubsection{Miwa-Jimbo's linearization}
The equation \eqref{5:16} has two regular singular points $x=0$ and $x=1$, and an irregular singular point
$x=\infty$ with the Poincar\'e rank 1. We will define the linear monodromy  
$\{M_0,~M_1,~\Gamma_{0\infty},~\Gamma_{1\infty},~G_1,~G_2,~e^{2\pi iT_\infty}\}$ 
of \eqref{5:16}
following \cite{MJ}.

\par\bigskip\noindent
1) At the regular singularity $x=\nu, (\nu=0,1)$, the local behavior of $\Psi(x)$ is given by 
\begin{eqnarray*}
\Psi^{(\nu)}(x)&=&\left(1+O((x-\nu))\right)x^{T_\nu},
\end{eqnarray*}
where
\begin{eqnarray*}
 T_\nu&=&\left(
    \begin{array}{cc}
    \theta_\nu/2&0 \\
    0 &-\theta_\nu/2 
    \end{array}\right).
\end{eqnarray*}
The local monodromy of $\Psi^{(\nu)}(x)$ around $x=\nu$ is 
\begin{eqnarray*}
 M_\nu&=&e^{2\pi iT_\nu}.
\end{eqnarray*}
2) At the irregular singularity $x=\infty$, a formal solution is given by
\begin{eqnarray*}
\Psi^{(\infty)}&=&\left(1+\frac{\Psi_1}{x}+\cdots\right)e^{T(x)}, \\
T(x)&=& 
\frac{1}{2}\left(
    \begin{array}{cc}
    t &0 \\
    0 &-t
    \end{array}
\right)x
+\frac{1}{2}\left(
    \begin{array}{cc}
    \theta_\infty &0 \\
    0 &-\theta_\infty
    \end{array}
\right)\log \frac{1}{x}, \\
\Psi_1&=&\left(
    \begin{array}{cc}
    -H_{V} & \frac{u}{t}\left[z+\theta_0-y\left(z+\frac{\theta_0-\theta_1+\theta_\infty}{2}\right)\right]\\
    \frac{1}{tu}\left[z-\frac{1}{y}\left(z+\frac{\theta_0+\theta_1+\theta_\infty}{2}\right)\right]& H_{V}
    \end{array}
\right),
\end{eqnarray*}
where 
\begin{eqnarray*}
H_{V}&=&\frac{-1}{t}\left[z-\frac{1}{y}\left(z+\frac{\theta_0+\theta_1+\theta_\infty}{2}\right)\right]
\left[z+\theta_0-y\left(z+\frac{\theta_0-\theta_1+\theta_\infty}{2}\right)\right]  \nonumber\\
&&-z-\frac{\theta_0+\theta_\infty}{2}.
\end{eqnarray*}
Since $x=\infty$ is an irregular singularity, the actual asymptotic behavior of $\Psi(x)$ changes the form
in the Stokes region of the complex $x$-plane: 
\begin{eqnarray*}  
S_j&=& \left\{x\ |\ \pi(j-1)-\varepsilon < \arg (xt) < \pi j+\varepsilon,|x|>R  \right\}\quad (j=1,2,3), 
\end{eqnarray*} 
where $\varepsilon$ is sufficiently small and $R$ is sufficiently large. 

We denote $\Psi_{\infty}^{(j)}$ is a holomorphic solution in $S_j.$  
According to the Stokes phenomenon, if 
\begin{eqnarray*} 
   \Psi_{\infty}^{(j)}&\sim&\Psi^{(\infty)}(x) \quad \textrm{as} \quad x \rightarrow \infty \quad \textrm{in} \quad S_j,\\
  \Psi_{\infty}^{(2)}&=&\Psi_{\infty}^{(1)}G_1, \quad \Psi_{\infty}^{(3)}=\Psi_{\infty}^{(2)}G_2, 
  \quad \Psi_{\infty}^{(3)}=\Psi_{\infty}^{(1)}e^{2\pi iT_{\infty}},
\end{eqnarray*}
 where the matrices $G_j,(j=1,2)$ are called the Stokes matrices
and $e^{2\pi iT_{\infty}}$ is a formal monodromy around $x=\infty$.

\par\noindent
3)  
 Since both $\Psi^{(0)}$, $\Psi^{(1)}$ and $\Psi_{\infty}^{(1)}$  satisfy \eqref{5:16}, 
 they are related by the connection matrix $\Gamma_{\nu\infty}$:
\begin{equation*}
\Psi^{(\nu)}=\Psi_{\infty}^{(0)}\Gamma_{\nu \infty},\quad (\nu=0,1).  
\end{equation*}
4) We have
\begin{equation*}
\Gamma_{0\infty}M_0\Gamma_{0\infty}^{-1}\Gamma_{1\infty}M_1\Gamma_{1\infty}^{-1} 
G_1G_2e^{2i\pi T_{\infty}}=I_2. 
\end{equation*}
Generally, we cannot calculate $G_j$ and $\Gamma_{\nu\infty}$.
By the isomonodromic condition, the linear monodromy is invariant for any $t$.
For the solution (I) of the fifth Painlev\'e equation, we can calculate the linear monodromy, 
because \eqref{5:16} is reduced to the hypergeometric equation when $t=0$.  

 \begin{theorem}\label{t:10}   For the the solution (I) of the fifth Painlev\'e equation,
the linear monodromy is
\begin{eqnarray*}
  M_0&=&\left(  \begin{array}{cc}  e^{-i\pi \theta_0} &0 \\  0 &e^{i\pi\theta_0}  \end{array} \right),  \quad
  M_1=\left(  \begin{array}{cc}  e^{-i\pi \theta_1} &0 \\  0 &e^{i\pi\theta_1}  \end{array} \right),    \\
\Gamma_{0\infty}&=& \left( \begin{array}{cc} 
\frac{e^{\frac{i\pi}{-2}(\theta_\infty-\theta_1-\theta_0)} 
\Gamma(1-\theta_0)\Gamma(1-\theta_\infty)}{\Gamma(1-\frac{\theta_\infty+\theta_1+\theta_0}{2})
\Gamma(1-\frac{\theta_\infty-\theta_1+\theta_0}{2})} & 
\frac{e^{\frac{i\pi}{-2}(\theta_\infty-\theta_1+\theta_0)}\Gamma(1+\theta_0)\Gamma(1-\theta_\infty)}
{\Gamma(1-\frac{\theta_\infty-\theta_1-\theta_0}{2})\Gamma(1-\frac{\theta_\infty+\theta_1-\theta_0}{2})} \\ 
\frac{-e^{\frac{i\pi}{2}(\theta_\infty+\theta_1+\theta_0)}\Gamma(1-\theta_0)\Gamma(\theta_\infty-1)}
{\Gamma(\frac{\theta_\infty-\theta_1-\theta_0}{2})\Gamma(\frac{\theta_\infty+\theta_1-\theta_0}{2})} & 
\frac{-e^{\frac{i\pi}{2}(\theta_\infty+\theta_1-\theta_0)}\Gamma(1+\theta_0)\Gamma(\theta_\infty-1)}
{\Gamma(\frac{\theta_\infty-\theta_1+\theta_0}{2})\Gamma(\frac{\theta_\infty+\theta_1+\theta_0}{2})} 
 \end{array} \right),   \\
\Gamma_{1\infty}&=& \left( \begin{array}{cc} 
\frac{\Gamma(1-\theta_1)\Gamma(1-\theta_\infty)}{\Gamma(1-\frac{\theta_\infty+\theta_1+\theta_0}{2})
\Gamma(1-\frac{\theta_\infty+\theta_1-\theta_0}{2})} & 
\frac{e^{i\pi\theta_1}\Gamma(1+\theta_1)\Gamma(1-\theta_\infty)}
{\Gamma(1-\frac{\theta_\infty-\theta_1+\theta_0}{2})\Gamma(1-\frac{\theta_\infty-\theta_1-\theta_0}{2})} \\ 
\frac{\Gamma(1-\theta_1)\Gamma(\theta_\infty-1)}
{\Gamma(\frac{\theta_\infty-\theta_1-\theta_0}{2})\Gamma(\frac{\theta_\infty-\theta_1+\theta_0}{2})} & 
\frac{e^{i\pi\theta_1}\Gamma(1+\theta_1)\Gamma(\theta_\infty-1)}
{\Gamma(\frac{\theta_\infty+\theta_1-\theta_0}{2})\Gamma(\frac{\theta_\infty+\theta_1+\theta_0}{2})} 
 \end{array} \right),   \\
 G_1&=&G_2=I_2 ,\\
  e^{2i\pi T_\infty}&=& \left( \begin{array}{cc} e^{i\pi \theta_\infty} &0 \\ 0 &e^{-i\pi \theta_\infty} 
\end{array}  \right).
\end{eqnarray*}
\end{theorem}

\noindent
{\it Remark. }\  While  $x=\infty$ is an irregular singularity with the Poincar\'e rank 1 in   \eqref{5:16},
 $x=\infty$ becomes the regular singularity when $t=0$. 
This  means that the formal solution around the irregular singularity $x=\infty$, which is expressed
in the form of an asymptotically expanded power series converges for any $t$
by the isomonodromic condition. Therefore, every Stokes matrix becomes the unit matrix. 
It is difficult to prove this fact directly but we prove this for the special value of parameters:
$\alpha+\beta=0,\, \gamma=0\, (\theta_0=\theta_1=1/2)$, in section five.

\par\bigskip\noindent
{\it Remark.}\
The Stokes multipliers become zero 
for our solutions, which are analytic around zero. Y.~Sibuya studied differential equations 
whose Stokes multipliers vanish at irregular singular points 
\cite{Sibuya} (Professor Okamoto taught us Sibuya's paper). Although he did not consider 
isomonodromic deformations,   we think that the isomonodromic deformation equations 
become simple when Stokes multipliers vanish. 



\subsection{Comparison with classical solutions}
 Umemura studied special solutions of the Painlev\'e equations \cite{HU1}, which are called 
 {\it classical solutions}. Umemura's classical 
solutions are either   algebraic solutions or the Riccati solutions \cite{YM},\cite{MN1},\cite{HU2}. 

We show that the new special solution (I) includes an algebraic solution $y \equiv -1$
and one point of the Riccati solution. Since (II) and (III) are the B\"acklund transforms of (I), 
(II) and (III) also contain classical solutions. 

We have the following Riccati solutions:
\begin{enumerate}
\item In case of $\theta_0+\theta_1+\theta_\infty=0$, we have $w \equiv 0$ from \eqref{2:z} and \eqref{I:w},  and
      $y$ satisfies the  Riccati equation.
\item In case of $\theta_0-\theta_1-\theta_\infty=0$, we have $y \equiv 0$ from \eqref{2:y} and \eqref{I:y} , and
      $w$ satisfies the  Riccati equation.
\item In case of $\theta_0+\theta_1-\theta_\infty=0$, all monodromy data become upper  
     half triangular matrices by  Theorem  \ref{t:10}. 
\item  In case of $\theta_0+\theta_1-\theta_\infty=-2$, all monodromy data become lower 
     half triangular matrices by   Theorem  \ref{t:10}. 
\end{enumerate} 
In every case above,  (I) includes one point of the Riccati solution.
We remark that the Riccati solutions have the same linear monodromy.

In case of $\alpha +\beta=0$ and $\gamma=0$ (i.e. $\theta_0=\theta_1=1/2$), 
the system  \eqref{2:y} and \eqref{2:z} has a special 
solution $y \equiv -1$ and $w=\frac{1+\theta_\infty}{-4}+\frac{t}{8}$, which is a 
rational solution of $P_V$. We will study this rational solution in the next section.

We remark that (III) contains the Riccati solution  $y=e^t$ 
for $\theta_0=\theta_1=\theta_{\infty}=0$.

\subsection{R.~Fuchs' observation for the solution $y \equiv -1$}
For an algebraic solution  of the Painlev\'e equation, R.~Fuchs 
observed that the associated linear equation can be transformed by 
an appropriate variable change to an equation which does not include 
the deformation parameter $t$.  He showed that
the linear equation for   special Picard's solutions, which correspond to
three, four and six divided points of elliptic curves, can be reduced to the hypergeometric equation 
\cite{RF2}.  

In this section we will show that  R.~Fuchs' observation is true  for  the rational solution 
 $y \equiv -1$ for $\theta_0=\theta_1=1/2$ of the fifth Painlev\'e equation. 
 This solution is a special case  of (I), as we claimed in the previous section. 
 For a generic parameter, we can directly calculate the linear monodromy of (I) only for $t=0$. 
 But in   case of $\theta_0=\theta_1=1/2$,  we can directly calculate the linear monodromy of  $y \equiv -1$
  for generic $t\in\mathbb{C}$. The authors learned the method in this section from Professor Kazuo Okamoto.

We substitute the solution $y \equiv -1$ into  Miwa-Jimbo's isomonodromic deformation equations 
\eqref{5:16} and \eqref{5:17}: 
\begin{eqnarray}\label{5:1}
\frac{\partial \Psi(x,t)}{\partial x}&=&A(x,t)\Psi(x,t),   \\ 
A(x,t)&=&
  \frac{1}{2}\left(
    \begin{array}{cc}
    t &0 \\
    0 &-t
    \end{array}
\right)
-
\frac{1}{x}\left(
    \begin{array}{cc}
    \frac{\theta_\infty}{4}+\frac{t}{8} & u\left(\frac{1-\theta_\infty}{4}-\frac{t}{8}\right) \\
    u^{-1}\left(\frac{1+\theta_\infty}{4}+\frac{t}{8}\right) & \frac{-\theta_\infty}{4}-\frac{t}{8}
    \end{array}
\right)\nonumber \\
& &-\frac{1}{x-1}
\left(
    \begin{array}{cc}
    \frac{\theta_\infty}{4}-\frac{t}{8} & u\left(\frac{-1+\theta_\infty}{4}-\frac{t}{8} \right) \\
    -\frac{1}{u}\left(\frac{1+\theta_\infty}{4}-\frac{t}{8}\right)
 & -\frac{\theta_\infty}{4}+\frac{t}{8}
    \end{array}
\right),    \nonumber \\
 \frac{\partial \Psi(x,t)}{\partial t}&=&\left(
    \begin{array}{cc}
    \frac{x}{2} &\frac{u}{4} \\
    -\frac{u^{-1}}{4} &-\frac{x}{2} 
    \end{array}
\right)\Psi(x,t).
\end{eqnarray}
Putting
\begin{eqnarray*}
\Psi=\left(\begin{array}{c}
\psi_1  \\ \psi_2
\end{array} \right),
\end{eqnarray*}
we have equations for $\psi_1$:
\begin{eqnarray}
\frac{\partial ^2\psi_1}{\partial x^2}&+&\left[\frac{1}{x}+\frac{1}{x-1}-\frac{2t}{2tx-t+2(1-\theta_\infty)}\right]\frac
{\partial \psi_1}{\partial x} \nonumber\\ 
&-&\biggl[\frac{t^2}{4}+\frac{1}{16x^2}+\frac{1}{16(x-1)^2}+\frac{4(1-\theta_\infty)t^2}{t^2-4(1-\theta_\infty)^2}
\cdot \frac{1}{2tx-t+2(1-\theta_\infty)}   \nonumber\\
&+&\left(\frac{t}{2}+\frac{\theta_\infty}{2} -(\frac{\theta_\infty}{4}+\frac{t}{8})
\frac{2t}{t-2(1-\theta_\infty)}+\frac{1}{8}-\frac{t^2}{16}-\frac{\theta_\infty^2}{4}-\frac{t\theta_\infty}{4}\right)
\frac{1}{x}  \nonumber\\
&+&\left(\frac{t}{2}-\frac{\theta_\infty}{2}+(\frac{\theta_\infty}{4}-\frac{t}{8})\frac{2t}
{t+2(1-\theta_\infty)}-\frac{1}{8}+\frac{t^2}{16}+\frac{\theta_\infty^2}{4}-\frac{t\theta_\infty}{4}\right)
\frac{1}{x-1}\biggr]\psi_1=0 ,\nonumber \label{p:1} \\
\end{eqnarray}
\begin{equation*}
\biggl[\frac{1-\theta_\infty}{x(x-1)}+\frac{t}{2}\cdot\frac{(2x-1)}{x(x-1)}\biggr]\frac{\partial 
\psi_1}{\partial t} -\frac{\partial \psi_1}{\partial x} 
=\left[\frac{1}{x}\left(\frac{\theta_\infty}{4}+\frac{t}{8}\right)-\frac{1}{x-1}\left(\frac{1}{2}-\frac{\theta_\infty}{4}
+\frac{t}{8}\right)\right]\psi_1 .\nonumber \label{p:11} \\
\end{equation*}
The equation \eqref{p:1} has three regular singularities; $x=0, \quad x=1$ and 
$x=\frac{1}{2}-\frac{1-\theta_\infty}{t}$ and
an irregular singularities at $x=\infty$ with the Poincar\'e rank 1.  We remark that 
$x=\frac{1}{2}-\frac{1-\theta_\infty}{t}$ is an apparent  singularity.

We take new variables:
\begin{eqnarray*}
    (\psi_1,x)&\longrightarrow&(\phi_1,\xi): \\
\psi_1&=&\phi_1\left[x(x-1)\right]^{\frac{-1}{4}}\left[x-\frac{1}{2}+\sqrt 
{x(x-1)}\right]e^{\frac{t}{4}+\frac{t}{1-\theta_\infty} \sqrt{x(x-1)}}, \\
\xi&=&\left(x-\frac{1}{2}+\sqrt{x(x-1)}\right)e^{\frac{t}{1-\theta_\infty}\sqrt{x(x-1)}}.
\end{eqnarray*}
Then \eqref{p:1} is reduced to 
\begin{eqnarray}
\frac{\partial ^2\phi_1}{\partial \xi^2}+\frac{3}{\xi}\frac{\partial \phi_1}{\partial \xi}+\left[1-\left(\frac{1-\theta_\infty}{2}\right)^2\right]
\frac{\phi_1}{\xi^2}=0,\label{p:2}
\end{eqnarray} 
which is independent of $t$. 
We can solve \eqref{p:2} easily:
\begin{eqnarray*}
\phi_1=c_1\xi^{\alpha_1}+c_2\xi^{\alpha_2},\quad \left(\alpha_1,\alpha_2=-1\pm\frac{1-\theta_\infty}{2}\right).
\end{eqnarray*}
Therefore solutions of \eqref{p:1} are given by
\begin{eqnarray*}
\psi_1=\left[x(x-1)e^{-t}\right]^{\frac{-1}{4}}\biggl[c_1\left(\sqrt{x}+\sqrt{(x-1)}\right)^{1-\theta_\infty} 
e^{\frac{t}{2}\sqrt{x(x-1)}}\nonumber \\
+c_2\left(\sqrt{x}+\sqrt{x-1}\right)^{-1+\theta_\infty}e^{\frac{-t}{2}\sqrt{x(x-1)}}\biggr],
\end{eqnarray*}
where $c_1,c_2$ are constants. 
\par\bigskip
We notice that if  $\theta_0=1/2, \theta_1=1/2$  and $y=-1$, 
we have $z=-(t+2\theta_\infty +2)/8$ and $u= -c^{-1} e^{t/2} $ for a constant $c$.
The fundamental solution $\Psi$ of  \eqref{5:1} is 
\begin{align}
 x^{\frac{1}{4}}&(x-1)^{\frac{1}{4}} \Psi    \nonumber\\
= &\left(\begin{array}{cc}
e^{\frac{t}{4}+ \frac{t}{2}\sqrt{x(x-1)}} \left(\sqrt{x}+\sqrt{x-1}\right)^{1-\theta_\infty} & 
e^{\frac{t}{4}- \frac{t}{2}\sqrt{x(x-1)}} \left(\sqrt{x}+\sqrt{x-1}\right)^{-1+\theta_\infty} \\ 
c\ e^{-\frac{t}{4}+ \frac{t}{2}\sqrt{x(x-1)}}  \left(\sqrt{x}+\sqrt{x-1}\right)^{-1-\theta_\infty} &
c\ e^{-\frac{t}{4}- \frac{t}{2}\sqrt{x(x-1)}}  \left(\sqrt{x}+\sqrt{x-1}\right)^{1+\theta_\infty}
\end{array}\right).  \label{p:3} 
\end{align}
The solution \eqref{p:3}  has  a regular singularity at $x=\infty$ if we put $t=0$.  Although  
\eqref{p:3} has irregular singularity at $x=\infty$ in case of $t\not=0$,
every Stokes matrix becomes a unit matrix since they give convergent series around $x=\infty$.

The linear monodromy of the fundamental solution \eqref{p:3} is
\begin{eqnarray*}
M_0&=&\left(\begin{array}{cc} 0 & ie^{i\pi \theta_\infty} \\ ie^{-i\pi\theta_\infty} & 0 \end{array}\right),\quad
M_1=\left(\begin{array}{cc} 0 & -i \\ -i & 0 \end{array}\right), \\
\Gamma_{0\infty}&=&\Gamma_{1\infty}=I_2,\quad G_1=G_2=I_2, \\
e^{2i\pi T_\infty}&=&\left(\begin{array}{cc} e^{-i\pi\theta_\infty} & 0 \\ 0 & e^{i\pi\theta_\infty}
 \end{array}\right).
\end{eqnarray*}
We have
\begin{equation*}
M_0M_1e^{2\pi i\theta_\infty}=I_2 . 
\end{equation*}

\subsection{Extension of deformation to $t=0$}
In this section, we will show that the fundamental solution $\Psi(x,t)$ of \eqref{5:16} 
exists for any $t \in \mathbb{C}$. The equation \eqref{5:16} has an irregular singularity at 
$x=\infty$ with the Poincar\'e rank 1, which turns out the regular singularity when $t=0$. 
For the special solution (I),  $B(x,t)$ in \eqref{5:17} is holomorphic at $t=0$.  Therefore
we have a fundamental solution $\Psi(x,t)$ which is analytic on $(x,t)$ and has 
a branch along $x=\infty$.  

We set the Pauli matrix
$$\sigma_3=\left(\begin{array}{cc} 1 & 0 \\ 0 & -1 \end{array}\right).$$
The following theorem assures   that the  isomonodromic deformation extends to $t=0$.
\begin{theorem}
For the special solution (I), we have a fundamental solution at $x=\infty$
\begin{equation}\label{6:sol}
\Psi(x,t) = \left(1+\frac{\Psi_1(t)}{x}+ \frac{\Psi_2(t)}{x^2}+  \cdots\right)e^{T(x)},
\end{equation}
where 
\begin{equation*}
T(x)= \left(\frac{t}{2}x-\frac{\theta_\infty}{2}\log x\right)\sigma_3. 
\end{equation*}
Here $\Psi_j(t)$ is holomorphic around $t=0$ for $j=1,2,3,\cdots$.
\end{theorem}
{\it Proof.}\quad 
We write equation \eqref{5:17} as follows:
\begin{eqnarray*}
\frac{\partial \Psi(x,t)}{\partial t}&=&\left[\frac{x}{2}\sigma_3+\frac{1}{t}
\left(\begin{array}{cc}
0 & b \\ c & 0 \end{array}\right)\right]\Psi(x,t),  
\end{eqnarray*} 
where
\begin{eqnarray*}
b&=&-u\left[z+\theta_0-y\left(z+\frac{\theta_0-\theta_1+\theta_\infty}{2}\right)\right], \\
c&=&u^{-1}\left[z-\frac{1}{y}\left(z+\frac{\theta_0+\theta_1+\theta_\infty}{2}\right)\right].
\end{eqnarray*}
For the solution (I),  we have $b(0)=c(0)=0$.

At the irregular singularity $x=\infty$, a formal solution \eqref{6:sol} exists for $t\not=0$.
Therefore, we have
\begin{eqnarray*}
\frac{\partial \Psi^{(\infty)}}{\partial t}&=&\biggl[\left(1+\frac{\Psi_1'}{x}+\frac{\Psi_2'}{x^2}+\cdots\right)
+\left(1+\frac{\Psi_1}{x}+\frac{\Psi_2}{x^2}+\cdots \right)\frac{x}{2}\sigma_3 \biggr]e^{T(x)}\nonumber\\
&=&\left[\frac{x}{2}\sigma_3+\frac{1}{t}
\left(\begin{array}{cc}
0 & b \\ c & 0 \end{array}\right)\right]\left(1+\frac{\Psi_1}{x}+\frac{\Psi_2}{x^2}+\cdots \right)e^{T(x)},
\end{eqnarray*}
where $'$ means a derivation by $t$.
\begin{eqnarray*}
\biggl[1&+&\frac{\Psi_1'}{x}+\frac{\Psi_2'}{x^2}+\cdots
+\left(1+\frac{\Psi_1}{x}+\frac{\Psi_2}{x^2}+\cdots \right)\frac{x}{2}\sigma_3 \biggr] \nonumber\\
&=&\left[\frac{x}{2}\sigma_3+\frac{1}{t}
\left(\begin{array}{cc}
0 & b \\ c & 0 \end{array}\right)\right]\left(1+\frac{\Psi_1}{x}+\frac{\Psi_2}{x^2}+\cdots \right).
\end{eqnarray*}
We put $\Psi_n=\left(\begin{array}{cc} a_n & b_n \\ c_n & d_n \end{array}\right)$ and compare 
the coefficients of the equal degree of $x$ in both sides:

\par\medskip\noindent 
1)\quad $x^{0}$:
We have
\begin{eqnarray*}
\frac{1}{2}\left[\Psi_1,\sigma_3\right]=\frac{1}{t}\left(
\begin{array}{cc} 0 & b \\ c & 0 \\ \end{array}\right).
\end{eqnarray*}
Therefore
\begin{eqnarray*}
\left(\begin{array}{cc} 0 & -b_1 \\ c_1 & 0 \\ \end{array}\right)=\frac{1}{t}\left(
\begin{array}{cc} 0 & b \\ c & 0 \\ \end{array}\right). 
\end{eqnarray*}
Since $b(0)=0, c(0)=0$, $b_1$ and $c_1$ are holomorphic around $t=0$.

\par\medskip\noindent 
2)\quad $x^{-1}$: We have
\begin{eqnarray*}
\frac{1}{2}\left[\Psi_2,\sigma_3\right]&=&\frac{1}{t}\left(
\begin{array}{cc} 0 & b \\ c & 0 \\ \end{array}\right)\Psi_1-\Psi_1'.
\end{eqnarray*}
Therefore 
\begin{eqnarray*}
\left(
\begin{array}{cc} 0 & -b_2 \\ c_2 & 0 \\ \end{array}\right)&=&\frac{1}{t}\left(
\begin{array}{cc} 0 & b \\ c & 0 \\ \end{array}\right)
\left(\begin{array}{cc} a_1 & b_1 \\ c_1 & d_1 \\ \end{array}\right)-
\left(\begin{array}{cc} a_1' & b_1' \\ c_1' & d_1' \end{array}\right) \nonumber\\
&=&\left(\begin{array}{cc}  c_1 b/t - a_1' & d_1 b/t- b_1' \\
\ a_1 c /t- c_1' & b_1 c/t -d_1'\end{array}\right).
\end{eqnarray*}
Compared with the diagonal components, $a_1'$ and $d_1'$ are holomorphic because 
$b(0)=0$ and $c(0)=0$. 
Therefore $\Psi_1 =\left(\begin{array}{cc} a_1  & b_1  \\ c_1  & d_1  \\ \end{array}\right)$
is holomorphic around $t=0$. 

Compared with the off-diagonal components, $b_2$ and $c_2$ are holomorphic. 

\par\medskip\noindent 
3)\quad $x^{-n}$: We have
\begin{eqnarray*}
\frac{1}{2}\left[\Psi_{n+1}, \sigma_3\right]&=&\frac{1}{t}
\left(\begin{array}{cc}0 & b \\ c & 0 \end{array}\right)\Psi_n-\Psi_n'.
\end{eqnarray*}
Therefore 
\begin{eqnarray*}
\left(
\begin{array}{cc} 0 & -b_{n+1} \\ c_{n+1} & 0 \\ \end{array}\right)&=&\frac{1}{t}\left(
\begin{array}{cc} 0 & b \\ c & 0 \\ \end{array}\right)
\left(\begin{array}{cc} a_n & b_n \\ c_n & d_n \\ \end{array}\right)-
\left(\begin{array}{cc} a_n' & b_n' \\ c_n' & d_n' \end{array}\right) \nonumber\\
&=&\left(\begin{array}{cc} c_n b/t- a_n' & d_n b/t- b_n' \\
                           a_n c/t- c_n' & b_n c/t- d_n'\end{array}\right).
\end{eqnarray*}
In the same way,  $a_n', d_n'$ and $b_{n+1}, c_{n+1}$ are holomorphic. 
Therefore $\Psi_n =\left(\begin{array}{cc} a_n  & b_n \\ c_n  & d_n  \\ \end{array}\right)$
is holomorphic around $t=0$.  \hfill $\boxed{}$

\par\bigskip
For $\alpha+\beta=0,\gamma=0$ we showed that we may put $t=0$ in section 4.2.

\section{The sixth Painlev\'e equation} \label{P6}
This section is based on the paper \cite{KK2}.
In this section, we will give four meromorphic solutions around each fixed singularity
$t=0, 1, \infty$, respectively, which are transformed each other by the action of the 
B\"acklund transformation group. We will calculate the linear monodromy for one of these 
meromorphic solutions at $t=0$ by Jimbo's method given in \cite{MJ2}. We take two 
confluences of singularities of the linear equation.
One is the confluence between $x=0$ and $x=t$ and the other is the confluence 
between $x=1$ and $x=\infty$. For the former, the linear equation is reduced to the Gauss
hypergeometric equation and for the latter, it is reduced to a Heun's
type equation whose general solution can be obtained as a linear combination of two
monomials.
From these two confluences we obtain the linear monodromy for our solution explicitly. 
We will give the comparison with Umemura's classical solutions.  

\subsection{Meromorphic solutions around the fixed  singularities} 
In this section we will classify all of the meromorphic solutions around a fixed singularity.
We consider a solution of \eqref{y:0} and \eqref{z:0} (and that of \eqref{S:1} and 
\eqref{S:2} 
simultaneously) around $t=0$:
\begin{eqnarray}
y(t)=t^l\sum_{i=0}^{\infty}a_i t^i, \quad\bar z(t)=t^k\sum_{i=0}^{\infty}b_i t^i,
\quad z(t)=t^n\sum_{i=0}^{\infty}c_i t^i\
  (l,k,n\in \mathbb {Z}).
\end{eqnarray}
%
\begin{theorem} \label{m0:1}
For generic values of parameters, the sixth Painlev\'e equation has the following four
meromorphic solutions around $t=0$: 
\begin{eqnarray}
{\rm (0{\textrm -}I)}:\quad y(t)&=&\frac{\alpha_4}{\alpha_4-\alpha_0}t
+\frac{\alpha_0\alpha_4\left[-1-\alpha_1^2+\alpha_3^2+(\alpha_4-\alpha_0)^2\right]}
{2\left[1-(\alpha_4-\alpha_0)^2\right](\alpha_4-\alpha_0)^2}t^2+O(t^3),\\
\bar z(t)&=&\frac{1-\alpha_1^2+\alpha_3^2-(\alpha_4-\alpha_0)^2}
{4\left[1-(\alpha_4-\alpha_0)^2\right]}+O(t),\\
z(t)&=&\frac{\alpha_4-\alpha_0}{t}+O(t^0),
\end{eqnarray}
\begin{eqnarray}
{\rm (0{\textrm -}II)}:\quad y(t)&=&\frac{\alpha_4}{\alpha_4+\alpha_0}t
+\frac{-\alpha_0\alpha_4\left[1+\alpha_1^2-\alpha_3^2-(\alpha_4+\alpha_0)^2\right]}
{2\left[1-(\alpha_4+\alpha_0)^2\right](\alpha_4+\alpha_0)^2}t^2+O(t^3),\\
\bar z(t)&=&\frac{1-\alpha_1^2+\alpha_3^2-(\alpha_4+\alpha_0)^2}
{4\left[1-(\alpha_4+\alpha_0)^2\right]}+O(t),\\
 z(t)&=&\frac{\alpha_2(\alpha_1+\alpha_2)}
{1-\alpha_4-\alpha_0}+O(t),\\
{\rm (0{\textrm -}III)}:\quad y(t)&=&\frac{\alpha_1+\alpha_3}{\alpha_1}
+\frac{-\alpha_3\left[1+\alpha_4^2-\alpha_0^2-(\alpha_1+\alpha_3)^2\right]}
{2\alpha_1\left[1-(\alpha_1+\alpha_3)^2\right]}t+O(t^2),\\
\bar z(t)&=&\frac{-\alpha_1}{2(\alpha_1+\alpha_3)}+O(t),\quad
z(t)=\frac{-\alpha_1\alpha_2}{\alpha_1+\alpha_3}+O(t),\\
 {\rm (0{\textrm -}IV)}:\quad y(t)&=&\frac{\alpha_1-\alpha_3}{\alpha_1}
+\frac{\alpha_3\left[1+\alpha_4^2-\alpha_0^2-(\alpha_1-\alpha_3)^2\right]}
{2\alpha_1\left[1-(\alpha_1-\alpha_3)^2\right]}t+O(t^2),\\
\bar z(t)&=&\frac{-\alpha_1}{2(\alpha_1-\alpha_3)}+O(t),\quad
z(t)=\frac{-\alpha_1(\alpha_1+\alpha_2)}{\alpha_1-\alpha_3}+O(t).
\end{eqnarray}
\end{theorem} 
These solutions satisfy the system \eqref{S:1}, \eqref{S:2} and 
\eqref{y:0}, \eqref{z:0} and they are convergent since \eqref{S:1} and \eqref{S:2} are
of the Briot-Bouquet type at $t=0$ \cite{BB}. We gave the proof in section 6.
For generic values of parameters, there are no meromorphic solutions around $t=0$ 
except for these four solutions.\\

\begin{remark}
(1) These four solutions exist for the following condition:
\begin{eqnarray}
&&{\rm (0{\textrm -}I)}:\alpha_1 \ne 0, \quad\alpha_4-\alpha_0 \notin \mathbb {Z},
\quad{\rm (0{\textrm -}II)}:\alpha_1 \ne 0, \quad\alpha_4+\alpha_0 \notin \mathbb {Z},\\
&&{\rm (0{\textrm -}III)}:\alpha_1\notin \mathbb {Z},\, 
\alpha_1+\alpha_3 \notin \mathbb {Z},
\quad{\rm (0{\textrm -}IV)}:\alpha_1\notin \mathbb {Z},\, 
\alpha_1-\alpha_3 \notin \mathbb {Z}.
\end{eqnarray}
(2)\, In the case of $\alpha_0=0$, $y(t)$ of the solution (0-I) coincide with (0-II) and 
$y(t) \equiv t$.\\
\qquad In the case of $\alpha_3=0$, $y(t)$ of the solution (0-III) coincide with (0-IV) and 
$y(t)\equiv1$.\\
Both $y(t) \equiv t$ and $y(t)\equiv1$ are Riccati solutions.\\

\noindent(3) In the case of $\alpha_1=0 ~(\alpha=0)$, 
the sixth Painlev\'e equation has the following
special solution around $t=0$:
\begin{eqnarray}   
y(t)&=&t^{\pm\alpha_3}(a_0+a_1t+a_2t^2+\cdots),\quad
\bar z(t)=t^{\mp\alpha_3}(b_0+b_1t+b_2t^2+\cdots),\nonumber\\
z(t)&=&t^{\mp\alpha_3}(c_0+c_1t+c_2t^2+\cdots)\qquad(a_i, b_i,c_i \in\mathbb {C}).
\end{eqnarray}
\end{remark}
\bigskip
The B\"acklund transformations for the sixth Painlev\'e equation are defined 
in subsection 2.2.8, (where $y=q, z=p$) \cite{NY}.
If we let $\sigma_1$ and $\sigma_2$   
act on the solutions (0-I), (0-II), (0-III) and (0-IV), we then obtain the meromorphic 
solutions of the system \eqref{S:1}, \eqref{S:2} and 
\eqref{y:0}, \eqref{z:0} 
which are meromorphic around $t=1$ and $t=\infty$.\\
\\

\begin{theorem}\label{m1:1} The sixth Painlev\'e equation has the following meromorphic solutions 
around $t=1$ and $t=\infty$.\\
(1) Around $t=1$:
\begin{eqnarray}
 {\rm(1{\textrm -}I)}:\quad y(t)&=&1
+\frac{\alpha_3}{\alpha_0-\alpha_3}(1-t)\nonumber\\
&+&\frac{\alpha_0\alpha_3\left[-1-\alpha_1^2+\alpha_4^2+(\alpha_0-\alpha_3)^2\right]}
{2\left[1-(\alpha_0-\alpha_3)^2\right](\alpha_0-\alpha_3)^2}(1-t)^2+O((1-t)^3),\\
\bar z(t)&=&\frac{1-\alpha_1^2+\alpha_4^2-(\alpha_0-\alpha_3)^2}
{4\left[1-(\alpha_0-\alpha_3)^2\right]}+O((1-t)),                     \\
z(t)&=&\frac{\alpha_0-\alpha_3}{1-t}+O((1-t)^0),
\end{eqnarray}
\begin{eqnarray}
{\rm(1{\textrm -}II)}:\quad y(t)&=&1+\frac{-\alpha_3}{\alpha_0+\alpha_3}(1-t)\nonumber\\
&+&\frac{\alpha_0\alpha_3\left[1+\alpha_1^2-\alpha_4^2-(\alpha_0+\alpha_3)^2\right]}
{2\left[1-(\alpha_0+\alpha_3)^2\right](\alpha_0+\alpha_3)^2}(1-t)^2+O((1-t)^3),\\
\bar z(t)&=&\frac{1-\alpha_1^2+\alpha_4^2-(\alpha_0+\alpha_3)^2}
{4\left[1-(\alpha_0+\alpha_3)^2\right]}+O((1-t)),                     \\
z(t)&=&\frac{\alpha_2(\alpha_1+\alpha_2)}{\alpha_0+\alpha_3-1}+O((1-t)),  \\
{\rm(1{\textrm -}III)}:\quad y(t)&=&-\frac{\alpha_4}{\alpha_1}
+\frac{\alpha_4\left[1+\alpha_3^2-\alpha_0^2-(\alpha_4+\alpha_1)^2\right]}
{2\alpha_1\left[1-(\alpha_4+\alpha_1)^2\right]}(1-t)+O((1-t)^2),\\
\bar z(t)&=&\frac{\alpha_1}{2(\alpha_1+\alpha_4)}+O((1-t)),\quad
 z(t)=\frac{\alpha_1\alpha_2}{\alpha_1+\alpha_4}+O((1-t)),  
\end{eqnarray}
\begin{eqnarray}
{\rm(1{\textrm -}IV)}:\quad y(t)&=&\frac{\alpha_4}{\alpha_1}
+\frac{-\alpha_4\left[1+\alpha_3^2-\alpha_0^2-(\alpha_4-\alpha_1)^2\right]}
{2\alpha_1\left[1-(\alpha_4-\alpha_1)^2\right]}(1-t)+O\left((1-t)^2\right),\\
\bar z(t)&=&\frac{-\alpha_1}{2(\alpha_4-\alpha_1)}+O((1-t)),\
z(t)=\frac{\alpha_1(\alpha_2+\alpha_4)}
{\alpha_1-\alpha_4}+O((1-t)).
\end{eqnarray}
(2) Around $t=\infty$:
\begin{eqnarray}
{\rm(\infty{\textrm -}I)}:\quad y(t)&=&\frac{\alpha_1-\alpha_0}{\alpha_1}t
+\frac{\alpha_1\left[\left(1+\alpha_3^2+\alpha_4^2-\alpha_1^2-(\alpha_0-\alpha_1)^2\right)
\left(\frac{\alpha_0}{\alpha_1}\right)^2+1+\alpha_0^2\right]}
{2\alpha_0\left[1-(\alpha_0-\alpha_1)^2\right]} \nonumber\\
&+&O((t^{-1})),\\
\bar z(t)&=&\frac{-\alpha_1}{2(\alpha_1-\alpha_0)}\frac{1}{t}+O(t^{-2}),\quad
  z(t)=-\frac{\alpha_1(\alpha_1+\alpha_2)}{(\alpha_1-\alpha_0)t}+O(t^{-2}),\\
{\rm(\infty{\textrm -}II)}:\quad y(t)&=&\frac{\alpha_1+\alpha_0}{\alpha_1}t
+\frac{-\alpha_1\left[\left(1+\alpha_3^2+\alpha_4^2-\alpha_1^2-(\alpha_0+\alpha_1)^2\right)
\left(\frac{\alpha_0}{\alpha_1}\right)^2+1+\alpha_0^2\right]}
{2\alpha_0\left[1-(\alpha_0+\alpha_1)^2\right]} \nonumber\\   
&+&O((t^{-1}),\\
 \bar z(t)&=&\frac{-\alpha_1}{2(\alpha_1+\alpha_0)}\frac{1}{t}+O(t^{-2}),\quad
  z(t)=-\frac{\alpha_1\alpha_2}{\alpha_1+\alpha_0}\cdot\frac{1}{t}+O(t^{-2}),
\end{eqnarray}
\begin{eqnarray}
{\rm(\infty{\textrm -}III)}:\quad y(t)&=&\frac{\alpha_4}{\alpha_4+\alpha_3}
+\frac{-\alpha_3\alpha_4\left[-1+\alpha_0^2-\alpha_1^2+(\alpha_3+\alpha_4)^2\right]}
{2\left[1-(\alpha_3+\alpha_4)^2\right](\alpha_3+\alpha_4)^2}\frac{1}{t}+O(t^{-2}),\\    
\bar z(t)&=&\frac{1-\alpha_1^2+\alpha_0^2-(\alpha_3+\alpha_4)^2}
{4\left[1-(\alpha_3+\alpha_4)^2\right]}\frac{1}{t}+O(t^{-2}),                     \\  
z(t)&=&\frac{\alpha_2(\alpha_1+\alpha_2)}{1-\alpha_3-\alpha_4}\cdot\frac{1}{t}+O(t^{-2}),\\
{\rm(\infty{\textrm -}IV)}:\quad y(t)&=&\frac{\alpha_4}{\alpha_4-\alpha_3}
+\frac{\alpha_3\alpha_4\left[-1+\alpha_0^2-\alpha_1^2+(\alpha_3-\alpha_4)^2\right]}
{2\left[1-(\alpha_3-\alpha_4)^2\right](\alpha_3-\alpha_4)^2}\frac{1}{t}+O(t^{-2}),\\     
\bar z(t)&=&\frac{1-\alpha_1^2+\alpha_0^2-(\alpha_3-\alpha_4)^2}
{4\left[1-(\alpha_3-\alpha_4)^2\right]}\frac{1}{t}+O(t^{-2}),\\  
z(t)&=&\alpha_4-\alpha_3+O(t^{-1}).
\end{eqnarray}
\end{theorem}
%
\begin{remark}
 If we assume the meromorphy of a solution 
around $t=0$ and $t=1$, $y(t)$ and $\bar z(t)$ inevitably
become holomorphic there.
\end{remark}
\begin {theorem} 
These twelve meromorphic solutions are invariant under the action of 
the B\"acklund transformation group.
\end{theorem}

\begin{figure}[htbp]
        \begin{picture}(300,100)(-90,0)
        \put(52,95){$s_1$}
        \put(134,95){$s_2$}
        \put(216,95){$s_3$}
        \put(52,50){$s_4$}
        \put(134,50){$s_2$}
        \put(216,50){$s_3$}
        \put(52,5){$s_3$}
        \put(134,5){$s_2$}
        \put(216,5){$s_4$}
        \put(15,22){$\sigma_1$}
        \put(15,70){$\sigma_2$}
        \put(95,22){$\sigma_1$}
        \put(95,70){$\sigma_2$}
        \put(175,22){$\sigma_1$}
        \put(175,70){$\sigma_2$}
        \put(255,22){$\sigma_1$}
        \put(255,70){$\sigma_2$}
        \put(0,90){\rm($\infty${\textrm -}I)}
        \put(80,90){\rm($\infty${\textrm -}II)}
        \put(160,90){\rm($\infty${\textrm -}III)}
        \put(240,90){\rm($\infty${\textrm -}IV)}
        \put(3,45){\rm($0${\textrm -}I)}
        \put(83,45){\rm($0${\textrm -}II)}
        \put(163,45){\rm($0${\textrm -}III)}
        \put(243,45){\rm($0${\textrm -}IV)}
        \put(3,0){\rm($1${\textrm -}I)}
        \put(83,0){\rm($1${\textrm -}II)}
        \put(163,0){\rm($1${\textrm -}III)}
        \put(243,0){\rm($1${\textrm -}IV)}
        \put(56,92){\vector(1,0){18}}
        \put(56,92){\vector(-1,0) {18}}
        \put(138,92){\vector(1,0){18}}
        \put(138,92){\vector(-1,0) {18}}
        \put(220,92){\vector(1,0){18}}
        \put(220,92){\vector(-1,0) {18}}
        \put(56,47){\vector(1,0){18}}
        \put(56,47){\vector(-1,0) {18}}
        \put(138,47){\vector(1,0){18}}
        \put(138,47){\vector(-1,0) {18}}
        \put(220,47){\vector(1,0){18}}
        \put(220,47){\vector(-1,0) {18}}
        \put(56,2){\vector(1,0){18}}
        \put(56,2){\vector(-1,0) {18}}
        \put(138,2){\vector(1,0){18}}
        \put(138,2){\vector(-1,0) {18}}
        \put(220,2){\vector(1,0){18}}
        \put(220,2){\vector(-1,0) {18}}
        \put(13,24){\vector(0,1){11}}
        \put(13,24){\vector(0,-1) {11}}
        \put(13,72){\vector(0,1){11}}
        \put(13,72){\vector(0,-1) {11}}
        \put(93,24){\vector(0,1){11}}
        \put(93,24){\vector(0,-1) {11}}
        \put(93,72){\vector(0,1){11}}
        \put(93,72){\vector(0,-1) {11}}
        \put(173,24){\vector(0,1){11}}
        \put(173,24){\vector(0,-1) {11}}
        \put(173,72){\vector(0,1){11}}
        \put(173,72){\vector(0,-1) {11}}
        \put(253,24){\vector(0,1){11}}
        \put(253,24){\vector(0,-1) {11}}
        \put(253,72){\vector(0,1){11}}
        \put(253,72){\vector(0,-1) {11}}
        \end{picture}
\caption{The B\"acklund transformations of the twelve solutions}
\end{figure}
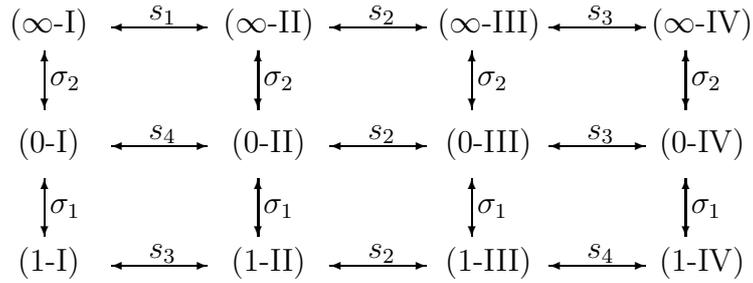

%
\subsection{The linear monodromy for the solution {\rm(0-I)}}
For a solution of the sixth Painlev\'e equation, let $M_j (j=0,t,1,\infty)$ be the 
monodromy matrices of the equation \eqref{J:3} along the path around $x=j$ shown in 
Figure 2.\\
\begin{figure}[htbp]
   \centering
    \includegraphics[scale=0.7]{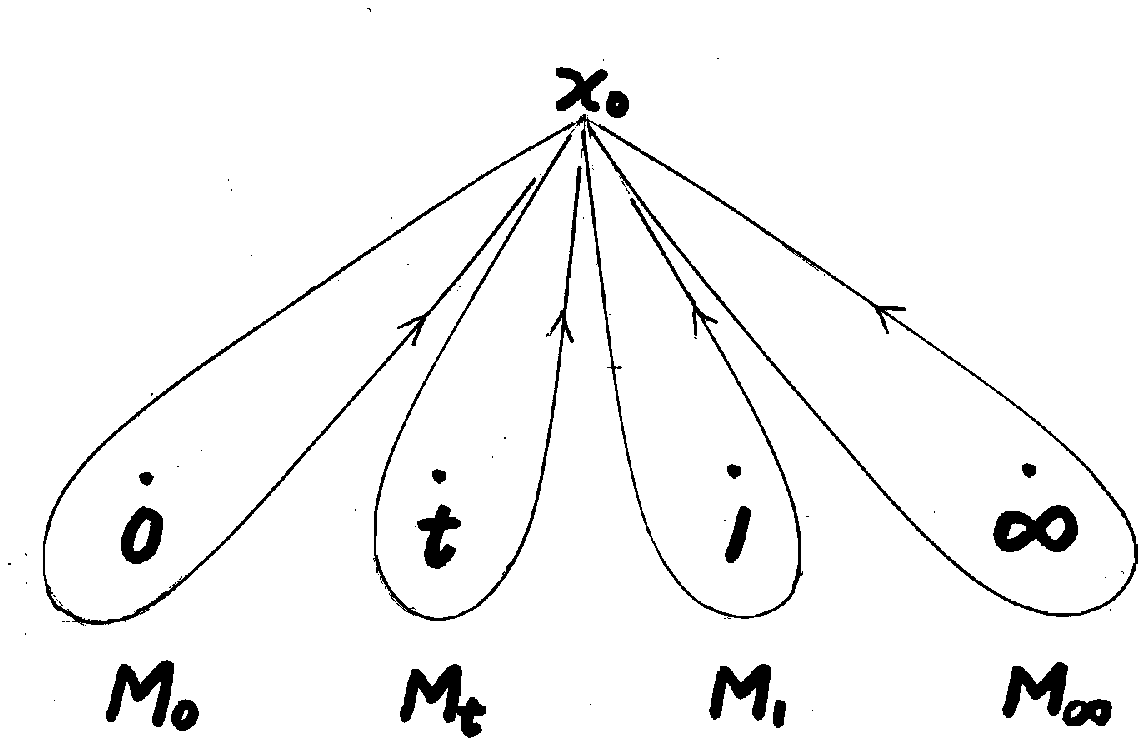}
   \caption{The paths going around regular singular points with the base point $x_0$.}
   \end{figure}

\noindent Note that $M_j  (j=0,t,1,\infty)$ satisfy
\begin{eqnarray}
     M_{\infty}M_1M_tM_0=I_2.
\end{eqnarray}
We can then calculate the linear monodromy $\{M_0,M_t,M_1,M_{\infty}\}$ 
explicitly for the solution (0-I) by the method given in \cite{MJ2}.\\
%

\subsection {The limit of \eqref{J:3}}
 We will take 
the limit $ t\rightarrow 0 $ after substituting the solution (0-I) into \eqref{J:3}.\\ 
{\bf5.3.1)}  Put $\bar Y={}^t(\bar \psi^{(1)},\bar \psi^{(2)} )$, then the limit
\begin{eqnarray}
\bar \psi_0^{(1)}(x)=\lim_{t\rightarrow 0}\bar\psi^{(1)}(x,t)
\end{eqnarray}
satisfies
\begin{eqnarray}
\frac{d^2 \bar\psi_0^{(1)}}{dx^2}
&+&\left(\frac{1}{x}+\frac{1}{x-1}\right)\frac{d \bar\psi_0^{(1)}}{dx}\nonumber\\
&-&\biggl[\frac{(\alpha_0-\alpha_4)^2}{4}\cdot\frac{1}{x^2}+\frac{\alpha_3^2}{4}\frac{1}
{(x-1)^2}-\frac{1-\alpha_1^2+\alpha_3^2+(\alpha_0-\alpha_4)^2}{4x(x-1)}
\biggr]\bar\psi_0^{(1)}=0.\nonumber\\\label{M:1}
\end{eqnarray}
The Riemann scheme of \eqref{M:1} is
\begin{eqnarray} 
&&P\left(
\begin{array}{cccc}
  x=0\cdot t &  x=1  & x= \infty & \\
\begin{matrix} \frac{\alpha_0-\alpha_4}{2} \\-\frac{\alpha_0- \alpha_4}{2} \end{matrix} &
   \begin{matrix} -\frac{\alpha_3}{2} \\ \frac{\alpha_3}{2} \end{matrix} &
   \begin{matrix} \frac{1}{2}(1+\alpha_1)  \\\frac{1}{2}(1-\alpha_1)      \end{matrix} &
   \begin{matrix} ;x \\   \end{matrix} \end{array} \right)  \\
&&=x^{\frac{\alpha_0-\alpha_4}{2}}(x-1)^{-\frac{\alpha_3}{2}}P\left(
\begin{array}{cccc}
  x=0\cdot t &  x=1  & x= \infty  \\
\begin{matrix} 0 \\ \alpha_4- \alpha_0 \end{matrix} &
   \begin{matrix} 0 \\ \alpha_3 \end{matrix} &
   \begin{matrix}\alpha_0+\alpha_1+\alpha_2  \\ \alpha_0+\alpha_2 \end{matrix} &
   \begin{matrix} ;x \\   \end{matrix} \end{array} \right).
\end{eqnarray} 
Therefore a fundamental system of solutions 
of \eqref{M:1} is
\begin{eqnarray}
\biggl(x^{\frac{\alpha_0-\alpha_4}{2}}(x-1)^{-\frac{\alpha_3}{2}}
{}_2F_1(\alpha_0+\alpha_1+\alpha_2,\alpha_0+\alpha_2,1+\alpha_0-\alpha_4;x),\nonumber\\
x^{\frac{\alpha_4-\alpha_0}{2}}(x-1)^{-\frac{\alpha_3}{2}}
{}_2F_1(\alpha_1+\alpha_2+\alpha_4,\alpha_2+\alpha_4,1+\alpha_4-\alpha_0;x)\biggr).
\end{eqnarray}
 \begin{figure}[htbp]
   \centering
    \includegraphics[scale=0.7]{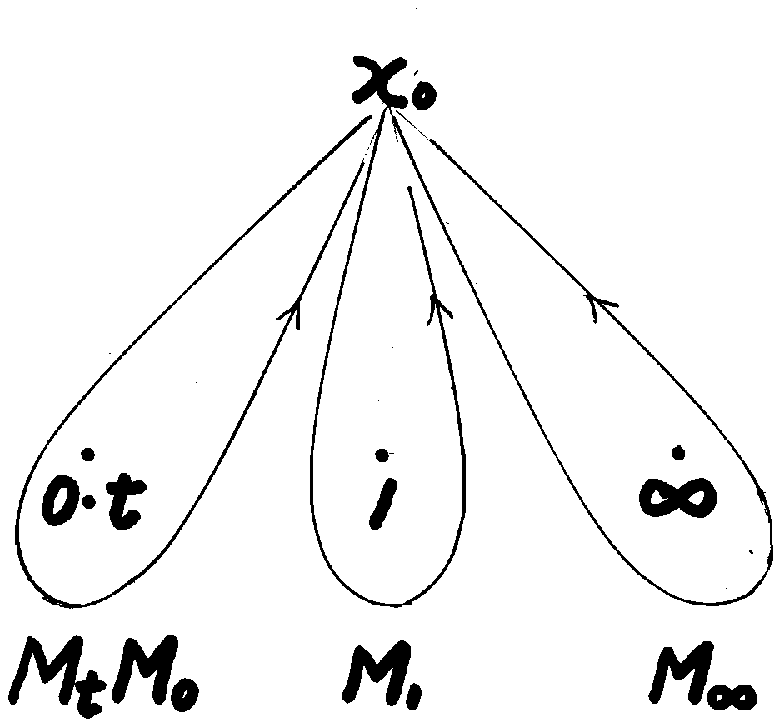}
   \caption{The paths used to calculate the linear monodromy of \eqref{M:1}}
\end{figure}
The linear monodromy of \eqref{M:1} is equivalent to $\{M_tM_0,M_1,M_{\infty}\}$.\\
The exponent matrices of \eqref{M:1} at $x=0, 1$ and $\infty$ are given by
\begin{eqnarray}
   T_0&=&\left(\begin{array}{cc} \frac{\alpha_0-\alpha_4}{2} & 0 \\
              0 & -\frac{\alpha_0-\alpha_4}{2}  \end{array}        \right),\quad
   T_1=\left(\begin{array}{cc} -\frac{\alpha_3}{2} & 0 \\
              0 & \frac{\alpha_3}{2}  \end{array}        \right),\\
T_{\infty}&=&\left(\begin{array}{cc} \frac{1+\alpha_1}{2} & 0 \\
              0  & \frac{1-\alpha_1}{2}  \end{array}        \right).
\end{eqnarray}
We may assume
\begin{eqnarray}
M_tM_0&=&e^{2\pi iT_0}=\left(\begin{array}{cc} e^{\pi i (\alpha_0-\alpha_4)} &  0  \\
               0  & e^{-\pi i (\alpha_0-\alpha_4)}\end{array}         \right)  \label{28:0},\\
M_1&=&\Gamma_{01}^{-1}e^{2 \pi i T_1}\Gamma_{01},\quad
M_{\infty}=\Gamma_{0\infty}^{-1}e^{2 \pi i T_{\infty}}\Gamma_{0\infty}, \label{7:0}
\end{eqnarray}
where
\begin{eqnarray}
 \Gamma_{01}&=&\biggl( \begin{array}{cc}
    {\Gamma(1+\alpha_0-\alpha_4)\Gamma(\alpha_3)}\over
         {\Gamma(1-\alpha_1-\alpha_2-\alpha_4)\Gamma(1-\alpha_2-\alpha_4)} &
   \frac {\Gamma(1+\alpha_4-\alpha_0)\Gamma(\alpha_3)}
         {\Gamma(1-\alpha_0-\alpha_1-\alpha_2)\Gamma(1-\alpha_0-\alpha_2)}\\
   \frac {\Gamma(1+\alpha_0-\alpha_4)\Gamma(-\alpha_3)}
         {\Gamma(\alpha_0+\alpha_1+\alpha_2)\Gamma(\alpha_0+\alpha_2)}&
   \frac {\Gamma(1+\alpha_4-\alpha_0)\Gamma(-\alpha_3)}
         {\Gamma(\alpha_1+\alpha_2+\alpha_4)\Gamma(\alpha_2+\alpha_4)}
                  \end{array}  \biggr), \label{g:1}\\
\Gamma_{0\infty}&=&\left( \begin{array}{cc}
   e^{(\alpha_0+\alpha_1+\alpha_2)\pi i} {\Gamma(1+\alpha_0-\alpha_4)\Gamma(-\alpha_1)}\over
         {\Gamma(\alpha_0+\alpha_2)\Gamma(1-\alpha_1-\alpha_2-\alpha_4)} &
  e^{(\alpha_1+\alpha_2+\alpha_4)\pi i} {\Gamma(1+\alpha_4-\alpha_0)\Gamma(-\alpha_1)}\over
         {\Gamma(\alpha_2+\alpha_4)\Gamma(1-\alpha_0-\alpha_1-\alpha_2)}\\
  e^{(\alpha_0+\alpha_2)\pi i} {\Gamma(1+\alpha_0-\alpha_4)\Gamma(\alpha_1)}\over
         {\Gamma(\alpha_0+\alpha_1+\alpha_2)\Gamma(1-\alpha_2-\alpha_4)}&
  e^{(\alpha_2+\alpha_4)\pi i} {\Gamma(1+\alpha_4-\alpha_0)\Gamma(\alpha_1)}\over
         {\Gamma(\alpha_1+\alpha_2+\alpha_4)\Gamma(1-\alpha_0-\alpha_2)}
                  \end{array}  \right).\label{g:2}
\end{eqnarray}

We should separate the monodromy data $M_tM_0$.\\
\noindent{\bf5.3.2)} In the following, we consider the confluence of $x=1$ and $x=\infty$
 in \eqref{J:3}.\\ 
Put $x=t\xi$, then the limit 
\begin{eqnarray}
\bar  \psi_1^{(1)}(\xi)=\lim_{t\rightarrow 0}\bar  \psi^{(1)}(t\xi,t)
\end{eqnarray}
satisfies
\begin{eqnarray}
\frac{d^2 \bar\psi_1^{(1)}}{d\xi^2}&+&\biggl(\frac{1}{\xi}+ \frac{1}{\xi-1}
-\frac{1}{\xi-s}\biggr)\frac{d \bar\psi_1^{(1)}}{d \xi}\nonumber\\
&-&\biggl[\frac{\alpha_4^2}{4}\frac{1}{\xi^2}+\frac{\alpha_0^2}{4}\frac{1}{(\xi-1)^2}
+\frac{-\alpha_4^2-\alpha_0^2+(\alpha_0-\alpha_4)^2}{4\xi(\xi-1)}\biggr]\bar\psi_1^{(1)}
=0,
\label{M:3}
\end{eqnarray}
where 
\begin{eqnarray}
s=\frac{\alpha_4}{\alpha_4-\alpha_0}.
\end{eqnarray}
This is a Heun's type equation with an apparent singularity at $\xi=s=\alpha_4/(\alpha_4-\alpha_0)$. 
The singularities $\xi=0,1$ and $\infty$ correspond to $x=0,t$ and $1\cdot\infty$, respectively.
The Riemann scheme of \eqref{M:3} is
\begin{equation} 
P\ \left(
\begin{array}{ccccc}
  \xi=0  & \xi=1&  \xi=s & \xi=\infty &    \\
   \begin{matrix}  -\frac{\alpha_4}{2}  \\ \frac{\alpha_4}{2} \end{matrix} &
   \begin{matrix}\frac{\alpha_0}{2} \\-\frac{\alpha_0}{2}  \end{matrix} &
   \begin{matrix} 0 \\ 2 \end{matrix}& 
   \begin{matrix} \frac{\alpha_4-\alpha_0}{2} \\ -\frac{\alpha_4-\alpha_0}{2}     \end{matrix}& 
   \begin{matrix} ;\xi \\  \end{matrix} 
\end{array} \right). 
\end{equation} 
A fundamental system of solutions of \eqref{M:3} is 
\begin{eqnarray}
 \biggl(\begin{array}{cc} \xi^{-\frac{\alpha_4}{2}}(\xi-1)^{\frac{\alpha_0}{2}}
,&
  \xi^{\frac{\alpha_4}{2}}(\xi-1)^{-\frac{\alpha_0}{2}}
 \end{array} \biggr).                                                     
\end{eqnarray}
    \begin{figure}[htbp]
   \centering
    \includegraphics[scale=0.7]{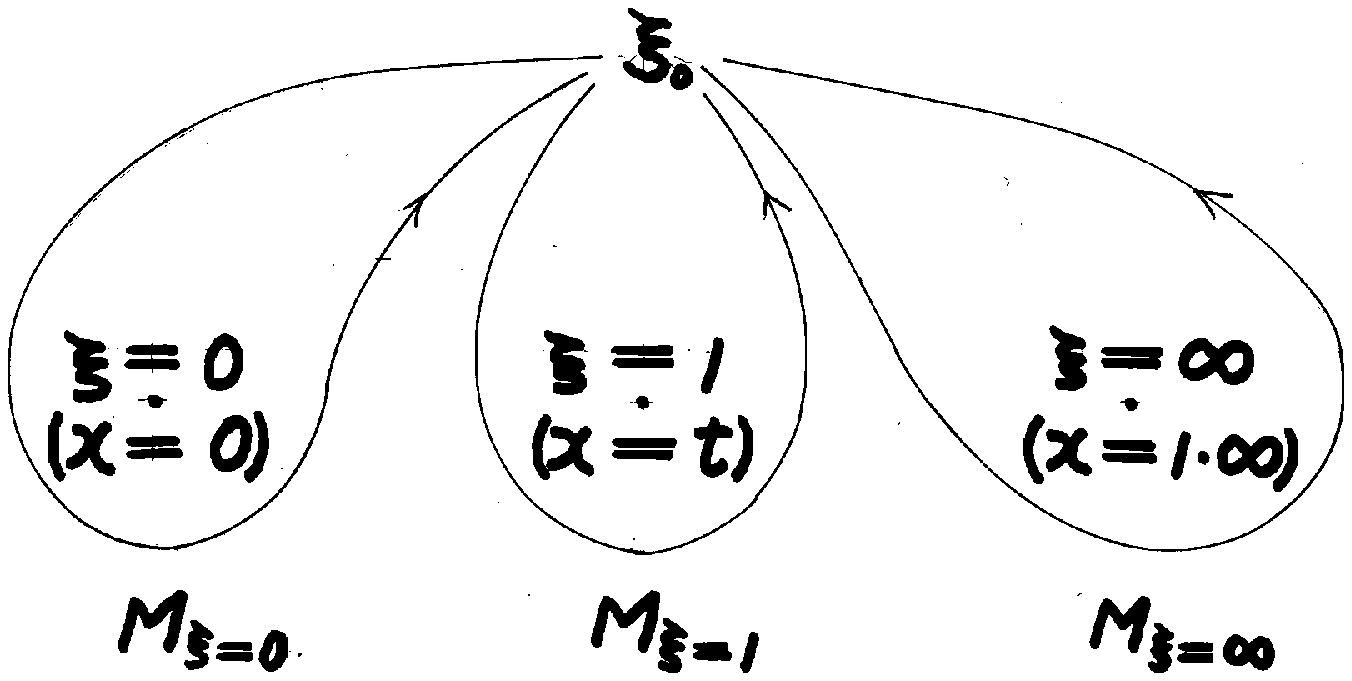}
   \caption{The paths use to calculate the linear monodromy of \eqref{M:3}}
   \end{figure}
The linear monodromy  
$\{L_0,L_1,L_{\infty}\}$ of \eqref{M:3} is equivalent to $\{M_0,M_t,M_{\infty}M_1\}$.\\
\begin{eqnarray}
M_0&=&P^{-1}L_0P,\quad M_t=P^{-1}L_1P,\quad M_{\infty}M_1=P^{-1}L_{\infty}P\label{28:1}
\end{eqnarray}
for a matrix $P\in SL(2,\mathbb {C})$.\\
The linear monodromy$\{L_0,L_1,L_{\infty}\}$ is given by
\begin{eqnarray} 
L_0&=&\left( \begin{array}{cc}  e^{-\pi i\alpha_4} &  0 \\
                 0 &  e^{\pi i\alpha_4}\end{array}  \right),\quad
L_1=\left( \begin{array}{cc}  e^{\pi i\alpha_0} &  0 \\
                 0 &  e^{-\pi i\alpha_0}\end{array}  \right),\label{28:2}\\
L_{\infty}&=&\left( \begin{array}{cc}  e^{-\pi i(\alpha_0-\alpha_4)} &  0 \\
                 0 &  e^{\pi i(\alpha_0-\alpha_4)}\end{array}  \right).\label{0:2}
\end{eqnarray}
Comparing \eqref{28:0} and \eqref{28:1}, \eqref{28:2}, we have
\begin{eqnarray} 
     M_{t}M_{0}=P^{-1}L_1L_0P,\quad M_tM_0=L_1L_0=\left(\begin{array}{cc} e^{\pi i (\alpha_0-\alpha_4)} &  0  \\
               0  & e^{-\pi i (\alpha_0-\alpha_4)}\end{array}         \right).
\end{eqnarray}
Therefore $P$ is a diagonal matrix, since $\alpha_0-\alpha_4 \notin \mathbb {Z}$ 
for the solution (0-I).\\


\begin{theorem} 
 The linear monodromy of \eqref{J:3} for the solution (0-I) is as follows:
\begin{eqnarray}
M_{0}&=&\left( \begin{array}{cc}  e^{-\pi i\alpha_4} &  0 \\
                 0 &  e^{\pi i\alpha_4}\end{array}  \right),\quad
M_{t}=\left( \begin{array}{cc}  e^{\pi i\alpha_0} &  0 \\
                 0 &  e^{-\pi i\alpha_0}\end{array}  \right),\label{0:4}\\
M_1&=&\Gamma_{01}^{-1}\left(\begin{array}{cc} e^{-\pi i \alpha_3} & 0 \\
0 & e^{\pi i\alpha_3} \end{array} \right)\Gamma_{01},\quad
M_{\infty}=\Gamma_{0\infty}^{-1}\left(\begin{array}{cc} -e^{\pi i\alpha_1} & 0 \\
0 & -e^{-\pi i\alpha_1} \end{array} \right)\Gamma_{0\infty}.\label{0:5}
\end{eqnarray}
where $\Gamma_{01}$ and $\Gamma_{0\infty}$ are given in \eqref{g:1} and \eqref{g:2}.
 We remark that $\alpha_0-\alpha_4 \notin \mathbb{Z}$ if the solution (0-I) 
exists.
\end{theorem}
\noindent In a similar way, we can calculate the linear monodromy
explicitly for all of the twelve solutions in Theorem \ref{m0:1} and Theorem \ref{m1:1}\\
\\
\begin{theorem} 
 The twelve solutions in Theorem \ref{m0:1} and Theorem \ref{m1:1}  
are all monodromy solvable.
\end{theorem}

 \subsection{Comparison with classical solutions}

Umemura studied special solutions of the Painlev\'e equations \cite{HU1}. Umemura's classical 
solutions are either rational solutions or the Riccati solutions. 
We show that some of our solutions include an algebraic solution
and one of the Riccati solutions.

\noindent1) In the case of $\alpha_1=\alpha_4$ and $\alpha_0=\alpha_3 ~(\alpha+\beta=0,\quad \gamma+\delta=\frac{1}{2})$,
the sixth Painlev\'e equation has an algebraic solution 
\begin{eqnarray}
y(t)&=&
\sqrt t=1+\frac{1}{2}(t-1)+\frac{1}{2!}\cdot\frac{-1}{4}(t-1)^2+\cdots,\nonumber\\
z(t)&=&\frac{1}{4}(2\alpha_3+2\alpha_4-1)
\frac {1}{\sqrt t}\label{5:1a}\\
&=&\frac{1}{4}(2\alpha_3+2\alpha_4-1)\biggl[1-\frac{1}{2}(t-1)+\frac{1}{2!}\cdot\frac{3}{4}(t-1)^2-\cdots\biggr].\nonumber
\end{eqnarray}
The solution \eqref{5:1a} is a special case of the solution \rm{(1-II)} for
 $\alpha_1=\alpha_4,\alpha_0=\alpha_3$.

\noindent 2) In the case of $\alpha_0=0~(\delta=\frac{1}{2})$, the sixth Painlev\'e equation has the Riccati solution
\begin{eqnarray}
   (1)\quad y(t)&=&t, \nonumber\\ z(t)&=&\frac{{}_2F_1'(\alpha_2,\alpha_1+\alpha_2,1-\alpha_4;t)}
{{}_2F_1(\alpha_2,\alpha_1+\alpha_2,1-\alpha_4;t)}
=\frac{\alpha_2(\alpha_1+\alpha_2)}{1-\alpha_4}+O(t),\\
   (2)\quad y(t)&=&t, \nonumber\\ z(t)&=&\frac{\left[t^{\alpha_4}{}_2F_1(
\alpha_2+\alpha_4,\alpha_1+\alpha_2+\alpha_4,1+\alpha_4;t)\right]'}
{t^{\alpha_4}{}_2F_1(
\alpha_2+\alpha_4,\alpha_1+\alpha_2+\alpha_4,1+\alpha_4;t)}=\frac{\alpha_4}{t}\left(1+O(t)\right).
\end{eqnarray} 
These are obtained by putting $\alpha_0=0$ in the solution \rm{(0-II)} and \rm{(0-I)}, 
respectively.\\
\noindent 3) In the case of $\alpha_2=0$, the system \eqref{S:1} and \eqref{S:2} has the Riccati solution
\begin{eqnarray}
(1)\quad z(t)&\equiv&0, \nonumber \\
y(t)&=&-\frac{t(t-1)}{\alpha_1}\frac{\left
[(t-1)^{\alpha_4}{}_2F_1(
\alpha_4,1-\alpha_3,\alpha_0+\alpha_4;t)\right]'}
{(t-1)^{\alpha_4}{}_2F_1(
\alpha_4,1-\alpha_3,\alpha_0+\alpha_4;t)}\nonumber\\
&=&\frac{\alpha_4}{\alpha_0+\alpha_4}t+O(t^2),\\
(2)\quad z(t)&\equiv&0, \nonumber \\
y(t)&=&-\frac{t(t-1)}{\alpha_1}\frac{\left
[t^{\alpha_1+\alpha_3}(t-1)^{\alpha_4}{}_2F_1(
1-\alpha_0,1+\alpha_0,1+\alpha_1+\alpha_3;t)\right]'}
{t^{\alpha_1+\alpha_3}(t-1)^{\alpha_4}{}_2F_1
(1-\alpha_0,1+\alpha_0,1+\alpha_1+\alpha_3;t)}\nonumber\\
&=&\frac{\alpha_1+\alpha_3}{\alpha_1}+O(t).
\end{eqnarray}
These are obtained by putting $\alpha_2=0$ in the solution \rm{(0-II)} and \rm{(0-III)}, 
respectively.

\noindent 4) In the case of $\alpha_3=0\,(\gamma=0)$, the system \eqref{S:1} 
and \eqref{S:2} has the Riccati solution
\begin{eqnarray}
(1)\quad y(t)&\equiv&1, \nonumber \\
z(t)&=&-t\frac{\left[t^{\alpha_2}{}_2F_1(
\alpha_2,\alpha_2+\alpha_4,1-\alpha_1;t)\right]'}
{ t^{\alpha_2}{}_2F_1(
\alpha_2,\alpha_2+\alpha_4,1-\alpha_1;t)}
=-\alpha_2+O(t),\\
(2)\quad y(t)&\equiv&1, \nonumber \\
z(t)&=&-t\frac{\left[t^{\alpha_1+\alpha_2}{}_2F_1(
\alpha_1+\alpha_2,\alpha_1+\alpha_2+\alpha_4,1+\alpha_1;t)\right]'}
{ t^{\alpha_1+\alpha_2}{}_2F_1(
\alpha_1+\alpha_2,\alpha_1+\alpha_2+\alpha_4,1+\alpha_1;t)}\nonumber\\
&=&-(\alpha_1+\alpha_2)+O(t).
   \end{eqnarray}
These are obtained by putting $\alpha_3=0$ in the solution \rm{(0-III)} and \rm{(0-IV)}, 
respectively.\\
\noindent 5) In the case of $\alpha_4=0~(\beta=0)$, the system \eqref{S:1} and \eqref{S:2}
 has the Riccati solution
\begin{eqnarray}
(1)\quad y(t)&\equiv&0, \nonumber \\
     z(t)&=&(t-1)\frac{[(t-1)^{\alpha_2}{}_2F_1(\alpha_2,\alpha_2+\alpha_3,1-\alpha_0;t)]'}
                    {(t-1)^{\alpha_2}{}_2F_1(\alpha_2,\alpha_2+\alpha_3,1-\alpha_0;t)}
\nonumber\\
&=&\frac{\alpha_2(\alpha_1+\alpha_2)}{1-\alpha_0}+O(t).\\
(2)\quad y(t)&\equiv&0, \nonumber \\
     z(t)&=&(t-1)\frac{[t^{\alpha_0}(t-1)^{\alpha_2}{}_2F_1
(\alpha_0+\alpha_2,\alpha_0+\alpha_2+\alpha_3,1+\alpha_0;t)]'}
                    {t^{\alpha_0}(t-1)^{\alpha_2}{}_2F_1
(\alpha_0+\alpha_2,\alpha_0+\alpha_2+\alpha_3,1+\alpha_0;t)}\nonumber\\
&=&-\frac{\alpha_0}{t}+O(t^0).
\end{eqnarray}
These are obtained by putting $\alpha_4=0$ in the solution \rm{(0-II)} and \rm{(0-I)}, 
respectively.\\

\section{ The Briot-Bouquet theorem}\label{Briot-Bouquet}
The Briot-Bouquet theorem \cite{BB} is well-known but 
we will explain the Briot-Bouquet theorem for a system and give a brief proof. 
From the Briot-Bouquet theorem   series expansions of solutions 
of the fifth and the sixth Painlev\'e equations in Theorem \ref{p5:mero}, \ref{m0:1}
and  \ref{m1:1} converge around the fixed singularities.
 We denote $(f)_0:=f(0)$ for a holomorphic function $f$.
 
\begin{theorem} (Briot-Bouquet)
For the simultaneous equations
\begin{eqnarray}
 t\frac{du}{dt}=f(u,v,t),\label{7:1} \\
 t\frac{dv}{dt}=g(u,v,t),\label{7:2}
\end{eqnarray}
where f and g are holomorphic functions of $u,v$ and $t$ near the origin. 
We assume that $f(0,0,0)=0,\ g(0,0,0)=0$. Then a 
holomorphic solution with the initial condition $u(0)=0, v(0)=0$ exists  if
\begin{eqnarray*}
\Delta_n:=\begin{vmatrix}
n-(f_u)_0 &(-f_v)_0 \\ (-g_u)_0 & n-(g_v)_0\end{vmatrix}\ne0, 
\end{eqnarray*}
for any non-negative integer $n$.
\end{theorem}

\par\noindent
{\it Proof.}\quad At first we will show the existence of a formal solution $u, v$ for \eqref{7:1} 
and \eqref{7:2}. 
Expand  $u, v$ as 
\begin{equation*}\label{bb:1}
u= a_1t+a_2t^2+\cdots,  \quad v = b_1t+b_2t^2+\cdots.
\end{equation*}
From \eqref{7:1} and \eqref{7:2} we have
\begin{eqnarray}
 t\frac{d^2u}{dt^2}+\frac{du}{dt}=f_t+f_u\frac{du}{dt}+f_v\frac{dv}{dt},\label{7:6}\\
 t\frac{d^2v}{dt^2}+\frac{dv}{dt}=g_t+g_u\frac{du}{dt}+g_v\frac{dv}{dt}.\label{7:7}
\end{eqnarray}
Putting $t=0$, we have
\begin{eqnarray*}
\left[1-(f_u)_0\right]\left(\frac{du}{dt}\right)_0-(f_v)_0\left(\frac{dv}{dt}\right)_0=(f_t)_0,\\
-(g_u)_0\left(\frac{du}{dt}\right)_0+\left[1-(g_v)_0\right]\left(\frac{dv}{dt}\right)_0=(g_t)_0.
\end{eqnarray*}
Solving the above equations, we obtain
\begin{eqnarray*}
a_1=\left(\frac{du}{dt}\right)_0=\frac{1}{\Delta_1}\begin{vmatrix}
(f_t)_0 &(-f_v)_0 \\ (g_t)_0 & 1-(g_v)_0\end{vmatrix} ,\quad
b_1=\left(\frac{dv}{dt}\right)_0=\frac{1}{\Delta_1}\begin{vmatrix}
1-(f_u)_0 &(f_t)_0 \\ (-g_u)_0 & (g_t)_0\end{vmatrix}.
\end{eqnarray*}
Here $\Delta_1 \not=0$ from the assumption.

By differentiating \eqref{7:6} and \eqref{7:7} with $t$ and putting $t=0$, $(a_2, b_2)$ is also determined uniquely as follows:
\begin{eqnarray*}
a_2=\frac{1}{2!}\left(\frac{d^2u}{dt^2}\right)_0=\frac{1}{2!}\frac{1}{\Delta_2}\begin{vmatrix}
A_1 &(-f_v)_0 \\ B_1 & 2-(g_u)_0-(g_v)_0\end{vmatrix} ,\\
b_2=\frac{1}{2!}\left(\frac{d^2v}{dt^2}\right)_0=\frac{1}{2!}\frac{1}{\Delta_2}\begin{vmatrix}
2-(f_u)_0-(f_v)_0    &A_1 \\ (-g_u)_0 & B_1\end{vmatrix} ,
\end{eqnarray*}
where 
\begin{eqnarray*}
A_1&=&\biggl[f_{tt}+2f_{tu}\frac{du}{dt}+f_{uu}(\frac{du}{dt})^2+2f_{tv}\frac{dv}{dt}+2f_{uv}\frac{du}{dt}
\frac{dv}{dt} 
+f_{vv}(\frac{dv}{dt})^2\biggr]_{t=0},\\
B_1&=&\biggl[g_{tt}+2g_{tu}\frac{du}{dt}+g_{uu}(\frac{du}{dt})^2+2g_{tv}\frac{dv}{dt}+2g_{uv}\frac{du}{dt}
\frac{dv}{dt} 
+g_{vv}(\frac{dv}{dt})^2\biggr]_{t=0},
\end{eqnarray*}
and so on. Thus coefficients $(a_n, b_n)$ can be uniquely determined.

In the second step, we will show the formal solutions $u=\sum_{k=1}^{\infty}a_kt^k$ 
and $v=\sum_{k=1}^{\infty}b_kt^k$  are convergent. 
We prepare the following  auxiliary  functions $p(t), q(t)$ defined as implicit functions:
\begin{eqnarray*}
p=f(t,p,q),\quad q=g(t,p,q).
\end{eqnarray*}
Since $\Delta_1\not= 0$,  the holomorphic functions $p$ and $q$
with  $p(0)=0 , q(0)=0$ exist by the implicit function theorem:
\begin{eqnarray*}
p=c_1t+c_2t^2+\cdots,\quad q=d_1t+d_2t^2+\cdots.
\end{eqnarray*}
In the similar way, we have
\begin{eqnarray*}
c_1=\left(\frac{dp}{dt}\right)_0=\frac{1}{\Delta_1'}\begin{vmatrix}
(f_t)_0 &(-f_q)_0 \\ (g_t)_0 & 1-(g_q)_0\end{vmatrix},\\ 
d_1=\left(\frac{dq}{dt}\right)_0=\frac{1}{\Delta_1'}\begin{vmatrix}
1-(f_p)_0 &(f_t)_0 \\ (-g_p)_0 & (g_t)_0\end{vmatrix}, 
\end{eqnarray*}
where
\begin{eqnarray*}
\Delta_1=\Delta_1'=\begin{vmatrix}
1-(f_p)_0 &(-f_q)_0 \\ (-g_p)_0 & 1-(g_q)_0\end{vmatrix}\ne0,
\end{eqnarray*}
and
\begin{eqnarray*}
a_2=\frac{1}{2!}\left(\frac{d^2p}{dt^2}\right)_0=\frac{1}{2!}\frac{1}{\Delta_2'}\begin{vmatrix}
A_1' &(-f_q)_0 \\ B_1' & 1-(g_p)_0-(g_q)_0\end{vmatrix} ,\\
d_2=\frac{1}{2!}\left(\frac{d^2q}{dt^2}\right)_0=\frac{1}{2!}\frac{1}{\Delta_2'}\begin{vmatrix}
1-(f_p)_0-(f_q)_0   &A_1' \\ (-g_p)_0 & B_1'\end{vmatrix} ,
\end{eqnarray*}
where 
\begin{eqnarray*}
\Delta_2'&=&\begin{vmatrix}
1-(f_p)_0-(f_q)_0   &(-f_q)_0 \\ (-g_p)_0 & 1-(g_p)_0-(g_q)_0\end{vmatrix}\ne0,\\ 
A_1'&=&\biggl[f_{tt}+2f_{tp}\frac{dp}{dt}+f_{pp}(\frac{dp}{dt})^2+2f_{tq}\frac{dq}{dt}+2f_{pq}\frac{dp}{dt}
\frac{dq}{dt}
+ f_{qq}(\frac{dq}{dt})^2\biggr]_{t=0},\\
B_1'&=&\biggl[g_{tt}+2g_{tp}\frac{dp}{dt}+g_{pp}(\frac{dp}{dt})^2+2g_{tq}\frac{dq}{dt}+2g_{pq}\frac{dp}{dt}
\frac{dq}{dt}
+ g_{qq}(\frac{dq}{dt})^2\biggr]_{t=0},
\end{eqnarray*}
and so on.

Comparing $(u, v)$ with $(p, q)$, we have 
\begin{eqnarray*}
a_n\Delta_n=c_n\Delta_n',\quad b_n\Delta_n=d_n\Delta_n',\quad
\left|\Delta_n'\right|\leq\left|\Delta_n \right|.
\end{eqnarray*}
Therefore  
\begin{eqnarray*}
\left|a_n\right|\leq\left|c_n \right|,\quad  \left|b_n\right|\leq\left|d_n \right|.
\end{eqnarray*}
Since $p$ and $q$ are dominant series  of $u$ and $v$,  $u=\sum_{k=1}^\infty a_kt^k$ and 
$v=\sum_{k=1}^\infty b_kt^k$ are convergent.  Thus the theorem is proved. \hfill$\boxed{}$


\begin{thebibliography}{3}
\bibitem{AK97} F.~V.~Andreev, and A.~V.~Kitaev, 
On connection formulas for the asymptotics of some special solutions of the fifth Painlev\'e equation, 
{\it Zap. Nauchn. Sem. S.-Peterburg. Otdel. Mat. Inst. Steklov.} {\bf   243} (1997),   19--29;
 translation in {\it J. Math. Sci. (New York)} {\bf 99} (2000),   808--815. 

\bibitem{AK00} F.~V.~Andreev, and A.~V.~Kitaev, 
Connection formulae for asymptotics of the fifth Painlev\'e transcendent 
on the real axis, {\it Nonlinearity} {\bf 13} (2000),  1801--1840. 


\bibitem {BB} C.~Briot and J.-C.~Bouquet,  
Recherches sur les propri\'et\'es des fonctions d\'efinies par des \'equations 
diff\'erentielles,
{\it Journal de l'Ecole Polytechnique.} {\bf 36e} Cahier, Tome {\bf 21} (1856), 133--198.

\bibitem {BR} A.~D.~Bryuno and~I.~V.~Goryuchkina, Expansion of solutions of the sixth
Painlev\'e equation, (Russian) {\it Dokl. Akad. Nauk} {\bf 395} (2004),  733--737.

\bibitem {FZ} A.~S.~Fokas and X.~Zhou, On the solvability of Painlev\'e II and IV, 
{\it Comm.~Math.~Phys.} $\boldsymbol{144}$, (1992), 601--622.

\bibitem {RF05} R.~Fuchs,  
Sur quelques \'equations diff\'erentielles lin\'earires du second order,
{\it C. R. Acad. Sci. Paris} {\bf 141} (1905), 555--558.

\bibitem {RF1} R.~Fuchs, \"Uber lineare homogene Differentialgleichungen zweiter Ordnung  
mit drei im endlichen gelegenen wesentlich singul\"aren Stellen. 
{\it Math.~Ann.} $\boldsymbol{63}$ (1907), 301--321. 

\bibitem {RF2} R.~Fuchs, \"Uber lineare homogene Differentialgleichungen zweiter Ordnung 
     mit drei im endlichen gelegenen wesentlich singul\"aren Stellen. 
{\it Math.~Ann.} $\boldsymbol{70}$   (1911), 525--549.

\bibitem {RF3} R.~Fuchs, \"Uber die analytische Natur der L\"osungen von Differential
gleichungen zweiter ordnung mit festen kritischen Punkten, 
{\it Math.~Ann.} $\boldsymbol{75}$     (1914), 469--496.

\bibitem {Gambier} B.~Gambier, 
 Sur les \'equations  diff\'erentielles du second ordre
 et du premier degr\'e dont l'int\'egrale g\'en\'erale est \`a points
 critiques fix\'es, {\it Acta Math.} {\bf 33}  (1909), 1--55.
     
\bibitem {Garnier:1902} R.~Garnier,
Sur des \'equations diff\'erentielles du troisi\`eme
ordre dont l'int\'egrale g\'en\'erale est uniforme et sur une classe d'\'equations
nouvelles d'order sup\'erieur dont l'int\'egrale g\'en\'erale 
a ses points critiques fixes, 
{\it Ann.~Sci.~\'Ecole Norm.~Sup.} (3) {\bf 29}  (1912) , 1--126.

\bibitem {RG} R.~Garnier, 
\'Etude de l'int\'egrale g\'en\'erale de l'\'equation VI 
de M.~Painlev\'e dans le voisinage de ses singularit\'es transcendantes. 
{\it Ann. Sci. de l'\'Ecole Normale Sup. S\'er 3},~$\boldsymbol{34}$, (1917), 239--353.

\bibitem {GLS} V.~I.~Gromak, I.~Laine and  S.~Shimomura, 
{\it Painlev\'e Equations in the Complex Plane},
Walter de Gruyter, Berlin-New York 2002.

\bibitem {DG} D.~Guzzetti,
 Matching procedure for the sixth Painlev\'e equation, 
{\it the Preprint RIMS-1541}, May 2006.

\bibitem {DG1} D.~Guzzetti, On the critical behavior, the connection problem 
and the elliptic representation of a Painlev\'e VI Equation, 
{\it  Math. Phys. Anal. Geom.} {\bf 4} (2001),  293--377. 

\bibitem {DG2} D.~Guzzetti, 
The elliptic representation of the general Painlev\'e VI equation, 
{\it Comm. Pure Appl. Math.} {\bf 55} (2002), 1280--1363.

\bibitem {JH} J.~Heading, 
The Stokes phenomenon and Whittaker function, 
{\it J.~London Math. Soc.} $\boldsymbol{37}$ (1962), 195--208.

\bibitem {Hum}
J.~E.~Humphreys,  {\it Reflection Groups and Coxeter Groups}, Cambridge (1992). 
     
     

\bibitem {MJ} M.~Jimbo and T.~Miwa, 
Monodromy preserving deformation of linear ordinary 
     differential equations with rational coefficients. II, 
{\it Phys.~D}, $\boldsymbol{2}$ (1981), 407--448.
     
\bibitem {MJ2} 
M.~Jimbo, 
Monodromy problem and the boundary condition for some Painlev\'e equations,
{\it Publ. Res. Inst. Math. Sci.} {\bf 18} (1982), 1137--1161.

\bibitem {JMMS} 
M. Jimbo, T. Miwa, Y. Mori and M. Sato, 
Density matrix of an impenetrable Bose gas and the fifth Painlev\'e transcendent, 
{\it Phys.~D} $\boldsymbol{1}$ (1980), 80--158. 

\bibitem {KK1} K.~Kaneko, 
A new solution of the fourth Painlev\'e equation with a solvable monodromy, 
{\it Proc.~Japan Acad.~Ser. A Math. Sci.} $\boldsymbol {81}$ (2005), 75--79.
    
\bibitem {KOH} K.~Kaneko, Y.~Ohyama, 
Fifth Painlev\'e transcendents which are analytic at the origin, 
to appear.
     
\bibitem {KK2} K.~Kaneko, 
Painlev\'e VI transcendents which are meromorphic at a fixed singularity,
{\it Proc.~Japan Acad.~Ser. A Math. Sci.} $\boldsymbol{82}$ (2006), 71--76.

\bibitem {KOK} K.~Kaneko, S.~Okumura,
Special solutions of the sixth Painlev\'e equation with solvable monodromy, 
to appear.     

\bibitem {AVK} A.~V.~Kitaev, 
Symmetric solutions for the first and second Painlev\'e equations, 
{\it Zap. Nauchn. Sem. Leningrad. Otdel. Mat. Inst. Steklov. (LOMI)}
$\boldsymbol{187}$ (1991),  129--138,
translation in {\it J. Math. Sci.} $\boldsymbol{73}$ (1995), 494-499. 

\bibitem {AVK1} A.~V.~Kitaev, 
Grothendieck's dessins d'enfants, their deformations, 
and algebraic solutions of the sixth Painlev\'e and Gauss hypergeometric equations,  
{\it St. Petersburg Math. J.}   $\boldsymbol{17}$, (2006), 169--206. 
    
\bibitem {AVKK} A.~V.~Kitaev and Korotkin, 
On solutions of the Schlesinger equations in terms of $\Theta$-functions, 
{\it Internat. Math. Res. Notices} $\boldsymbol{17}$ (1998), 877--905.
     
\bibitem {MMA} M.~Mazzocco, 
Picard and Chazy solutions to the Painlev\'e VI equation, 
{\it Math. Ann.} $\boldsymbol{321}$ (2001), 157--195.
     
\bibitem {YM} Y.~Murata, 
Rational solutions of the second and the fourth Painlev\'e equations, 
{\it Funkcial. Ekvac.} {\bfseries 28} (1985), 1--32.
     
\bibitem {MN1} M.~Noumi and K.~Okamoto, 
Irreducibility of the second and the fourth Painlev\'e equations, 
{\it Funkcial. Ekvac.} {\bfseries 40} (1997), 139--163. 

\bibitem {NY} M.~Noumi and Y.~Yamada, 
A new Lax pair for the sixth Painlev\'e equation associated with $s\hat{}o(8)$, 
in ``{\it Microlocal Analysis and Complex Fourier Analysis}'' 
  (Eds. T.~Kawai and K.~Fujita), 238--252, World Scientific, 2002.

\bibitem {KO} K.~Okamoto, 
Isomonodromic deformation and Painlev\'e equations, and the Garnier system. 
{\it J. Fac. Sci. Univ. Tokyo Sect. 1A, Math.} $\boldsymbol{33}$ (1986), 575--618.
     

\bibitem {OO} Y.~Ohyama, S.~Okumura,  
 R. Fuchs' problem of the Painlev\'e equations from the first to the fifth, 
{\texttt {math.CA/0512243}}.


\bibitem{P:06} P. Painlev\'e,  
Sur les \'equations diff\'erentielles du second ordre \`a points critiques fixes,
{\it C. R. Acad. Sci. Paris} {\bf 143} (1906), 1111--1117.

\bibitem {EP} E.~Picard, 
M\'emoire sur la th\'eorie des functions alg\'ebriques de deux variables,
{\it J. de Liouville} $\boldsymbol{5}$ (1889), 135--319.  
     
\bibitem {Sibuya} Y.~Sibuya,  
Convergence of power series solutions of a linear Pfaffian system at an irregular singularities,
{\it Keio Engrg.~Rep.} {\bf 31} (1978),  79--86.  

\bibitem {SS1}  S.~Shimomura, 
Painlev\'e Transcendents in the neighbourhood of fixed sjngular points,
{\it Funkcial. Ekvac.}  $\boldsymbol{25}$ (1982), 163--184.

\bibitem {SS2}  S.~Shimomura, 
Series expansions of Painlev\'e transcendents in the neighbourhood of a fixed singular Point,
{\it Funkcial. Ekvac.}    $\boldsymbol{25}$ (1982), 185--197.

\bibitem {KT} K.~Takano, 
Reduction for Painlev\'e equations at the fixed singular points of the first kind, 
{\it Funkcial. Ekvac.}  $\boldsymbol{29}$ (1986), 99--119.

\bibitem {HU1} H.~Umemura, 
Birational automorphism groups and differential equations.
{\it Equations differentielles dans le champ complex,} Vol.~$\boldsymbol{II}$ 
(Strasbourg,1985), 119--227, Publ. Inst. Rech. Math. Av., Univ. Louis Pasteur, Strasbourg.
(also in {\it Nagoya Math. J.} {\bfseries 119} (1990), 1--80)

\bibitem{HU2} H.~Umemura and H.~Watanabe, 
{Solutions of the second and the fourth Painlev\'e equations}, 
{\it Nagoya Math. J.} {\bfseries 148} (1997), 151--198.


     
\bibitem {ETW} E.~T.~Whittaker and G.~N.~Watson, 
{\it A course of modern analysis}, 
Cambridge at the university press, (1927).


\end{thebibliography}
\end{document}